\documentclass[11pt]{article}
\usepackage{array,amssymb,amsmath,verbatim,amsthm,color}
\usepackage{graphicx}
\usepackage{booktabs}
\usepackage{footnote}
\usepackage[numbers]{natbib} \setcitestyle{open={},close={}}

\newtheorem{theorem}{Theorem}[section]
\newtheorem{lemma}[theorem]{Lemma}

\newtheorem{corollary}[theorem]{Corollary}

\newtheorem{example}[theorem]{Example}
\usepackage{longtable}
\usepackage[figuresright]{rotating}
\usepackage{float}
\newtheorem{construction}[theorem]{Construction}

\newcommand{\ignore}[1]{}

\usepackage{caption}

\begin{document}
\title{Maximum $w$-cyclic holely group divisible packings with block size three and applications to optical orthogonal codes\footnote{Supported by NSFC under grant 11971053 (J. Zhou), NSFC under grant 11771119 and NSFHB under A2019507002 (L. Wang).}
}
\author{{Zenghui Fang}$^{1}$,  Junling Zhou$^{1}$, Lidong Wang$^{2}$\\
\small{${}^{1}$Department of Mathematics,
Beijing Jiaotong University,
Beijing 100044, P. R. China}\\
\small{${}^{2}$Department of Smart Policing,
China People's Police University,
Langfang 065000, P. R. China}\\
\small{17118435@bjtu.edu.cn,
jlzhou@bjtu.edu.cn, lidongwang@bjtu.edu.cn}
}
\date{ }
\maketitle
\begin{abstract}
In this paper we investigate combinatorial constructions for $w$-cyclic holely group divisible packings with block size three (briefly by $3$-HGDPs). For any positive integers  $u,v,w$ with $u\equiv0,1~(\bmod~3)$, the exact number of base blocks of a maximum  $w$-cyclic $3$-HGDP of type $(u,w^v)$ is determined.
This result is used to determine the exact number of codewords in a maximum three-dimensional $(u\times v\times w,3,1)$ optical orthogonal code with 
at most one optical pulse per spatial plane and per wavelength plane.


\end{abstract}
\section{Introduction and preliminaries}
Let $K$ be a set of positive integers. A {\it group divisible packing}, denoted by $K$-GDP, is a triple $(X,\mathcal{G},\mathcal{B})$, where $X$ is a finite set of points, $\mathcal{G}$ is a partition of $X$ into subsets (called {\it groups}) and $\mathcal{B}$ is a collection  of subsets (called {\it blocks}) of $X$, each block has cardinality from $K$, such that any pair of $X$  from two different groups is contained in at most one block.
The multi-set $T=\{|G|:G\in\mathcal{G}\}$ is called the {\it type}  of the $K$-GDP.
If $\mathcal{G}$ contains $u_i$ groups of size $g_i$ for $1\leq i\leq r$, we also denote the  type  by $g_1^{u_1}g_2^{u_2}\cdots g_r^{u_r}$. If $K=\{k\}$, we write a $\{k\}$-GDP as  $k$-GDP. If every  pair from different groups appears in exactly one block, then the $K$-GDP is referred to as a {\it group divisible design} and denoted by $K$-GDD.
A $K$-GDD of type $1^v$ is commonly called a {\it pairwise balanced design}, denoted by PBD$(v,K)$.

A sub-GDD $(Y,\mathcal{H},\mathcal{A})$ of a GDD $(X,\mathcal{G},\mathcal{B})$ is a GDD satisfying that $Y\subseteq X, \mathcal{A}\subseteq\mathcal{B}$, and every group of $\mathcal{H}$ is contained in some group of $\mathcal{G}$. If a GDD has a missing sub-GDD, then we say that the GDD has a {\it hole}. (In fact, the sub-GDD as a hole may not exist.) If a GDD has several holes which partition the point set of the GDD, then we call it a double GDD. Next we generalize the concept of double GDD to a formal definition of double group divisible packing.

Let $u,v$ be positive integers and $K$ be a set of positive integers. A {\it double group divisible packing}, denoted by $K$-DGDP, is a quadruple $(X,\mathcal{G},\mathcal{H},\mathcal{B})$ satisfying that
\begin{description}
\item{(1)} $X$ is a finite set of points;

 \item{(2)} $\mathcal{G}=\{G_1,G_2,\ldots, G_u\}$ is a partition of $X$ into $u$ subsets (groups);

 \item{(3)} $\mathcal{H}=\{H_1,H_2,\ldots, H_v\}$ is another partition of $X$ into $v$ subsets (called {\it holes});

 \item{(4)} $\mathcal{B}$ is a collection of subsets (blocks)  of $X$ with cardinality from $K$, such that
each block meets each group and each hole in at most one point, and any two points from different groups and different holes appear together in at most one block.
 \end{description}
 \noindent
Let $G_i\cap H_j=R_{ij}$ for $1\leq i\leq u$ and $1\leq j\leq v$. The $u\times v$ matrix $T=(|R_{ij}|)$ is called the {\it type} of this design.
If $K=\{k\}$, we write a $\{k\}$-DGDP as $k$-DGDP. If every two points from different groups and different holes appear in exactly one block, then the $K$-DGDP
is a {\it double group divisible design} and denoted by $K$-DGDD.

For $0\leq e\leq v-1$, let $\mathbf{c}_e=(w,w,\ldots,w)^T$ be a column vector of length $u$. A $K$-DGDP  of type $(\mathbf{c}_0,\mathbf{c}_1,\ldots,\mathbf{c}_{v-1})$ is usually called a {\it holely group divisible packing} and denoted by $K$-HGDP of type $(u,w^v)$.  Note that a $K$-HGDP of type $(u,w^v)$ is equivalent to a $K$-HGDP of type $(v,w^u)$ because we can exchange the expression of groups and holes. For $0\leq e\leq v-1$, let $\mathbf{c}_e=(w,w,\ldots,w,wt)^T$ be a column vector of length $u-t+1~(0\leq t<u)$. We call a $K$-DGDP of type $(\mathbf{c}_0,\mathbf{c}_1,\ldots,\mathbf{c}_{v-1})$  an {\it incomplete holely group divisible packing} and denote by $K$-IHGDP of type $(u,t,w^v)$, by which we mean that  the group of  size $tvw$ is regarded as the union of $t$ disjoint parts of equal size $vw$ such that each pair of points from the $t$ parts does not appear in any block. Especially, if $t=0,1$, a $K$-IHGDP of type $(u,t,w^v)$  is the same as a $K$-HGDP of type $(u,w^v)$.
When $w=1$,
a $K$-HGDP of type $(u,1^v)$ is often said to be a {\it modified group
divisible packing} and denoted by $K$-MGDP of type $v^u$.

Naturally we have the notions of a holely group divisible design (HGDD), an incomplete holely group divisible design (IHGDD) and a modified group divisible design (MGDD).
Wei~[\cite{wei93}] first introduced the concept of HGDDs and gave a complete solution to the existence of $3$-HGDDs.
\begin{theorem}\label{hgdd}{\rm [\cite{wei93}]}
There exists a $3$-HGDD of type $(u,w^v)$ if and only if $u,v\geq3,$ $(u-1)(v-1)w\equiv 0~(\bmod~2)$ and $uv(u-1)(v-1)w^2\equiv 0~(\bmod~6).$
\end{theorem}

 In what follows we always denote $I_{u}=\{0,1,\ldots,u-1\}$ to be a set and denote by $Z_{v}$ the additive group of integers in $I_v$ (implying that the arithmetic is taken modulo $v$). Suppose that $(X,\mathcal{G},\mathcal{H},\mathcal{B})$ is a $k$-DGDP of type $(\mathbf{d}_0,\mathbf{d}_1,\ldots,\mathbf{d}_{v-1})$, where $\mathbf{d}_e=(wg_1,wg_2,\ldots,$ $wg_u)^T$ for $0\leq e\leq v-1$. Let $\pi$ be a permutation on $X$. If $\pi$ keeps $\mathcal{G},$ $\mathcal{H}$ and $\mathcal{B}$ invariant, then $\pi$ is an {\it automorphism} of the $k$-DGDP.   The block set $\mathcal{B}$ is partitioned into equivalence classes called {\it block orbits} under the action of any automorphism $\pi$. A set of {\it base blocks} is an arbitrary set of representatives for these block orbits of $\mathcal{B}$. In particular, if the  $k$-DGDP admits an automorphism which is the product of $v\sum_{i=1}^{u}g_i$ disjoint $w$-cycles, fixes every group and every hole, and leaves $\mathcal{B}$ invariant, then the design is said to be {\it $w$-cyclic}.  
Without loss of generality, we always construct a  $w$-cyclic $k$-DGDP  on the point set $X=(\bigcup_{i=1}^{u}R_i)\times I_v\times Z_w$ with group set $\mathcal{G}=\{R_i\times I_v\times Z_w:1\leq i\leq u\}$ and hole set $\mathcal{H}=\{(\bigcup_{i=1}^{u}R_i)\times \{j\}\times Z_w: 0\leq j\leq v-1\}$, where $R_i$, $1\leq i\leq u$, are $u$ pairwise disjoint sets with $|R_i|=g_i$. In this case, the automorphism will be taken as $(a,b,c)\mapsto (a,b,c+1)\bmod (-,-,w)$ for $a\in \bigcup_{i=1}^{u}R_i,~b\in I_v$ and $c\in Z_{w}$.

For positive integers $u,v,w,k$ with $u,v\geq k$, let $\Psi(u\times v\times w,k,1)$ denote the largest possible number of base blocks of a $w$-cyclic $k$-HGDP of type $(u,w^v)$.
 A $w$-cyclic $k$-HGDP of type $(u,w^v)$ is called {\it maximum} if it contains $\Psi(u\times v\times w,k,1)$ base blocks.

 Based on the Johnson bound~[\cite{s}] for constant weight codes, Dai $et~al.$~[\cite{DCW}] displayed
an upper bound of $\Psi(u\times v\times w,k,1)$, as stated in the following lemma.
\begin{lemma}{\rm[\cite{DCW}, Theorem 5.2]}
\label{Johnson bound}
$\Psi(u\times v\times w,k,1)\leq J(u\times v\times w,k,1)$ holds for any positive integers $u,v \geq k,$ where $$J(u\times v\times w,k,1)=\lfloor\frac{uv}{k}\lfloor\frac{(u-1)(v-1)w}{k-1}\rfloor\rfloor.$$
\end{lemma}

 This paper is devoted to research on maximum $w$-cyclic $3$-HGDPs.  The main theorem  is as follows.
\begin{theorem}\label{all}
For all positive integers $u,v,w$ with $u\equiv 0,1~(\bmod~3)~(u,v\geq3)$, we have
$$\Psi(u\times v\times w,3,1)=\left\{
  \begin{array}{ll}
   6(w-1), & \hbox{if $(u,v)=(3,3)$ and $w\equiv 0~(\bmod~2)$,} \\
   J(u\times v\times w,3,1), & \hbox{otherwise.}
  \end{array}
\right.$$
\end{theorem}

 In the following of this section we prepare some preliminary results of close relevance to $k$-HGDPs.

 If $\Psi(u\times v\times w,k, 1)=\frac{uvw(u-1)(v-1)}{k(k-1)},$ then we can check that a $w$-cyclic $k$-HGDP of type $(u,w^v)$ with $\Psi(u\times v\times w,k, 1)$ base blocks is equivalent to  a $w$-cyclic $k$-HGDD of type $(u,w^v)$. For $k=3,$
Dai $et~al.$~[\cite{DCW}] solved the existence problem of $w$-cyclic 3-HGDDs of type $(u,w^v)$.
\begin{theorem}{\rm{[\cite{DCW},~Theorem 4.11]}}\label{Yt11}
There exists a $w$-cyclic $3$-HGDD of type $(u,w^v)$  if and only if $u,~v\geq3,$ $(u-1)(v-1)w\equiv0~(\bmod~2)$ and $uv(u-1)(v-1)w\equiv0~(\bmod~6)$ with the exceptions of $u=v=3$ and $w\equiv0~(\bmod~2)$.
\end{theorem}
Combining Lemma~\ref{Johnson bound} and Theorem~\ref{Yt11} yields the following corollary.
\begin{corollary}\label{co1}
$ \Psi(u\times v\times w,3, 1)=J(u\times v\times w,3,1)$ if any of the following conditions is satisfied: $(1)~u\equiv 1,3~(\bmod~6)$, $u,v\geq3$ with the exceptions of $u=v=3$ and $w\equiv0~(\bmod~2);$ $(2)~u\equiv0,4~(\bmod~6)$, $v\equiv1~(\bmod~2)$ and $w\equiv 1~(\bmod~2);$
$(3)~u\equiv0,4~(\bmod~6)$,~$v\geq4$ and $w\equiv0~(\bmod~2).$
\end{corollary}

Suppose that $(X,\mathcal{G},\mathcal{H},\mathcal{B})$ is a $K$-HGDD of type $(u,w^v)$.
It is said to be {\it semi-cyclic}, briefly by $K$-SCHGDD of type $(u,w^v)$, if it admits an automorphism which is the product of $u$ disjoint $vw$-cycles, fixes every group, leaves $\mathcal{H}$ and $\mathcal{B}$ invariant.
Without loss of generality, for a $K$-SCHGDD of type $(u,w^v)$, we always identify $X$ with $I_{u}\times Z_{vw}$, $\mathcal{G}$ with $\{\{i\}\times Z_{vw}:i\in I_u\}$  and $\mathcal{H}$ with $\{I_u\times\{i,v+i,\ldots,(w-1)v+i\}:0\leq i\leq v-1\}.$ In this case, the automorphism will be taken as $(i,x)\mapsto (i,x+1)\bmod (-,vw)$ for $ i\in I_{u}$ and $x\in Z_{vw}$. The existence of a $k$-SCHGDD of type $(u,w^v)$ implies that of a $w$-cyclic $k$-HGDD of type $(u,w^v)$, as indicated by the following result.
\begin{lemma}\rm{[\cite{WC15}, Construction~4.1]}\label{c15}
If there exists a $k$-SCHGDD of type $(u,(gw)^v)$, then there exists a $w$-cyclic $k$-HGDD of type $(u,(gw)^v)$.
\end{lemma}
 For $k=3$, we record the existence results of  $3$-SCHGDDs which were studied by [\cite{fwc,fww,WWF}].
\begin{theorem}\label{1}
{\rm [\cite{WWF}, Theorem 1.4]}
There exists a $3$-SCHGDD of type $(u,w^v)$ if and only if $u,v\geq3$, $(u-1)(v-1)w\equiv0~(\bmod~2)$ and $u(u-1)(v-1)w\equiv 0~(\bmod~6)$ except when
 $(1)~u\equiv3,7~(\bmod~12), w\equiv1~(\bmod~2)$ and $v\equiv2~(\bmod~4);$
  $(2)~u=3$, $w\equiv1~(\bmod~2)$ and $v\equiv 0~(\bmod~2);$
  $(3)~u=v=3$ and $w\equiv 0~(\bmod~2);$
$(4)~(u,w,v)\in\{(5,1,4),(6,1,3)\};$
 $  (5)~u\equiv11~(\bmod~12)$,~either $w\equiv3~(\bmod~6)$ and $v\equiv2~(\bmod~4)$,~or $w\equiv1,5~(\bmod~6)$ and $v\equiv 10~(\bmod~12)$.
\end{theorem}

In this paper, we concentrate on  combinatorial constructions for maximum $w$-cyclic $3$-HGDPs of type $(u,w^v)$. Our main purpose is to determine the value of $\Psi(u\times v\times w,3,1)$ for all positive integers $u,v,w$ where $u\equiv 0,1~(\bmod~3)$ $(u,v\geq 3)$.
The rest of the paper is organized as follows. In Section 2, we present a few recursive constructions for $w$-cyclic $k$-DGDPs and show some special applications to produce maximum $w$-cyclic $3$-HGDPs.  In Section 3, we shall determine the number of base blocks of maximum $w$-cyclic 3-HGDPs of type $(u,w^v)$ when $u\equiv 0~(\bmod~3)$. In Section 4, we construct two types of auxiliary designs, namely incomplete HGDPs and semi-cyclic $K$-GDDs. In section $5$ we pay attention to the case  $u\equiv1~(\bmod~3)$ and finally prove the main theorem. The existence of maximum $w$-cyclic 3-HGDPs of type $(u,w^v)$ not only is interesting in its own right, but also has nice applications in three-dimensional optical orthogonal codes with certain restrictions. We consider this application in Section 6. In the Appendix we provide Supporting Information on many direct constructions for small designs utilized in the proofs of the paper.\footnote{The Appendix for Supporting Information can be found in the full version of the paper: arXiv}

\section{ Recursive constructions}
In this section, some recursive constructions for $w$-cyclic $k$-DGDPs  will be given. Then several special constructions yielding maximum $w$-cyclic $3$-HGDPs will be displayed.
\begin{construction}
\label{xc1}
Suppose that the following designs exist:

 $(1)$ a $w$-cyclic $k$-DGDP of type $(\mathbf{d}_0,\mathbf{d}_1,\ldots,\mathbf{d}_{v-1})$ with $b$ base blocks where $\mathbf{d}_e=(wm_1,wm_2,\\\ldots,wm_s)^T$ for $0\leq e\leq v-1;$

 $(2)$ a $w$-cyclic $k$-IHGDP of type $(m_i+t,t,w^v)$
 with $c_i$ base blocks for each $ 1\leq i\leq s-1$.

\noindent
Then there exists a $w$-cyclic $k$-IHGDP of  type $(\sum_{i=1}^{s}m_i+t,m_s+t,w^v)$ with $b+\sum_{i=1}^{s-1}c_i$ base blocks.
Furthermore, if there exists a $w$-cyclic $k$-IHGDP of type $(m_s+t,t,w^v)$ with $c_s$ base blocks, then  there exists a
  $w$-cyclic $k$-IHGDP of  type $(\sum_{i=1}^{s}m_i+t,t,w^v)$ with $b+\sum_{i=1}^sc_i$ base blocks.
\end{construction}
\proof
Let $R_i$, $1\leq i\leq s$, be $s$ pairwise disjoint sets with $|R_i|=m_i$.
By assumption we have a $w$-cyclic $k$-DGDP $(X,\mathcal{G},\mathcal{H},\mathcal{B})$  of type $(\mathbf{d}_0,\mathbf{d}_1,\ldots,\mathbf{d}_{v-1})$ with $b$ base blocks, where $X=(\bigcup_{i=1}^{s}R_{i})\times I_{v}\times Z_{w}$, $\mathcal{G}=\{R_{i}\times I_{v}\times Z_{w}:1\leq i\leq s\}$ and $\mathcal{H}=\{(\bigcup_{i=1}^{s}R_{i})\times \{j\}\times Z_{w}:j\in I_{v}\}.$
Let $\mathcal{F}$ be the collection  of all base blocks of this DGDP.

Let $E$ be a set of size $t$ such that $E\times I_v\times Z_w$ is disjoint from $X.$ For each group $R_i\times I_v\times Z_w, 1\leq i\leq s-1$, construct a $w$-cyclic $k$-IHGDP of type $(m_i+t,t,w^v)$ with $c_i$ base blocks 
on $(R_i\cup E)\times I_v\times Z_w$ with
group set $\{\{x\}\times I_v\times Z_w:x\in R_i\}\cup \{E\times I_v\times Z_w\}$ and hole set $\{(R_i\cup E)\times\{j\}\times Z_w:j\in I_v\}$.
Denote by $\mathcal{A}_i$ the collection of all base blocks of this design.
Let $\mathcal{A}=(\bigcup_{i=1}^{s-1}\mathcal{A}_i)\cup \mathcal{F}$. It is readily checked that $\mathcal{A}$ forms   $b+\sum_{i=1}^{s-1}c_i$  base blocks of a $w$-cyclic $k$-IHGDP of type $(\sum_{i=1}^{s}m_i+t,m_s+t,w^v)$
on $((\bigcup_{i=1}^{s}R_i)\cup E)\times I_v\times Z_w$ with
 group set $\{\{x\}\times I_v\times Z_w:x\in \bigcup_{i=1}^{s-1}R_i\}\cup\{(R_s\cup E)\times I_v\times Z_w\}$ and hole set $\{((\bigcup_{i=1}^{s}R_i)\cup E)\times\{j\}\times Z_w:j\in I_v\}$.

Furthermore, construct a $w$-cyclic $k$-IHGDP of type $(m_s+t,t,w^v)$  with $c_s$ base blocks 
 on $(R_s\cup E)\times I_v\times Z_w$ with group set $\{\{x\}\times I_v\times Z_w:x\in R_s\}\cup\{E\times I_v\times Z_w\}$ and hole set $\{(R_s\cup E)\times\{j\}\times Z_w:j\in I_v\}$.
Denote $\mathcal{A}_s$ the collection of all base blocks of this design. Let $\mathcal{C}=(\bigcup_{i=1}^{s}\mathcal{A}_i)\cup \mathcal{F}$. It is readily checked that $\mathcal{C}$ forms $b+\sum_{i=1}^sc_i$ base blocks of a $w$-cyclic
$k$-IHGDP of type $(\sum_{i=1}^{s}m_i+t,t,w^v)$.
\qed

Sometimes we do not consider if a design admits an automorphism group. So we list a corollary to Construction~\ref{xc1} as follows.
\begin{corollary}\label{col2}
Suppose that the following designs exist:

 $(1)$ a $k$-DGDP of type $(\mathbf{d}_0,\mathbf{d}_1,\ldots,\mathbf{d}_{v-1})$ with $b$ blocks where $\mathbf{d}_e=(wm_1,wm_2,\ldots,wm_s)^T$ for $0\leq e\leq v-1;$

 $(2)$ a $k$-IHGDP of type $(m_i+t,t,w^v)$ with $c_i$ blocks for each $ 1\leq i\leq s-1$.

\noindent
 Then there exists a $k$-IHGDP of  type $(\sum_{i=1}^{s}m_i+t,m_s+t,w^v)$ with $b+\sum_{i=1}^{s-1}c_i$ blocks. Furthermore, if there exists a $k$-IHGDP of type $(m_s+t,t,w^v)$ 
 with $c_s$ blocks, then  there exists a
 $k$-IHGDP of  type $(\sum_{i=1}^{s}m_i+t,t,w^v)$
  with $b+\sum_{i=1}^sc_i$ blocks.
\end{corollary}

Suppose that $(X,\mathcal{G},\mathcal{B})$ is a
$K$-GDP of type $\{wm_1,wm_2,\ldots,wm_r\}$. 
If there is a permutation $\pi$ on $X$ that is the product of $\sum_{i=1}^{r}m_i$ disjoint $w$-cycles, fixes every group and leaves $\mathcal{B}$ invariant, then the $K$-GDP is said to be {\it $w$-cyclic}. A $w$-cyclic $K$-GDP of type $w^u$ is called  {\it semi-cyclic} and denoted by $K$-SCGDP of type $w^u.$
Without loss of generality, we may construct a $K$-SCGDP on the point set $X=(\bigcup_{i=1}^{r}R_{i})\times Z_w$ with group set $\mathcal{G}=\{R_i\times Z_w:1\leq i\leq r\}$, where $|R_i|=m_i$. And the automorphism will be taken as $(a,b)\mapsto (a,b+1)\bmod (-,w)$ for $a\in \bigcup_{i=1}^{r}R_i$ and $b\in Z_{w}$.
  In some  cases a semi-cyclic $K$-GDD can be obtained from a semi-cyclic $L$-HGDD; we show this by an obvious lemma.

\begin{lemma}\label{eq}
Let $K$ be a set of positive integers and $u\not\in K$.
A $K$-SCHGDD of type $(u,1^v)$ can be viewed as a $K\cup\{u\}$-SCGDD of type $v^u$ with one base block of size $u$.
\end{lemma}

 \begin{construction}\label{xc2}
Suppose that there exists a $w$-cyclic $K$-GDD of type $\{wm_1,wm_2,\ldots,wm_r\}$ with $b_k$ base blocks for each $k\in K$. If there exists an $h$-cyclic $l$-HGDP of type $(k,h^v)$ with $c_k$ base blocks for each $k\in K,$ then there exists an $hw$-cyclic $l$-DGDP of type
$(\mathbf{d}_0,\mathbf{d}_1,\ldots,\mathbf{d}_{v-1})$ with $\sum_{k\in K} b_kc_k$ base blocks, where $\mathbf{d}_e=(hwm_1,hwm_2,\ldots,hwm_r)^T$ for $0\leq e\leq v-1$.
\end{construction}
\proof
 Start from the given $w$-cyclic $K$-GDD of type $\{wm_1,wm_2,\ldots,wm_r\}$ defined on $X=(\bigcup_{i=1}^{r}R_{i})\times Z_w$  with group set $\{R_i\times Z_w:1\leq i\leq r\}$, where $|R_i|=m_i$. Let $\mathcal{F}$ be the collection of the base blocks.

 For each block $B\in\mathcal{F}$ and $|B|=k$, construct an $h$-cyclic $l$-HGDP of type $(k,h^v)$ with $c_k$ base blocks on $Y=B\times I_v\times Z_h$ with group set $\{\{x\}\times I_v\times Z_h:x\in B\}$ and hole set $\{B\times \{j\}\times Z_h:j\in I_v\}$. Denote the family of base blocks by $\mathcal{A}_B$. Let $\mathcal{A}=\bigcup_{B\in\mathcal{F}}\mathcal{A}_B$.

For 
each $(a,b,c,d)\in Y$, define a mapping $\tau:(a,b,c,d)\mapsto (a,c,b+wd)$.
Then define $\tau(A)=\{\tau(x):x\in A\}$ and  $\mathcal{C}=\bigcup_{A\in\mathcal{A}}\tau(A)$. It is readily checked that  $\mathcal{C}$ forms all base blocks of an $hw$-cyclic $l$-DGDP of type $(\mathbf{d}_0,\mathbf{d}_1,\ldots,\mathbf{d}_{v-1})$
on $(\bigcup_{i=1}^{r}R_{i})\times I_v\times Z_{hw}$ with
group set $\{R_i\times I_v\times Z_{hw}:1\leq i\leq r\}$
 and hole set $\{(\bigcup_{i=1}^{r}R_{i})\times \{j\}\times Z_{hw}: j\in I_v\}$. It is immediate that the number of base blocks equals $\sum_{k\in K} b_kc_k$.
\qed

\begin{corollary}\label{col4}
Suppose that there exists a $K$-SCGDD of type $w^v$ with $b_k$ base blocks for $k\in K$. If there exists an $h$-cyclic $l$-HGDP of type $(u,h^k)$ with $c_k$ base blocks for each $k\in K,$ then there exists an $hw$-cyclic $l$-HGDP of type $(u,(hw)^v)$
with $\sum_{k\in K} b_kc_k$ base blocks. 
\end{corollary}
\proof
Start from the given $K$-SCGDD, also a $w$-cyclic $K$-GDD, of type $w^v$.
Apply Construction~\ref{xc2} with  an $h$-cyclic $l$-HGDP of type $(u,h^k)$, also of type $(k,h^u)$, with $c_k$ base blocks, for each $k\in K$. Then we obtain an $hw$-cyclic $l$-HGDP of type $(v,(hw)^u)$, also of type $(u, (hw)^v)$, with $\sum_{k\in K} b_kc_k$ base blocks.
\qed
\begin{construction}\label{xc3}
Suppose that the following designs exist:

   $(1)$ a $K$-HGDD of type $(u,1^v)$ which has $b_k$  blocks of size $k$ for each $k\in K;$

   $(2)$ an $l$-SCGDP of type $w^{k}$ with $e_{k}$ base blocks of size $k$ for each $k\in K$.

\noindent
Then there exists a $w$-cyclic $l$-HGDP of  type $(u,w^{v})$ with $\sum_{k\in K}{b_{k}e_{k}}$ base blocks.
\end{construction}
\proof
 Start from the given $K$-HGDD of type $(u,1^{v})$ defined on $X=I_{u}\times I_{v}$ with group set $\mathcal{G}=\{\{i\}\times I_{v}:i\in I_{u}\}$ and hole set $\mathcal{H}=\{I_{u}\times\{j\}:j\in I_{v}\}.$ Denote the family of blocks by $\mathcal{F}$.
 For each $B\in\mathcal{F}$ and $|B|=k$, construct on $B\times Z_{w}$ an $l$-SCGDP of type $w^{k}$ with $e_k$ base blocks and with group set $\{\{x\}\times Z_{w}:x\in B\}$. Let $\mathcal{A}_{B}$ be the family of base blocks of this design. Define $\mathcal{A}=\bigcup_{B\in\mathcal{F}}{\mathcal{A}_{B}}$. It is readily checked that $\mathcal{A}$ forms  $\sum_{k\in K}b_ke_k$ base blocks of a $w$-cyclic $l$-HGDP of type $(u, w^{v})$ on $I_{u}\times I_{v}\times Z_{w}$ with group set $\mathcal{G}=\{\{i\}\times I_{v}\times Z_{w}:i\in I_u\}$ and hole set $\mathcal{H}=\{I_{u}\times\{j\}\times Z_{w}:j\in I_{v}\}$.
\qed

\begin{construction}\label{xc4}
Let  $\mathbf{c}_e=(g_1,g_2,\ldots,g_n)^T$  for $0\leq e\leq vw-1$ and $\mathbf{d}_f=(wg_1,wg_2,\ldots,wg_n)^T$ for $0\leq f\leq v-1$.
If there exist a $k$-DGDD of type $(\mathbf{d}_0,\mathbf{d}_1,\ldots,\mathbf{d}_{v-1})$ with b blocks and a $k$-DGDP of
type $(\mathbf{c}_0,\mathbf{c}_1,\ldots,\mathbf{c}_{w-1})$ with c blocks,
then there exists a  $k$-DGDP of type $(\mathbf{c}_0,\mathbf{c}_1,\ldots,\mathbf{c}_{vw-1})$ with $b+cv$ blocks.
\end{construction}
\proof
Let $R_i, 1\leq i\leq n$, be $n$ pairwise disjoint sets with $|R_i|=g_i$.~By assumption, we may construct a $k$-DGDD $(X,\mathcal{G},\mathcal{H},\mathcal{B})$
of type $(\mathbf{d}_0,\mathbf{d}_1,\ldots,\mathbf{d}_{v-1})$ with $b$ blocks,  where $X=(\bigcup_{i=1}^{n}R_i)\times I_{v}\times I_{w}$, $\mathcal{G}=\{R_i\times I_{v}\times I_w:1\leq i\leq n\}$,
and $\mathcal{H}=\{(\bigcup_{i=1}^{n}R_i)\times \{j\}\times I_w: j\in I_{v}\}$.

For each hole $H=(\bigcup_{i=1}^{n}R_i)\times \{j\}\times I_w,~0\leq j\leq v-1$, we construct a $k$-DGDP of type $(\mathbf{c}_0,\mathbf{c}_1,\ldots,\mathbf{c}_{w-1})$ with $c $ blocks on $H$ 
with group set $\{R_i\times\{j\}\times I_w: 1\leq i\leq n\} $ and hole set \{$(\bigcup_{i=1}^{n}R_i)\times\{j\}\times \{l\}: l\in I_w$\}. Denote the set of all blocks by $\mathcal{A}_{H}$. Let $\mathcal{A}=\bigcup_{H\in\mathcal{H}}\mathcal{A}_{H}.$
 It is easily checked that $\mathcal{A}\cup\mathcal{B}$ forms $b+cv$  blocks of the desired design on $(\bigcup_{i=1}^{n} R_i)\times I_{vw}$ with group set $\{R_i\times I_{vw}:1\leq i\leq n\}$ and hole set $\{(\bigcup_{i=1}^{n} R_i)\times\{j\}:j\in I_{vw}\}$.\qed

If we apply the previous recursive constructions and choose appropriate ingredient designs, then we  produce maximum $w$-cyclic $3$-HGDPs. To end  this section we show three special cases of constructions to produce $3$-HGDPs of type $(u,w^v)$ with $J(u\times v\times w,3,1)$ base blocks.

\begin{lemma}\label{zz}
Let  $u,v,g$ be even, $w$ be odd, $gu\equiv0~(\bmod~3)$ and $v\equiv s~(\bmod~g)$.
Suppose that
a $w$-cyclic $3$-DGDD of type $(\mathbf{d}_0,\mathbf{d}_1,\ldots,\mathbf{d}_{u-1})$ exists,
 where $\mathbf{d}_e=(gw,gw,\ldots,gw,sw)^T$ is
a column vector of length $\frac{v-s+g}{g}$ for $0\leq e\leq u-1$. If there exists a $w$-cyclic $3$-HGDP of type $(k,w^u)$ with $J(k\times u\times w,3,1)$ base blocks for each $k\in\{g,s\}$, then there  exists a $w$-cyclic $3$-HGDP of type $(v,w^u)$ with $J(v\times u\times w,3,1)$ base blocks.
\end{lemma}
\proof
 By assumption, there exists a $w$-cyclic $3$-DGDD of type $(\mathbf{d}_0,\mathbf{d}_1,\ldots,\mathbf{d}_{u-1})$, which has
$\frac{uw(u-1)(v-s)(v+s-g)}{6}$ base blocks. Apply Construction~\ref{xc1} with a $w$-cyclic $3$-HGDP of type $(k,w^u)$ with $J(k\times u\times w,3,1)$ base blocks for $k\in\{g,s\}$ to produce a $w$-cyclic $3$-HGDP of type $(v,w^u)$.
It has $$ \frac{uw(u-1)(v-s)(v+s-g)}{6}+\frac{v-s}{g}\times \frac{gu(g-1)(u-1)w-gu}{6}+\lfloor
 \frac{su(s-1)(u-1)w-su}{6}\rfloor$$
 $$~~~~~=\lfloor\frac{vu((v-1)(u-1)w-1)}{6}\rfloor=J(v\times u\times w,3,1)$$
base blocks.\qed

\begin{lemma}\label{z7}
Let $u,v$ be even, $mg$ be odd and $u(u-1)m \equiv0~(\bmod~3).$
If there exist

 $(1)$ a $\{3,5,v\}$-SCGDD of type $g^v$ containing one base block of size $v;$  and

 $(2)$ an $m$-cyclic $3$-HGDP of type $(u,m^v)$ with $J(u\times v\times m,3,1)$ base blocks.

\noindent
 Then there exists an $mg$-cyclic $3$-HGDP of type $(u,(mg)^v)$ with $J(u\times v\times mg,3,1)$ base blocks.
\end{lemma}
\proof
By assumption, there exists a $\{3,5,v\}$-SCGDD of type $g^v$ with one base block of size $v$. Let the SCGDD have $a$ base blocks of size $3$ and $b$ base blocks of size $5$.
It is immediate  that $3a+10b=\frac{v(v-1)(g-1)}{2}$. Apply Corollary~\ref{col4} with an $m$-cyclic $3$-HGDP of type $(u,m^v)$ with $J(u\times v\times m,3,1)$ base blocks and an $m$-cyclic $3$-HGDD of type $(u,m^k)$ for $k\in\{3,5\}$  which exists by Theorem~\ref{Yt11}. Then we obtain an $mg$-cyclic $3$-HGDP of type $(u,(mg)^v)$. It has
$$\lfloor\frac{uvm(u-1)(v-1)-uv}{6}\rfloor+a\times mu(u-1)+b\times\frac{10mu(u-1)}{3}$$
$$=\lfloor\frac{uvmg(u-1)(v-1)-uv}{6}\rfloor=J(u\times v\times mg,3,1)$$ base blocks.
\qed
\begin{lemma}\label{z8}
Let $u,v,g$ be even, $w$ be odd,  $ug\equiv 0~(\bmod~3)$ and $v\equiv t~(\bmod~g)$.
Suppose that there exists a $3$-GDD of type $g^{\frac{v-t}{g}}t^1$. If there exists  a $w$-cyclic $3$-HGDP of type $(u,w^k)$ with $J(u\times k\times w,3,1)$ base blocks for each $k\in\{g,t\}$, then there exists a $w$-cyclic $3$-HGDP of type $(u,w^v)$ with $J(u\times v\times w,3,1)$ base blocks.
\end{lemma}
\proof
By assumption, there exists a $3$-GDD of type $g^{\frac{v-t}{g}}t^1$ with $\frac{(v-t)(v+t-g)}{6}$ blocks, which is also a $\{3,g,t\}$-SCGDD of type $1^v$ containing $\frac{(v-t)(v+t-g)}{6}$ blocks of size 3, ${v-t\over g}$ blocks of size $g$ and one block of size $t$. Apply Corollary~\ref{col4} with a $w$-cyclic $3$-HGDP of type $(u,w^k)$ with $J(u\times k\times w,3,1)$ base blocks for $k\in\{g,t\}$ and a $w$-cyclic $3$-HGDD of type $(u,w^3)$  which exists from Theorem~\ref{Yt11}. Then we obtain a $w$-cyclic $3$-HGDP of type $(u,w^v)$. The number of base blocks equals
$$\frac{(v-t)(v+t-g)}{6}\times u(u-1)w+\frac{v-t}{g}\times \frac{ugw(u-1)(g-1)-ug}{6}+\lfloor\frac{ut((u-1)(t-1)w-1)}{6}\rfloor$$
$$=\lfloor\frac{uv((u-1)(v-1)w-1)}{6}\rfloor=J(u\times v\times w, 3,1).$$\qed

\section{Case $u\equiv 0~(\bmod~3)$}
In this section,
we shall determine the number of base blocks of a maximum $w$-cyclic $3$-HGDP of type $(u,w^v)$~when $u\equiv 0~(\bmod~3)$.
\begin{lemma}{\rm[{\cite{GJL}}, Theorem 2]}\label{scgdd}
A $3$-SCGDD of type $m^n$ exists if and only if $n\geq 3$ and
$(1)~n\equiv1,3~(\bmod~6)$ when $m\equiv1,5~(\bmod~6);$
$(2)~n\equiv1~(\bmod~2)$ when $m\equiv3~(\bmod~6);$
$(3)~n\equiv0,1,4,9~(\bmod~12)$ when $m\equiv2,10~(\bmod~12);$
$(4)~n\equiv0,1~(\bmod~3),~n\neq3$ when $m\equiv4,8~(\bmod~12);$
$(5)~n\neq3$ when $m\equiv0~(\bmod~12);$
$(6)~n\equiv0,1~(\bmod~4)$ when $m\equiv6~(\bmod~12)$.
\end{lemma}

\begin{lemma}\label{33w}
$\Psi(3\times 3\times w,3,1)\leq 6(w-1)$ for  $w\equiv0~(\bmod~2)$.
\end{lemma}
\proof
Let $X=I_3\times I_3\times Z_w$, $\mathcal{G}=\{\{i\}\times I_3\times  Z_w:i\in I_3\}$ and $\mathcal{H}=\{I_3\times\{j\}\times Z_w:j\in I_3\}$.
Suppose that $(X,\mathcal{G},\mathcal{H},\mathcal{B})$ is a $w$-cyclic 3-HGDP of type $(3,w^3)$ with at least $6w-5$ base blocks. Obviously, its base blocks can be partitioned into the following six possible types:

$(1)~\{(0,0,\ast),(1,1,\ast),(2,2,\ast)\}$, $(2)~\{(0,0,\ast),(1,2,\ast),(2,1,\ast)\},$
$(3)~\{(0,1,\ast),(1,0,\ast),(2,2,\ast)\}$,

$(4)~\{(0,1,\ast),(1,2,\ast),(2,0,\ast)\}$,
$(5)~\{(0,2,\ast),(1,0,\ast),(2,1,\ast)\}$, $(6)~\{(0,2,\ast),(1,1,\ast),(2,0,\ast)\}$.

\noindent
Using Pigeonhole principle guarantees that there exists a type containing at least $w$ base blocks. However, it is easy to check that
 every type  contains at most $w$ base blocks since the base blocks of this type omitting the first coordinate produce
a 3-SCGDP of type $w^3$. As a result, we have a type containing exactly $w$ base blocks which form a $3$-SCGDD of type $w^3$.
But Lemma~\ref{scgdd} shows that there is no 3-SCGDD of type $w^3$ for any $w\equiv0~(\bmod~2).$ A contradiction occurs. The conclusion then follows.
\qed
\begin{lemma}\label{33}
$\Psi(3\times3\times w,3,1)=6(w-1)$ for $w\equiv0~(\bmod~2)$.
\end{lemma}
\proof
By Theorem~\ref{hgdd}, there exists a 3-HGDD of type $(3,1^3)$ with six blocks. Apply Construction~\ref{xc3} with a 3-SCGDP of type $w^3$ with $w-1$ base blocks which exists from  [\cite[Lemma~5.1]{WC15}] to  obtain  a $w$-cyclic 3-HGDP of type $(3,w^3)$ with $6(w-1)$ base blocks. So the conclusion follows by Lemma~\ref{33w}.
\qed

\begin{lemma}\label{4,1^6}
There exists a $3$-HGDP of type $(4,1^6)$ with $J(4\times6\times1,3,1)$ blocks.
\end{lemma}

\proof We construct the desired design on $X=Z_{4}\times I_6$ with group set $\mathcal{G}=\{\{i\}\times I_6:i\in Z_4\}$ and hole set $\mathcal{H}=\{Z_4\times\{j\}:j\in I_6\}$.
All $56$ blocks can be obtained by developing the following $14$ base blocks by $(+1,-)$ mod $(4,-)$ successively.

\begin{minipage}[t]{0.45\textwidth}
\centering
\begin{tabular}{llllll}
$\{(0,0),(1,2),(2,1)\}$,&
$\{(0,0),(1,3),(2,2)\}$,&
$\{(0,0),(1,4),(3,1)\}$,&
$\{(0,0),(1,5),(3,2)\}$,\\
$\{(0,0),(2,3),(3,4)\}$,&
$\{(0,0),(2,4),(3,5)\}$,&
$\{(0,0),(2,5),(3,3)\}$,&
$\{(0,1),(1,2),(3,4)\}$,\\
$\{(0,1),(1,3),(2,5)\}$,&
$\{(0,1),(1,4),(2,3)\}$,&
$\{(0,1),(1,5),(3,3)\}$,&
$\{(0,1),(2,2),(3,5)\}$,\\
$\{(0,2),(1,3),(3,4)\}$,&
$\{(0,2),(1,4),(3,5)\}$.& &
\end{tabular}
\end{minipage}

\qed
\begin{lemma}\label{4}
Let $w\equiv1~(\bmod~2)$. Then there exists a $w$-cyclic $3$-HGDP of type $(u,w^v)$ with $J(u\times v\times w,3,1)$
base blocks for $(1)$ $u\in\{6,12\}$ and $v\in\{4,6,10\}$, $(2)$ $(u,v)\in\{(8,12),(14,12)\}.$
\end{lemma}
\proof When $(u,v,w)=(6,6,3)$, the desired design can be obtained from Example B.8 in Supporting Information 
directly.

For other parameters, by Theorem~\ref{1} there exists a 3-SCHGDD of type $(v,1^w)$ which can be viewed as a $\{3,v\}$-SCGDD of type $w^{v}$ with $\frac{v(v-1)(w-1)}{6}$ base blocks of size $3$ and one base block of size $v$ by Lemma \ref{eq}. From Lemma~\ref{4,1^6}, there exists a $3$-HGDP of type $(4,1^6)$, also of type $(6,1^4)$, with $J(6\times 4\times 1,3,1)$ blocks. For the rest parameters, Lemmas B.1-B.5 in Supporting Information provide  3-HGDPs of type $(u,1^v)$ with $J(u\times v\times 1,3,1)$ base blocks.
 Apply Lemma~\ref{z7} by taking parameters $(g,v,u,m)=(w,v,u,1)$ to yield a $w$-cyclic 3-HGDP of type $(u,w^v)$ with $J(u\times v\times w,3,1)$ base blocks.
\qed
\begin{lemma}\label{5}
Let $w\equiv1~(\bmod~2)$. There exists a $w$-cyclic $3$-HGDP of type $(6,w^v)$ with $J(6\times v\times w,3,1)$ base blocks for $v\in\{8,14\}$.
\end{lemma}
\proof
For $v\in\{8,14\}$, by Lemma A.1 in Supporting Information, there exists a $\{3,6\}$-HGDD of type $(6,1^v)$ with $ 5v(v-2)$ blocks of size $3$ and $v$ blocks of  size $6$,
 where the blocks of size 6 form a parallel class. Apply Construction~\ref{xc3} with a 3-SCGDD of type $w^{3}$ with $w$ base blocks which exists from Lemma~\ref{scgdd} and a $3$-SCGDP of type $w^6$ with 
 $5w-1$ base blocks to obtain  a $w$-cyclic 3-HGDP of type $(6,w^v)$ with $w\times (5v(v-2))+v\times(5w-1)=J(u\times v\times w,3,1)$ base blocks, where the required 3-SCGDP of type $w^{6}$ can be obtained from  [\cite[Lemma 5.1]{WC15}].
\qed
\begin{theorem}\label{YL1}{\rm[\cite{CDR}, Main Theorem]}
Let $g,t$ and $w$ be positive integers. Then there exists a $3$-GDD of type $g^{t}w^{1}$ if and only if all the following conditions are satisfied:
  $(1)$~$t\geq3$, or $t=2$ with $w=g;$
  $(2)$~$w\leq g(t-1);$
  $(3)$~$g(t-1)+w\equiv 0~(\bmod~2);$
  $(4)$~$gt\equiv 0~(\bmod~2);$
  $(5)$~$g^{2}t(t-1)+2gtw\equiv0~(\bmod~6)$.
\end{theorem}

\begin{lemma}\label{t1}
Let
$u\in\{6,12\}$, $v\equiv0~(\bmod~2)~(v\geq4)$ and $w\equiv1~(\bmod~2)$.
Then there exists a $w$-cyclic $3$-HGDP of type $(u,w^v)$ with $J(u\times v\times w, 3,1)$ base blocks.
\end{lemma}
\proof
 For $v\in\{4,6,8,10,14\}$, the desired designs exist from~Lemmas~\ref{4} and~\ref{5}. Note that a 3-HGDP of type $(u,1^{v})$ is also that of type $(v,1^{u})$.

Let $g=4$ if $v\equiv0~(\bmod~4)~(v\geq12)$ and let $g=6$ if $v\equiv6~(\bmod~12)~(v\geq18)$. From Theorem~\ref{1} and Lemma~\ref{c15} there exists a $w$-cyclic $3$-HGDD of type $(\frac{v}{g},(gw)^{u})$ for $u\in\{6,12\}$ and odd $w$. By Lemma~\ref{4} there exists a $w$-cyclic $3$-HGDP of type $(u,w^g)$, also of type $(g,w^u)$, with $J(g\times u\times w,3,1)$ base blocks.
Apply Lemma~\ref{zz} by taking  $g=s$ to obtain a $w$-cyclic $3$-HGDP of  type $(v,w^u)$, also of type $(u,w^v)$,  with $J(u\times v\times w,3,1)$ base blocks.

 Let $(g,t)=(6,8)$ if $v\equiv2~(\bmod~12)~(v\geq26)$  and let $(g,t)=(4,6)$ if $v\equiv10~(\bmod~12)~(v\geq22)$. By Theorem~\ref{YL1} there exists a $3$-GDD of type $g^{\frac{v-t}{g}}t^1$. For $u\in\{6,12\}$ and odd $w$,
apply Lemma~\ref{z8} to obtain a $w$-cyclic $3$-HGDP of  type $(u,w^v)$ with $J(u\times v\times w,3,1)$ base blocks, where the required maximum $w$-cyclic 3-HGDPs of types $(u,w^g)$ and $(u,w^t)$ exist  by the previous arguments. \qed
\begin{lemma}\label{3.8}
$$\Psi(u\times v\times w,3,1)=\left\{
  \begin{array}{ll}
   6(w-1), & \hbox{if  $(u,v)=(3,3)$ and $w\equiv 0~(\bmod~2)$, } \\
   J(u\times v\times w,3,1), & \hbox{if $u\equiv 0~(\bmod~3)$,~$u,v\geq3$~and~$w\geq1$.}
  \end{array}
\right.$$
\end{lemma}
\proof
The assertion holds for $(u,v)=(3,3)$ and $w\equiv0~(\bmod~2)$ by Lemma~\ref{33}. Then, by Corollary~\ref{co1}, we only need to handle the case that $u\equiv0~(\bmod~6)$, $v\equiv0~(\bmod~2)$ and $w\equiv1~(\bmod~2)$. Lemma~\ref{Johnson bound}  shows  $\Psi(u\times v\times w,3,1)\leq J(u\times v\times w,3,1)$.

For $u\in\{6,12\}$, by Lemma~\ref{t1} we have $\Psi(u\times v\times w,3,1)= J(u\times v\times w,3,1)$.

For $u\equiv0~(\bmod~6)$ and $u\geq18$,
by Theorem~\ref{1}, there exists a 3-SCHGDD of type $(\frac{u}{6},(6w)^{v})$, which is also a $w$-cyclic 3-HGDD of type $(\frac{u}{6},(6w)^{v})$ by Lemma~\ref{c15}. There exists a $w$-cyclic $3$-HGDP of type $(6,w^v)$ by Lemma~\ref{t1}.
Apply Lemma~\ref{zz}  by taking parameters $(g,s,v,u,w)=(6,6,u,v,w)$
to obtain a $w$-cyclic $3$-HGDP of type $(u,w^v)$ with $J(u\times v\times w,3,1)$ base blocks. So we have $\Psi(u\times v\times w,3,1)= J(u\times v\times w,3,1).$
\qed
\section{Auxiliary designs}

In this section, we study two types of auxiliary designs, namely incomplete HGDPs and semi-cyclic $K$-GDDs, which will play important roles in constructing $w$-cyclic $3$-HGDPs of type $(u,w^v)$ when $u\equiv1~(\bmod~3)$.
\subsection{Incomplete holely group divisible packings}

In this subsection we focus on incomplete $3$-HGDPs of type $(u,t,w^v)$ which contain the largest possible number of blocks.
Let $u,v,t$ be even. It is not difficult to check that the number of blocks of a $3$-IHGDP of type $(u,t,w^v)$ does not exceed
\begin{equation}
\Theta(u,t, v, w)=\left\{
  \begin{array}{ll}
   \frac{vw^2(u-t)(u+t-1)(v-1)}{6}, & \hbox{if $w$ is even,} \\
   \frac{vw(u-t)((u+t-1)(v-1)w-1)}{6}, & \hbox{if $w$ is odd.}
   \end{array}
   \right.
  \end{equation}
Furthermore, if $w$ is even, then a $3$-IHGDP of type $(u,t,w^v)$ with $\Theta(u,t, v, w)$ blocks is actually a $3$-IHGDD. In section 5 we shall employ $3$-IHGDPs of type $(u,t,1^v)$ with $\Theta(u,t, v, 1)$ blocks to produce maximum $3$-HGDPs of type $(u,1^v)$, which is $w$-cyclic with $w=1.$


\begin{lemma}\label{z5}
Let $u,v,t$ be even. If there exist a $3$-IHGDP of type $(u,t,1^v)$ with $\Theta(u, t, v, 1)$ blocks and a $3$-HGDP of type $(t,1^v)$ with $J(t\times v\times 1,3,1)$ blocks, then there exists a $3$-HGDP of type $(u,1^v)$ with $J(u\times v\times 1,3,1)$ blocks.
\end{lemma}
\proof
By assumption,
there exists a $3$-IHGDP of $(u,t,1^v)$ with $\Theta(u,t,v,1)$ blocks. Apply Corollary~\ref{col2} with a $3$-HGDP of type $(t,1^v)$ with $J(t\times v\times 1,3,1)$ blocks to obtain a $3$-HGDP of type $(u,1^v)$. It has $$\frac{v(u-t)((u+t-1)(v-1)-1)}{6}+\lfloor\frac{tv((t-1)(v-1)-1)}{6}\rfloor$$
$$=\lfloor\frac{uv(u-1)(v-1)-uv}{6}\rfloor=J(u\times v\times 1,3,1)$$ blocks.
\qed

\begin{lemma}\label{z2}
Let $m,t,v$ be even.
Suppose that there exists a $3$-HGDD of type $(u,(mw)^v)$.

$(1)$ Let $w$ be odd. If there exists a $3$-IHGDP of type $(m+t,t,w^v)$ with $\Theta(m+t,t,v,w)$ blocks, then there exist a $3$-IHGDP of type $(um+t,m+t,w^v)$ with $\Theta(um+t,m+t, v, w)$ blocks and a $3$-IHGDP of type $(um+t,t,w^v)$ with $\Theta(um+t,t, v, w)$ blocks.

$(2)$ Let $w$ be even. If there exists a $3$-IHGDD of type $(m+t,t,w^v)$, then there exist a $3$-IHGDD of type $(um+t,m+t,w^v)$ and a $3$-IHGDD of type $(um+t,t,w^v)$.
\end{lemma}
\proof
By assumption, there exists a 3-HGDD of type $(u, (mw)^v)$ with  $\frac{uv(u-1)(v-1)m^2w^2}{6}$ blocks. Apply Corollary~\ref{col2} with a
$3$-IHGDP of type $(m+t,t,w^v)$ which has $\Theta(m+t,t,v,w)$ blocks to obtain a $3$-IHGDP of type $(um+t,m+t,w^v)$. It is easy to check from Equation $(1)$ that the number of blocks equals $$\frac{uv(u-1)(v-1)m^2w^2}{6}+(u-1)\times\Theta(m+t,t,v,w)
=\Theta(um+t,m+t,v,w).$$
Also, we further obtain a $3$-IHGDP of type $(um+t,t,w^v)$ with $\Theta(um+t,t, v, w)$
 blocks.

In the case that $w$ is even, we can readily check that the produced $3$-IHGDPs become $3$-IHGDDs. This completes the proof.\qed
\begin{lemma}\label{z3}
Let $u,t,v$ be even.
If there exist a $3$-IHGDD of type $(u,t,w^v)$ and a $3$-IHGDP of type $(u,t,1^w)$ with $\Theta(u,t,w,1)$ blocks, then there exists a $3$-IHGDP of type $(u,t,1^{vw})$ with $\Theta(u,t, vw, 1)$ blocks.
\end{lemma}
\proof
Since there exists a $3$-IHGDD of type $(u,t,w^v)$ with $\frac{vw^2(u-t)(u+t-1)(v-1)}{6}$ blocks, apply Construction~\ref{xc4} with a $3$-IHGDP of type $(u,t,1^w)$ with $\Theta(u,t, w, 1)$ blocks to  obtain a $3$-IHGDP of type $(u,t,1^{vw})$. It has $$\frac{vw^2(u-t)(u+t-1)(v-1)}{6}+v\times\frac{w(u-t)((u+t-1)(w-1)-1)}{6}=\Theta(u,t,vw,1)$$ blocks.
\qed

  Now we construct two small examples of $3$-IHGDPs of type $(u,t,w^v)$ with $\Theta(u,t, v, w)$ blocks.
\begin{example}\label{e2}
We construct a $3$-IHGDD of type $(8,2,4^4)$ with $\Theta(8,2,4,4)$ blocks
 on $(Z_{6}\cup\{a,b\})\times Z_{4}\times Z_{4}$ with group set $\{\{i\}\times Z_{4}\times Z_{4}:i\in Z_{6}\}\cup\{\{a,b\}\times Z_{4}\times Z_{4}\}$ and hole set $\{(Z_{6}\cup\{a,b\})\times \{j\}\times Z_{4}:j\in Z_{4}\}$. All $1728$ blocks are obtained by developing the following $36$ base blocks by $(+2,+1,+1)\bmod~ (6,4,4)$ successively, where $a+2=a$ and $b+2= b$.

\begin{longtable}{lllllllll}
$\{( 0, 0, 0 ), ( 2, 1, 0 ), ( a, 2, 0 )\}$, & $\{( 0, 0, 0 ), ( a, 1, 1 ), ( 1, 2, 1 )\}$, &
 $\{( 0, 0, 0 ), ( 1, 3, 0 ), ( a, 2, 1 )\}$, \\ $\{( 0, 0, 0 ), ( 1, 2, 0 ), ( a, 1, 2 )\}$, &
 $\{( 0, 0, 0 ), ( 1, 1, 0 ), ( a, 2, 3 )\}$, & $\{( 0, 0, 0 ), ( 3, 2, 0 ), ( a, 1, 3 )\}$,\\
 $\{( 0, 0, 0 ), ( 3, 1, 0 ), ( a, 2, 2 )\}$, & $\{( 0, 0, 0 ), ( 5, 1, 0 ), ( a, 3, 1 )\}$, &
 $\{( 0, 0, 0 ), ( 1, 1, 2 ), ( a, 3, 2 )\}$, \\ $\{( 0, 0, 0 ), ( a, 3, 0 ), ( 1, 1, 1 )\}$, &
 $\{( 0, 0, 0 ), ( 1, 2, 3 ), ( a, 3, 3 )\}$, & $\{( 1, 0, 0 ), ( 3, 1, 1 ), ( a, 2, 2 )\}$, \\
 $\{( 0, 0, 0 ), ( 4, 1, 0 ), ( 5, 3, 2 )\}$, & $\{( 0, 0, 0 ), ( b, 1, 0 ), ( 5, 3, 3 )\}$, &
 $\{( 0, 0, 0 ), ( 2, 2, 0 ), ( 4, 3, 1 )\}$, \\ $\{( 0, 0, 0 ), ( 5, 2, 0 ), ( 2, 1, 2 )\}$, &
 $\{( 0, 0, 0 ), ( b, 2, 0 ), ( 5, 3, 0 )\}$, & $\{( 0, 0, 0 ), ( 3, 1, 1 ), ( 5, 3, 1 )\}$, \\
 $\{( 0, 0, 0 ), ( 2, 2, 1 ), ( 3, 1, 3 )\}$, & $\{( 0, 0, 0 ), ( 5, 2, 1 ), ( 1, 3, 1 )\}$, &
 $\{( 0, 0, 0 ), ( 3, 2, 1 ), ( 1, 1, 3 )\}$, \\ $\{( 1, 0, 0 ), ( b, 1, 0 ), ( 5, 3, 1 )\}$, &
 $\{( 0, 0, 0 ), ( 4, 1, 1 ), ( b, 3, 2 )\}$, & $\{( 0, 0, 0 ), ( 3, 3, 0 ), ( 4, 1, 2 )\}$, \\
 $\{( 0, 0, 0 ), ( 3, 3, 1 ), ( b, 2, 3 )\}$, & $\{( 0, 0, 0 ), ( 2, 3, 1 ), ( 5, 2, 3 )\}$, &
 $\{( 0, 0, 0 ), ( b, 2, 2 ), ( 5, 1, 3 )\}$, \\ $\{( 0, 0, 0 ), ( 4, 2, 1 ), ( b, 1, 2 )\}$, &
 $\{( 0, 0, 0 ), ( 3, 2, 2 ), ( b, 1, 3 )\}$, & $\{( 0, 0, 0 ), ( b, 3, 0 ), ( 3, 2, 3 )\}$, \\
 $\{( 0, 0, 0 ), ( 5, 1, 2 ), ( 1, 3, 3 )\}$, & $\{( 1, 0, 0 ), ( b, 2, 0 ), ( 5, 1, 2 )\}$, &
 $\{( 1, 0, 0 ), ( 5, 1, 1 ), ( b, 3, 3 )\}$, \\ $\{( 0, 0, 0 ), ( 5, 1, 1 ), ( 3, 3, 3 )\}$, &
 $\{( 0, 0, 0 ), ( b, 1, 1 ), ( 4, 2, 2 )\}$, & $\{( 1, 0, 0 ), ( 5, 1, 0 ), ( 3, 2, 3 )\}.$
\end{longtable}
\qed
\end{example}

\begin{example}\label{e3}
We construct a $3$-IHGDP of type $(8,2,1^4)$ with $\Theta(8,2,4,1)$ blocks on $X=I_{8}\times Z_{4}$ with group set $\mathcal{G}=\{\{i\}\times Z_{4}:i\in I_{6}\}\cup\{\{6,7\}\times Z_{4}\}$ and hole set
$\mathcal{H}=\{I_{8}\times\{j\}:j\in Z_{4}\}$. All $104$ blocks are obtained by developing the following
$26$ base blocks by $(-,+1)\bmod(-,4)$ successively.
\begin{longtable}{llllllll}

$\{(0, 0 ), ( 6, 1 ), ( 1, 2 )\}$, & $\{(0, 0 ), ( 6, 2 ), ( 2, 3 )\}$,&  $\{(0, 0 ), ( 4, 2 ), ( 6, 3 )\}$,&
 $\{(1, 0 ), ( 6, 1 ), ( 3, 3 )\}$,  \\$\{(2, 0 ), ( 6, 1 ), ( 5, 3 )\}$,&  $\{(1, 0 ), ( 6, 2 ), ( 5, 3 )\}$,&
 $\{(2, 0 ), ( 6, 2 ), ( 4, 3 )\}$, & $\{(3, 0 ), ( 6, 1 ), ( 4, 3 )\}$,  \\ $\{(3, 0 ), ( 5, 2 ), ( 6, 3 )\}$,&
 $\{(1, 0 ), ( 3, 1 ), ( 2, 3 )\}$, & $\{(0, 0 ), ( 2, 1 ), ( 1, 3 )\}$, & $\{(1, 0 ), ( 2, 1 ), ( 3, 2 )\}$,\\
 $\{(1, 0 ), ( 7, 2 ), ( 4, 3 )\}$,&  $\{(2, 0 ), ( 7, 1 ), ( 3, 3 )\}$,&  $\{(2, 0 ), ( 4, 1 ), ( 7, 2 )\}$,&
 $\{(2, 0 ), ( 5, 1 ), ( 4, 2 )\}$,\\  $\{(0, 0 ), ( 7, 1 ), ( 2, 2 )\}$,&  $\{(0, 0 ), ( 1, 1 ), ( 4, 3 )\}$,&
 $\{(1, 0 ), ( 5, 1 ), ( 7, 3 )\}$,&  $\{(1, 0 ), ( 7, 1 ), ( 5, 2 )\}$,\\  $\{(3, 0 ), ( 4, 2 ), ( 5, 3 )\}$,&
 $\{(0, 0 ), ( 5, 1 ), ( 7, 2 )\}$,&  $\{(0, 0 ), ( 3, 2 ), ( 7, 3 )\}$,&  $\{(0, 0 ), ( 4, 1 ), ( 5, 3 )\}$,\\
 $\{(0, 0 ), ( 3, 1 ), ( 5, 2 )\}$,&  $\{(3, 0 ), ( 4, 1 ), ( 7, 3 )\}$.&&
\end{longtable}\qed
\end{example}
\begin{lemma}\label{sc2}
There exists a $3$-IHGDP of type $(20,8,1^4)$ with $\Theta(20,8,4,1)$ blocks.
\end{lemma}
\proof
There exist a 3-HGDD of type
  $(3,6^4)$ and a 3-IHGDP of type $(8,2,1^4)$ from Theorem \ref{hgdd} and Example~\ref{e3}, respectively.
 Apply Lemma~\ref{z2} by taking parameters $(u,v,w,m,t)=(3,4,1,6,2)$ to obtain a $3$-IHGDP of type $(20,8,1^4)$ with $\Theta(20,8,4,1)$ blocks.
\qed
\begin{lemma}\label{z4}
Let $(x,y,z)\in\{(10,4,16),(16,4,16),(20,8,16),(10,4,20)\}$. Then there exists a $3$-IHGDP of type $(x,y,1^z)$ with $\Theta(x,y, z,1)$ blocks.
\end{lemma}
\proof
By~Theorem~\ref{hgdd}, there exist $3$-HGDDs of types $(4,4^4)$ and $(4,4^5)$ which are also $3$-IHGDDs of types $(4,1,4^4)$ and $(4,1,4^5)$, respectively. There exists a $3$-IHGDD of type $(8,2,4^4)$ from Example~\ref{e2}. 
Apply Lemma~\ref{z2} by taking parameters $(u,v,w,m,t)$ in the second column of Table I. Then we obtain a 3-IHGDD of type $(x,y,w^v)$, where $x=um+t$, $y=m+t$ and $z=vw$. The needed $3$-HGDDs of type $(u,(mw)^v)$ exist from Theorem~\ref{hgdd}.
 Furthermore, $3$-IHGDPs of type $(x,y,1^w)$ with $\Theta(x,y, w,1)$ blocks exist, for which we list the sources  in the third column of Table I (see Lemma C.1 in Supporting Information). Then by Lemma~\ref{z3}, we obtain the desired $3$-IHGDPs of type $(x,y,1^z)$.
\begin{longtable}{ccccc}
\caption*{Table I}\\
\hline
$(x,y,z)$ & $(u,v,w,m,t)$  & $(x,y,1^w)$\\
\hline
(10,4,16) & (3,4,4,3,1) &Lemma~C.1\\
(16,4,16) & (5,4,4,3,1) &Lemma~C.1\\
(20,8,16) & (3,4,4,6,2) &Lemma~\ref{sc2}\\
(10,4,20) & (3,5,4,3,1)&Lemma C.1\\
\hline
\end{longtable}

\qed

\subsection{$K$-SCGDDs}

In this subsection we  establish the existence of a class of $\{3,5,v\}$-SCGDDs of  type $w^v$, which will be applied in Lemma~\ref{z7} to produce maximum $w$-cyclic $3$-HGDPs. Note that Lemma \ref{eq} shows the relationship between $K$-SCHGDDs and $L$-SCGDDs, so we form $\{3,5,v\}$-SCGDDs of  type $w^v$ by constructing $\{3,5\}$-SCHGDDs of type $(v,1^w)$.

Cyclotomic cosets play an important role in direct constructions for SCHGDDs. Let $p\equiv 1~(\bmod~n)$ be a prime and $w$ be a primitive element of $Z_p$. Let $C_{0}^n$ denote the multiplicative subgroup $\{w^{in}:0\leq i<(p-1)/n\}$ of the $n$-th powers in $Z_p$ and $C_{j}^{n}$ denote the coset of $C_{0}^{n}$ in $Z_p{\setminus}\{0\},$ i.e., $C_{j}^{n}=w^j\cdot C_{0}^{n},~0\leq j\leq n-1.$

\begin{lemma}{\rm[\cite{fwc}, Lemma 6.8]}\label{squre}
Let $p\geq5$ be a prime. There exists an element $x\in Z_p$ such that $x\in C_{1}^2$, $x+1\in C_{1}^{2}$ and $x-1\in C_{0}^{2}.$
\end{lemma}

\begin{lemma}\label{81p}
There exists a $\{3,5\}$-SCHGDD of type $(8,1^p)$ for any prime $p\geq5.$
\end{lemma}

\proof
 Let $\omega$ be a primitive element of $Z_p$ where $p\geq 5$.
The desired design will be constructed on $X=I_{8}\times Z_{p}$ with group set $\mathcal{G}=\{\{i\}\times Z_{p}:i\in I_{8}\}$ and hole set $\mathcal{H}=\{I_8\times\{j\}:j\in Z_{p}\}$. By Lemma~\ref{squre}, we  may take $x\in Z_p$ such that $x\in C_1^2,~x+1 \in C_1^2$ and $x-1\in C_0^2$.

When $p\equiv1~(\bmod~4),$ we construct twelve initial base blocks of size $3$ and two initial base blocks of size $5$.  All base blocks are generated by multiplying the second coordinate of these initial base blocks by $\omega^{2r}$, $0\leq r \leq (p-3)/2.$ The initial base blocks are listed in four cases.

(1) $2, x+2 \in C_0^2.$
\begin{longtable}{lll}
$\{(0, 0), (1, x), (2, 1)\}$,&
$\{(1,0), (2,x), (5,x+1)\}, $&
$\{(1,0), (3,x-1), (7,x)\}, $\\
$\{(4,0), (1,1), (7,x)\},$
& $\{(2,0), (6,x), (7,x+1)\}, $&
$\{(3,0), (5,x), (7,x+1)\},$\\
$\{(0, 0), (3, 1), (5, x)\},$&
$\{(4,0), (5,x), (6,x+1)\},$&
$\{(1,0), (5,1), (6,x+1)\}, $\\
$\{(0, 0), (6, -1), (7, x)\},$&
$\{(3,0), (2,x), (6,x+1)\}$,& $\{(4,0), (2,x), (3,x+1)\};$\smallskip\\
\multicolumn{3}{l}{$\{(0,-x),(4,0),(5,1-x),(2,1),(7,2)\}$,
$\{(0,0), (1,1),(3,x+1),(4,x+2), (6,x)\}.$}

\end{longtable}

(2) $2, x+2 \in C_1^2.$
\begin{longtable}{lll}
$\{(3,1), (5,x), (7,0)\},$
&$\{(1,0), (5,1), (6,x+1)\},$&
$\{(3,0), (2,x), (6,x+1)\},$\\
$\{(4,0), (5,x), (6,1)\},$&
$\{(4,1-x), (2,0), (3,1)\},$&
$\{(1,0), (2,x), (5,x+1)\},$\\
$\{(2,0), (6,x), (7,1)\},$&
$\{(0,0), (1,x),\ (2,x+1)\},$& $\{(0,0), (3,1-x), (5,1)\},$  \\
$\{(0,0), (6,-1), (7,x-1)\},$&
$\{(1,0), (3,x-1), (7,-1)\},$&
$\{(4,0), (1,1), (7,1-x)\};$\smallskip\\
\multicolumn{3}{l}{$\{(0,0), (1,1), (3,x+1), (4,x+2), (6,x)\}$, $
\{(0,1), (4,0), (5,1-x), (2,2), (7,-x)\}.$}
\end{longtable}
(3) $2\in C_1^2,~x+2 \in C_0^2.$
\begin{longtable}{llll}
$\{(1,0), (5,1), (6,x+1)\},$&
$\{(3,1), (5,x), (7,0)\},$&$\{(3,0), (2,x), (6,x+1)\},$\\
$\{(4,0), (1,1), (7,1-x)\},$&
$\{(4,0), (5,x), (6,1)\},$&
$\{(4,1+x), (2,0), (3,1)\},$\\
$\{(1,0), (2,x), (5,x+1)\},$&
$\{(2,0), (6,x), (7,1)\},$&
$\{(1,0), (3,x-1), (7,-1)\},$\\
$\{(0,0),(1,x),(2,x+1)\},$&$\{(0,0), (3,1-x),(5,1)\},$&
$\{(0,0),(6,-1),(7,-x-1)\};$\smallskip\\
\multicolumn{3}{l}{$\{(0,0), (1,1), (3,x+1), (4,x+2), (6,x)\},
 \{(0,-2), (4,0), (5,x-1), (2,-1), (7,x)\}.$}
\end{longtable}
(4) $2\in C_0^2, x+2 \in C_1^2.$
\begin{longtable}{lll}
$\{(0,0), (1,x), (2,x+1)\},$& $\{(0,0), (3,1), (5,2)\}, $& $\{(0,0), (6,-1), (7,x-1)\},$\\
$\{(1,0), (2,x), (5,x+1)\},$&
$\{(2,0), (6,x), (7,1)\},$&
$\{(1,0), (3,x-1), (7,-1)\},$\\
$\{(1,0), (5,1), (6,x+1)\},$&
$\{(3,-1), (5,x), (7,0)\},$&
$\{(3,0), (2,x), (6,x+1)\},$\\
$\{(4,0), (1,1), (7,1-x)\},$&
$\{(4,0), (2,x), (3,x+1)\},$&
$\{(4,0), (5,x), (6,x+1)\};$\smallskip\\
\multicolumn{3}{l}{$\{(0,0), (1,1), (3,x+1), (4,x+2), (6,x)\},
\{(0,1), (4,0), (5,1-x), (2,2), (7,-x)\}.$}
\end{longtable}

When $p\equiv3~(\bmod~4),$ we construct six initial base blocks of size $3$ and one initial base block of size $5$.
All base blocks can be obtained by successively developing these initial base blocks by $(+4,-)\bmod (8,-)$ and multiplying the second coordinate by $\omega^{2r}$, $0\leq r \leq (p-3)/2.$ We list the initial base blocks by considering two cases.

Case $1.$ $x+2\neq p$.  Four subcases are distinguished.

(1) $-2\in C_1^2, x+2 \in C_0^2.$
\begin{longtable}{lll}
$\{(0, 0), (1, x), (2, x-1)\}, $& $\{(0, 0), (3, x), (5, -1)\},$ & $\{(1, 0), (2, 1), (5, x)\},$\\
$\{(1, 0), (3, -x-1), (7, -x)\},$& $\{(2, 0), (3, x), (6, x+1)\},$& $\{(0, 0), (6, 1), (7, 1-x)\};$\smallskip\\
\multicolumn{3}{l}{$\{(0, 0), (1, 1), (3, x+2), (6, x), (4, x+1)\}.$}
\end{longtable}
(2) $-2,~x+2\in C_1^2.$
\begin{longtable}{lll}
\centering
$\{(0, 0), (1, x-1), (2, x+1)\},$&
 $\{(0, 0), (3, x), (5, x+1)\},$ &
$\{(0,0),(6,-x-1),(7,1)\},$\\
$\{(1, 0), (2, x+1), (5, x-1)\},$
&
$\{(1, 0), (3, -1), (7, x-1)\}$,& $\{(2, 0), (3, x-1 ), (6, -1)\};$\smallskip\\
\multicolumn{3}{l}{$\{(0, 0), (1, x+2), (3, 1), (6, x+1), (4, x)\}.$ }
\end{longtable}
(3) $-2\in C_0^2, x+2 \in C_1^2.$
\begin{longtable}{lll}
\centering
$\{(1, 0),(2, 1),(5, x)\},$
&
$\{(0, 0),(3, x-1),(5, -1)\},$&
$\{(0, 0),(6, 1 ),(7, 1-x)\},$\\
$\{(0, 0),(1, x),(2, x-1)\},$&
$\{(1, 0),(3, -x-1),(7, -x)\},$&
$\{(2, 0),(3, x),(6, x-1)\};$\smallskip\\
\multicolumn{3}{l}{$\{(0, 0),(1, 1),(3, x+2),(6, x),(4, x+1)\}.$}
\end{longtable}
(4) $-2,~x+2\in C_0^2.$
\begin{longtable}{lll}
$\{(0, 0),(1, x),(2, 1)\},$&
$\{(2, 0),(3, -1),(6, x)\},$&
$\{(1, 0),(2, x+2),(5, x+1)\}$,\\
$\{(1, 0),(3, 1),(7, -x)\}$,&
$\{(0, 0),(3, x+2),(5, 1)\},$&
$\{(0, 0),(6, x),(7, x+1)\};$\smallskip\\
\multicolumn{3}{l}{$\{(0, 0),(1, 1),(3, x+1),(6, x+2),(4, x)\}.$}
\end{longtable}

Case $2.$ $x+2=p.$
It yields that $-2\equiv x~(\bmod~p)$ and  $-2\in C_1^2.$ We consider two subcases.

(1) $x-2\in C_1^2.$
\begin{longtable}{lll}
$\{(0, 0),(1, 1-x),(2, -1)\},$&
$\{(0, 0 ),(3, x),(5, 1)\},$&
$\{(0, 0),(6, x-1),(7, 1)\},$\\
$\{(1, 0),(2, x-1),(5, x+1)\},$&
$\{(1, 0 ),(3, -1),(7, x)\},$&
$\{(2, 0),(3, x),(6, x+1)\};$\smallskip\\
\multicolumn{3}{l}{$\{(0, 0),(1, x-1),(3, 1),(6, x+1),(4, x)\}.$}
\end{longtable}
(2) $x-2\in C_0^2.$
\begin{longtable}{lll}
$\{(1, 0),(2, 1),(5, x)\},$
&
$\{(0, 0),(3, x), (5, -1)\},$&
$\{(0, 0),(6, 1 ),(7, x-1)\},$\\
$\{(1, 0),(3, x ),(7, 1)\},$&
$\{(0, 0),(1, x),(2, x-1)\},$
&
$\{(2, 0),(3, x),(6, x-1)\};$\smallskip\\
\multicolumn{3}{l}{$\{(0, 0),(1, 1),(3, x-1),(6, x),(4, x+1)\}.$}
\end{longtable}\qed

Semi-cyclic HGDDs are closely related to {\it cyclic holely difference matrices} (CHDMs). A $(k,wt;w)$-CHDM is a $k\times w(t-1)$ matrix $D=(d_{ij})$ with entries from $Z_{wt}$ such that for any two distinct rows $x$ and $y,$ the difference list $L_{xy}=\{d_{xj}-d_{yj}:j\in I_{w(t-1)}\}$ contains each integer of $Z_{wt}{\setminus}S$ exactly once, where $S=\{0,t,\ldots,(w-1)t\}$ is a subgroup of order $w$ in $Z_{wt}$. A $(k, t)$-CDM ({\it cyclic difference matrix}) is a $k\times t$ matrix $D=(d_{ij})$ with entries from $Z_t$ such that for any two distinct rows $x$ and $y,$ the difference list $\{d_{xj}-d_{yj}:j\in I_t\}$ contains each integer of $Z_t$ exactly once. Note that a $(k,t;1)$-CHDM is equivalent to a $(k,t)$-CDM. 
\begin{lemma}{\rm[\cite{WY,Y}]}\label{CDM}
A $k$-SCGDD of type $m^k$ is equivalent to a  $(k,m)$-CDM. And a $k$-SCHGDD of type $(k,w^t)$ is equivalent to a $(k,wt;w)$-CHDM.
\end{lemma}

\begin{lemma}\rm{[\cite{AG}, Theorem 4.2]}\label{g5}
A $(5,g)$-CDM exists for any integer $g$ which satisfies any of the following conditions: (1)~$\gcd (g,6)$=1, (2) $g\in\{15,27,39,51\}$.
\end{lemma}

\begin{lemma}\label{col3}
Let $v\equiv 1,5~(\bmod~6)$ and $v=3^i$ for  $i\geq 3$. Then there exists a $(5,v;1)$-CHDM.
\end{lemma}
\proof
  When $v\equiv 1,5~(\bmod~6)$, there exists a $(5,v)$-CDM by  [\cite[Theorem 2.1]{CC}]. When $v=3^i$ and $i\geq3$, by [\cite[Lemma 3.10]{G}] there exists a $(5,3^i)$-CDM. Thus the conclusion follows.
\qed

The following constructions are slight generalizations of several constructions in~[\cite{fwc}]. Here we omit their proofs and refer the reader to parallel proofs in [\cite{fwc}].  A $K$-GDP of type $w^u$ is said to be {\it cyclic}, if it admits an automorphism forming a cycle of length $wu$. If the length of each orbit in the cyclic $K$-GDP of type $w^u$ is $wu$, then the $K$-GDP is called {\it strictly cyclic}.

\begin{construction}\label{c6}{\rm[\cite{fwc}, Construction 3.1]}
If there exist a $K$-SCHGDD of type $(n,(gw)^t)$ and a $K$-SCHGDD of type $(n,g^w)$, then there exists a $K$-SCHGDD of type $(n,g^{wt})$.
\end{construction}
\begin{construction}\label{c5}{\rm[\cite{fwc}, Construction 3.4]}
If there exist a $K$-SCHGDD of type $(n,w^t)$ and a $(k,v)$-CDM for each $k\in K,$ then there exists a $K$-SCHGDD of type $(n,(wv)^t)$.
\end{construction}
\begin{construction}\label{pbd}{\rm[\cite{fwc}, Construction 3.2]}
Suppose that there exist a $K$-SCGDD of type $g^n$ and an $L$-SCHGDD
of type $(k,w^t)$ for each $k\in K.$ Then there exists an $L$-SCHGDD of type
$(n,(gw)^t).$
\begin{construction}\label{c7}{\rm[\cite{fwc}, Construction 3.3]}
If there exist a strictly cyclic $k$-GDD of type $w^v$ and a $K$-MGDD of type $k^u$, then there exists a $K$-SCHGDD of type $(u,w^v)$.
\end{construction}
\end{construction}
\begin{theorem}\label{5.1}{\rm[\cite{A1}]}
For any $v\equiv 2~(\bmod~6)$ and $v\geq 14$, there exists a PBD $(v,\{3,4,5\})$.
\end{theorem}

\begin{lemma}\label{3581w}
Let $w\equiv 1,5~(\bmod~6)$ and $w\geq5.$
Then there exists a $\{3,5\}$-SCHGDD of type $(v,1^w)$ for $v\in\{8,14,20\}$.
\end{lemma}
\proof
 First we let $v=8$. Let $w=p_1^{a_{1}}p_2^{a_{2}}\cdots p_{t}^{a_{t}}$ be its prime factorization, where $a_{i}\geq1$ and $p_i\geq 5$ is a prime, $1\leq i\leq t.$ Start from a $\{3,5\}$-SCHGDD of type $(8,1^{p_{1}})$ which exists from Lemma~\ref{81p}. Apply Construction~\ref{c5} with  $(k,q)$-CDMs for $k=3,5$ and $q\in\{p_2,\ldots, p_{t}\}$ which exist from Lemmas~\ref{scgdd}, \ref{CDM} and \ref{g5}. Then  we obtain a $\{3,5\}$-SCHGDD of type $(8,q^{p_1})$.
 Apply Construction~\ref{c6} with a $\{3,5\}$-SCHGDD of type $(8,1^q)$ which exists by Lemma~\ref{81p} to obtain a $\{3,5\}$-SCHGDD of type $(8,1^{p_1q})$. Repeating this process will produce a $\{3,5\}$-SCHGDD of type $(8,1^w)$.

  Then let $v\in\{14,20\}$.  By Theorem~\ref{5.1} there exists a PBD$(v,\{3,4,5\})$ which is also a $\{3,4,5\}$-GDD of type $1^v$. When $w\equiv 1,5~(\bmod~6)$ and $w\geq 5$, by Lemma~\ref{col3} there exists a $(5,w;1)$-CHDM, which is a 5-SCHGDD of type $(5,1^w)$ from Lemma~\ref{CDM}. A 3-SCHGDD of type $(k,1^w)$ for $k\in\{3,4\}$ exists by Theorem \ref{1}. Then apply Construction~\ref{pbd} to obtain a $\{3,5\}$-SCHGDD of type $(v,1^w)$.
\qed
\begin{example}\label{519}
We construct a $\{3,5\}$-SCHGDD on $X=I_5\times Z_9$ with group set $\mathcal{G}=\{\{i\}\times Z_9:i\in I_5\}$ and hole set $\mathcal{H}=\{I_5\times\{j\}:j\in Z_9\}$. We have two base blocks of size $5$: $\{(0,0),(1,1),(2,2),(3,3),(4,4)\}$,
$\{(0,4),(1,8),(2,7),(3,6),(4,5)\}$.
Then we list $20$ base blocks of size $3$ as follows.
\begin{longtable}{llll}
$\{(0,0),(2,1),(4,2)\}$,&
$\{(0,0),(3,1),(2,8)\}$,&
$\{(0,0),(1,2),(2,5)\}$,&
$\{(0,0),(1,3),(3,7)\}$,\\
$\{(0,0),(4,3),(1,7)\}$,&
$\{(0,0),(2,4),(1,8)\}$,&
$\{(0,0),(3,4),(4,8)\}$,&
$\{(0,0),(1,5),(3,6)\}$,\\
$\{(0,0),(3,5),(2,7)\}$,&
$\{(0,0),(4,5),(2,6)\}$,&
$\{(0,0),(1,6),(4,7)\}$,&
$\{(0,0),(4,6),(3,8)\}$,\\
$\{(1,0),(2,2),(3,6)\}$,&
$\{(1,0),(4,2),(2,6)\}$,&
$\{(1,0),(3,3),(2,7)\}$,&
$\{(1,0),(2,4),(4,8)\}$,\\
$\{(1,0),(4,4),(3,8)\}$,&
$\{(1,0),(3,5),(4,7)\}$,&
$\{(2,0),(3,3),(4,6)\}$,&
$\{(2,0),(4,3),(3,6)\}$.
\end{longtable}\qed
\end{example}

\begin{lemma}\label{83r}
Let $w=3^r$ with $r\geq2$. Then there exists a $\{3,5\}$-SCHGDD of type $(v,1^w)$ for $v\in\{8,14,20\}$.
\end{lemma}
\proof
(1) Let $v=8$.
When $r\in\{2,3\}$, 
the required~$\{3,5\}$-SCHGDDs of type $(8,1^w)$ are constructed directly
in Lemma~A.2 in Supporting Information.

When $r=4,$ by  [\cite[Lemma~2.5]{A}] there exists a strictly cyclic 5-GDD of type $1^{81}.$
By Lemma~\ref{81p} there exists a $\{3,5\}$-SCHGDD of type $(8,1^5)$ which is also a $\{3,5\}$-MGDD of type $5^8$. Then apply Construction~\ref{c7} to obtain a $\{3,5\}$-SCHGDD of type $(8,1^{81})$.

 When $r\geq 5$, we use induction on $r$. Assume that there exists a $\{3,5\}$-SCHGDD of type $(8,1^{3^{r-3}})$. By Lemma~\ref{scgdd} there exists a 3-SCGDD of type $27^3$, which is also a $(3,27)$-CDM from Lemma~\ref{CDM}. By Lemma~\ref{g5}, there exists a $(5,27)$-CDM. Then apply Construction~\ref{c5} with the given $\{3,5\}$-SCHGDD of type $(8,1^{3^{r-3}})$ to obtain a $\{3,5\}$-SCHGDD of type $(8,27^{3^{r-3}})$. Then apply Construction~\ref{c6} with a $\{3,5\}$-SCHGDD of type $(8,1^{27})$ to get a $\{3,5\}$-SCHGDD of type $(8,1^{3^r})$.

(2) Let $v\in\{14,20\}$.
By Theorem~\ref{5.1} there exists a PBD$(v,\{3,4,5\})$ which is a $\{3,4,5\}$-GDD of type $1^v$.
For $k\in\{3,4\}$, there exists a 3-SCHGDD of type $(k,1^{3^r})$ from
Theorem~\ref{1}. Apply Construction~\ref{pbd} with $\{3,5\}$-SCHGDD of type $(5,1^{3^r})$ to obtain a $\{3,5\}$-SCHGDD of type $(v,1^{3^{r}})$, where the needed $\{3,5\}$-SCHGDDs of type $(5,1^{3^r})$ exist from Example~\ref{519} and Lemma~\ref{col3}.\qed

\begin{lemma}\label{sc}
Let $v\in\{8,14,20\}$. For $w=3^r$ with $r\geq2$ and $w\equiv1,5~(\bmod~6)$ with $w\geq5$, there exists a $\{3,5,v\}$-SCGDD of type $w^v$ containing one base block of size $v$.
\end{lemma}
\proof
The conclusion follows from Lemmas ~\ref{eq},~\ref{3581w} and~\ref{83r}.
\qed

\begin{lemma}\label{38}
There exists a $\{3,4\}$-SCGDD of type $3^8$ in which the blocks of size $4$ form a parallel class.
\end{lemma}
\proof
The design is constructed on $X= I_{8}\times Z_{3}$ with group set  $\{\{i\}\times Z_{3}:i\in I_{8}\}$.  We have two base blocks of size $4$: $\{( 0, 0 ), ( 1, 0 ), ( 2, 0 ), ( 3, 0 )\}$,~$\{( 4, 0 ), ( 5, 0 ), ( 6, 0 ), ( 7, 0 )\}$.
Next we list $24$ base blocks of size $3$  as follows.

\begin{longtable}{llllllllll}
$\{( 0, 0 ), ( 4, 0 ), ( 1, 1 )\}$, & $\{( 0, 0 ), ( 6, 0 ), ( 5, 1 )\}$,  & $\{( 0, 0 ), ( 5, 0 ), ( 7, 2 )\}$,&
  $\{( 0, 0 ), ( 7, 0 ), ( 4, 1 )\}$,\\  $\{( 0, 0 ), ( 2, 1 ), ( 5, 2 )\}$, & $\{( 0, 0 ), ( 6, 1 ), ( 3, 2 )\}$,&
  $\{( 0, 0 ), ( 1, 2 ), ( 4, 2 )\}$, & $\{( 0, 0 ), ( 3, 1 ), ( 6, 2 )\}$, \\ $\{( 0, 0 ), ( 7, 1 ), ( 2, 2 )\}$,&
  $\{( 3, 0 ), ( 7, 0 ), ( 4, 2 )\}$, & $\{( 2, 0 ), ( 4, 1 ), ( 6, 2 )\}$,&  $\{( 2, 0 ), ( 4, 0 ), ( 5, 2 )\}$,\\
  $\{( 2, 0 ), ( 6, 0 ), ( 7, 1 )\}$,&  $\{( 1, 0 ), ( 5, 0 ), ( 6, 1 )\}$, & $\{( 1, 0 ), ( 7, 1 ), ( 3, 2 )\}$,&
  $\{( 1, 0 ), ( 2, 1 ), ( 6, 2 )\}$,\\  $\{( 1, 0 ), ( 6, 0 ), ( 7, 2 )\}$,  &$\{( 3, 0 ), ( 6, 0 ), ( 4, 1 )\}$,&
  $\{( 2, 0 ), ( 5, 0 ), ( 3, 1 )\}$,  &$\{( 1, 0 ), ( 7, 0 ), ( 5, 2 )\}$, \\ $\{( 3, 0 ), ( 4, 0 ), ( 5, 1 )\}$,&
  $\{( 1, 0 ), ( 4, 1 ), ( 2, 2 )\}$, & $\{( 2, 0 ), ( 7, 0 ), ( 3, 2 )\}$,  &$\{( 1, 0 ), ( 3, 1 ), ( 5, 1 )\}$.
\end{longtable}\qed

\section{Case $u\equiv1~(\bmod~3)$ and proof of the main result}

In this section we  consider the maximum $w$-cyclic $3$-HGDPs of type $(u,w^v)$ when $u\equiv 1~(\bmod~3)$. Then we prove our main result in Theorem~\ref{all}.

\begin{example}\label{e1}
We  construct a $3$-HGDP of type $(4,1^4)$ with $J(4\times 4\times 1,3,1)$ blocks on $I_{4}\times I_{4}$ with  group set $\{\{i\}\times I_{v}:i\in I_{4}\}$ and hole set $\{I_{4}\times \{j\}:j\in I_{4}\}.$ All $21$ blocks are listed as follows.
\begin{longtable}{lllll}
$\{(0,0),(1,1),(3,2)\}$,&
$\{(0,0),(2,1),(1,2)\}$,&
$\{(0,0),(3,1),(2,3)\}$,&
$\{(0,0),(2,2),(1,3)\}$,\\
$\{(1,0),(3,1),(0,3)\}$,&
$\{(1,0),(0,2),(2,3)\}$,&
$\{(1,0),(2,2),(3,3)\}$,&
$\{(2,0),(0,1),(1,3)\}$,\\
$\{(2,0),(1,1),(3,3)\}$,&
$\{(2,0),(3,1),(1,2)\}$,&
$\{(2,0),(3,2),(0,3)\}$,&
$\{(3,0),(0,1),(2,2)\}$,\\
$\{(3,0),(1,1),(0,2)\}$,&
$\{(3,0),(2,1),(0,3)\}$,&
$\{(3,0),(1,2),(2,3)\}$,&
$\{(0,1),(1,2),(3,3)\}$,\\
$\{(0,1),(3,2),(2,3)\}$,&
$\{(1,1),(2,2),(0,3)\}$,&
$\{(2,1),(0,2),(3,3)\}$,&
$\{(2,1),(3,2),(1,3)\}$,\\
$\{(3,1),(0,2),(1,3)\}$.& &&
\end{longtable}\qed
\end{example}

\begin{lemma}\label{z6}
Let $u\in\{4,10,16\}$ and $v\in\{4,8,10,14,16,20\}$.~There exists a $3$-HGDP of type $(u,1^v)$ with $J(u\times v\times1,3,1)$ blocks.
\end{lemma}
\proof
 For $u=4$ and $v=4,8,14$, the desired designs exist from Example~\ref{e1} and Lemmas B.6, B.7 in Supporting Information.

For other parameters, apply Lemma~\ref{z5} by taking parameters $(u,t,v)$ in Table II to yield a $3$-HGDP of type $(u,1^v)$ with $J(u\times v\times1,3,1)$ blocks, where the needed 3-IHGDPs of type $(u,t,1^v)$ and 3-HGDPs of type $(t,1^v)$ exist from the sources listed in the fourth and fifth column, respectively. Again note that  a $3$-HGDP of type $(u,1^v)$ is equivalent to a $3$-HGDP of type $(v,1^u)$.

\begin{longtable}{cccccc}
\caption*{Table II}\\
\hline
Row&$(u,v)$&$(u,t,v)$& $(u,t,1^v)$& $(t,1^v)$&\\
\hline
1&$(10,4)$& $(10,4,4)$ & Lemma C.1 & Example~\ref{e1}\\
2&$(16,4)$& $(16,4,4)$& Lemma C.1&Example~\ref{e1}\\
3&$(20,4)$& $(20,8,4)$& Lemma~\ref{sc2}&Lemma~B.6\\
4&$(10,8)$&$(10,4,8)$& Lemma C.3 &Lemma~B.6\\
5&$(16,8)$&$(16,4,8)$& Lemma C.4&Lemma~B.6\\
6&$(10,10)$& $(10,4,10)$& Lemma C.2&Row 1\\
7&$(10,14)$&$(10,4,14)$& Lemma C.3&Lemma~B.7\\
8&$(16,14)$& $(16,4,14)$& Lemma C.4&Lemma~B.7\\
9&$(10,16)$ &$(10,4,16)$& Lemma~\ref{z4}&Row 2\\
10&$(16,16)$&$(16,4,16)$& Lemma~\ref{z4}&Row 2\\
11&$(20,16)$&$(20,8,16)$& Lemma~\ref{z4}&Row 5\\
12&$(10,20)$&$(10,4,20)$& Lemma~\ref{z4}& Row 3\\
\hline
\end{longtable}
\qed
\begin{lemma}\label{391420}
Let  $u\in\{4,10,16\}$ and $v\in\{8,14,20\}.$ There exists a $3$-cyclic $3$-HGDP of type $(u,3^v)$ with $J(u\times v\times 3,3,1)$ base blocks.
\end{lemma}
\proof
For $v\in\{8,14,20\}$, by Lemma~\ref{38} and Lemma D.1 in Supporting Information, there exists a $\{3,4,6\}$-SCGDD of type $3^v$ with $\frac{3v^2-8v+16}{6}$ base blocks of size $3$, two base blocks of size $4$ and $\frac{v-8}{6}$ base blocks of size $6$. Apply Corollary~\ref{col4} by taking $h=1$ and $w=3$ to yield a desired $3$-cyclic $3$-HGDP, where the required $3$-HGDPs of types $(u,1^3)$, $(u,1^4)$ and $(u,1^6)$ exist by Theorem~\ref{hgdd}, Lemmas~\ref{z6} and \ref{3.8}, respectively. The number of base blocks equals
$$
\frac{3v^2-8v+16}{6}\times u(u-1)+2\times\frac{6u^2-8u-1}{3}+\frac{v-8}{6}\times u(5u-6)$$
$$=
\frac{3uv(u-1)(v-1)-uv-4}{6}=J(u\times v\times3,3,1).
$$
\qed


\begin{lemma}\label{4.4}
Let $u\in\{4,10,16\}$,  $v\in\{4,8,10,14,16,20\}$ and $w\equiv1~(\bmod~2)$. Then there exists a $w$-cyclic $3$-HGDP of type $(u,w^v)$ with $J(u\times v\times w,3,1)$ base blocks.
\end{lemma}
\proof
(a) For $w=1$ and the assumed $u,v$, by Lemma~\ref{z6}, there exists a $3$-HGDP of type $(u,1^v)$ with $J(u\times v\times 1,3,1)$ blocks.

(b) For $w=3$ and the assumed $u$,
when $v\in\{8,14,20\}$, a $3$-cyclic $3$-HGDP of type $(u,3^v)$ with $J(u\times v\times 3,3,1)$ base blocks exists from~Lemma~\ref{391420}.

(c) For other assumed parameters,
apply Lemma~\ref{z7} by taking parameters $(g,v,u,m)$ in the second column of Table III.  The required $\{3,5,v\}$-SCGDDs of type $g^v$ containing one base block of size $v$ exist, for which the sources are listed in the fourth column.  Maximum $m$-cyclic 3-HGDPs of type $(u,m^v)$ exist from the sources  in the last column.  Finally we produce $w$-cyclic 3-HGDPs of type $(u,w^v)$ with $J(u\times v\times w, 3,1)$ base blocks.
\begin{longtable}{ccccc}
\caption*{Table III}\\
\hline
Row & $(u,w^v)$ & $(g,v,u,m)$&SCGDD & $(u,m^v)$\\
\hline
1&$v\in\{4,10,16\},w\equiv1(\bmod~2),w\geq 3$ & $(w,v,u,1)$ &Theorem~\ref{1} & Lemma~\ref{z6}\\
2&$v\in\{8,14,20\},w\equiv1,5~(\bmod~6),w \geq 5$&$(w,v,u,1)$& Lemma~\ref{sc} & Lemma~\ref{z6}\\
3&$v\in\{8,14,20\},w=3^r,r\geq 2$ &$(w,v,u,1)$& Lemma~\ref{sc} & Lemma~\ref{z6}\\
4&$v\in\{8,14,20\},w=3^rw',w'\equiv1,5~(\bmod~6)$&$(w',v,u,3^r)$ & Lemma~\ref{sc}&(b), Row 3
 \\  & ($r\ge 1, w'\ge 5$) & & & \\
\hline
\end{longtable}

\qed
\begin{lemma}\label{48}
Let $m\geq3$, $w\equiv 1~(\bmod~2)$, $v\equiv 0~(\bmod~2)$ and $s\in\{4,8\}$.
There exists a $w$-cyclic $3$-DGDD of type $(\mathbf{d}_0,\mathbf{d}_1,\ldots,\mathbf{d}_{v-1})$
where  $\mathbf{d}_e=(6w,6w,\ldots,6w,sw)^T$ is
a column vector of length $m+1$ for $0\leq e\leq v-1$.
\end{lemma}
\proof
For the assumed parameters, apply Construction~\ref{xc2} with a 3-GDD of type $6^ms^1$ whose existence is guaranteed by Theorem~\ref{YL1} and a $w$-cyclic 3-HGDD of type $(3,w^v)$  from Theorem~\ref{Yt11}. It is not difficult to learn that we obtain a $w$-cyclic $3$-DGDD of type $(\mathbf{d}_0,\mathbf{d}_1,\ldots,\mathbf{d}_{v-1})$.
\qed

\begin{lemma}\label{4wv}
Let $u\in\{4,10,16\}$,~$v\equiv0~(\bmod~2)~(v\geq4)$ and $w\equiv1~(\bmod~2)$. Then
there exists a $w$-cyclic $3$-HGDP of type $(u,w^v)$ with $J(u\times v\times w,3,1)$ base blocks.
\end{lemma}
\proof
For $v\in\{4,8,10,14,16,20\},$ the required designs exist by Lemma~\ref{4.4}.

Next note that $u\in\{4,10,16\}$ and $w\equiv1~(\bmod~2)$. For $v\equiv 0~(\bmod~6)$ and $v\geq6$, by Lemma~\ref{3.8} there exists a $w$-cyclic $3$-HGDP of type $(v,w^u)$, also of type $(u,w^v)$, with $J(u\times v\times w,3,1)$ base blocks.

 For $v\geq22$, let $s=8$ if $v\equiv 2~(\bmod~6)$ and let $s=4$ if $v\equiv 4~(\bmod~6)$. By Lemma~\ref{48} there exists a $w$-cyclic $3$-DGDD of type $(\mathbf{d}_0,\mathbf{d}_1,\ldots,\mathbf{d}_{u-1})$, where
 $\mathbf{d}_e=(6w,6w,\ldots,6w,sw)^T$ is
a column vector of length $\frac{v-s+6}{6}$ for $0\leq e\leq u-1$. For $k\in\{6,s\}$, there exists a $w$-cyclic $3$-HGDP of type $(k,w^u)$ with $J(k\times u\times w,3,1)$ base blocks from Lemmas~\ref{t1} and~\ref{4.4}.
Apply Lemma~\ref{zz} by taking parameters $(g,s,v,u,w)$ in Table IV to obtain a desired $w$-cyclic $3$-HGDP of type $(v,w^u)$,  also of type $(u,w^v)$.
\begin{table}[!ht]\caption*{Table IV}
\centering
\begin{tabular}{cc}
\hline
~~~$(u,w^v)$& $(g,s,v,u,w)$\\
\hline
~~~$v\equiv2~(\bmod~6)$, $v\geq26$ & $(6,8,v,u,w)$\\
~~~$v\equiv4~(\bmod~6)$, $v\geq22$ &
$(6,4,v,u,w)$\\
\hline
\end{tabular}
\end{table}

\qed
\begin{lemma}\label{u4}
$\Psi(u\times v\times w,3,1)=J(u\times v\times w,3,1)$
where $u\equiv 1~(\bmod~3)$ and $u, v\geq3$.
\end{lemma}
\proof
By Corollary~\ref{co1}, we only need to treat the case that $u\equiv4~(\bmod~6)$, $w\equiv1~(\bmod~2)$ and $v\equiv0~(\bmod~2)$, $v\geq4$.

For $u\in\{4,10,16\}$, the required designs exist by Lemma~\ref{4wv}.

 For $u\geq 22$, by Lemma~\ref{48} there exists a $w$-cyclic $3$-DGDP of type $(\mathbf{d}_0,\mathbf{d}_1,\ldots,\mathbf{d}_{v-1})$, where
 $\mathbf{d}_e=(6w,6w,\ldots,6w,4w)^T$ is
a column vector of length $\frac{u+2}{6}$ for $0\leq e\leq v-1$.
 And there exist maximum $w$-cyclic 3-HGDPs of type $(4,w^v)$ and $(6,w^v)$ from Lemmas~\ref{4wv} and~\ref{3.8}, respectively. Apply Lemma~\ref{zz}  by taking parameters $(g,s,v,u,w)=(6,4,u,v,w)$ to obtain a $w$-cyclic $3$-HGDP of type $(u,w^v)$ with $J(u\times v\times w,3,1)$ base blocks.
 \qed

Now we  prove our main result.

\noindent
{\em{Proof of Theorem~$\ref{all}$.}}
 Combining  Lemmas~\ref{3.8} and~\ref{u4} completes the proof. \qed
\section{Application}
An {\it optical orthogonal code} (OOC) is a family of $(0,1)$ sequences with good correlation properties.
It was introduced  as signature sequences to facilitate multiple access in optical fibre networks (see, e.g. [\cite{CSW89,J,JC}]).~OOCs have been found applications in wide ranges such as mobile radio, frequency-hoping spread-spectrum communications, radar, sonar signal designs~[\cite{G64}],  collision channels without feedback~[\cite{MM85}], and neuromorphic networks~[\cite{VS88}].

Let $u,v,w,k$ and $\lambda$ be positive integers, where $u,v$ and $w$ are the number of spatial channels, wavelengths, and time slots, respectively. A {\it three-dimensional $(u\times v\times w,k,\lambda)$  optical orthogonal code} (briefly 3-D $(u\times v\times w,k,\lambda)$-OOC), $\mathcal{C}$, is a family of $u\times v\times w$ $(0,1)$ arrays (called {\it codewords}) of Hamming weight $k$ satisfying that for any two arrays $A=[a(i,j,l)],$ $B=[b(i,j,l)]\in\mathcal{C}$ and any integer $\tau,$
$$\sum_{i=0}^{u-1}\sum_{j=0}^{v-1}\sum_{l=0}^{w-1}a(i,j,l)b(i,j,l+\tau)\leq \lambda,$$
where either $A\neq B$ or $\tau\not\equiv0~(\bmod~w)$, and the arithmetic $l+\tau$ is reduced modulo $w.$

A $v/w$ plane is usually called a {\it spatial plane}, and a $u/w$ plane is called a {\it wavelength plane}. There are several classes of 3-D OOCs of particular interest, for which the following three additional restrictions on the placement of pulses are often placed within the arrays of a 3-D OOC to simplify practical implementation.
\begin{itemize}
  \item One-pulse per plane (OPP) restriction: each codeword contains exactly one optical pulse per spatial plane.

  \item  At most one-pulse per spatial plane (AM-OPPSP) restriction: for any array in $\mathcal{C}$, the element $1$ appears at most once
in each spatial plane.
  \item At most one-pulse per wavelength plane (AM-OPPWP) restriction: for any array in $\mathcal{C}$, the element $1$ appears at most
once in each wavelength plane.
\end{itemize}


A 3-D OOC with both AM-OPPSP and AM-OPPWP properties is simply denoted by AM-OPPS/WP 3-D OOC. The number of codewords of a $3$-D OOC is called its {\it size}.
 An AM-OPPS/WP 3-D $(u\times v\times w,k,1)$-OOC of the largest possible size is said to be maximum.
Dai et al.  established the equivalence between  a maximum  AM-OPPS/WP 3-D $(u\times v\times w,k,1)$-OOC and a maximum $w$-cyclic $k$-HGDP of type $(u,w^{v})$.

\begin{theorem}{\rm [\cite{DCW},~Theorem 5.1]}\label{HGDP}
A maximum AM-OPPS/WP $3$-D$~(u\times v\times w,k,1)$-OOC is equivalent to a maximum $w$-cyclic $k$-HGDP of  type $(u,w^{v})$.
\end{theorem}
By Theorems~\ref{all} and~\ref{HGDP}, the size of a maximum AM-OPPS/WP $3$-D $(u\times v\times w,3,1)$-OOC is determined for any integers $u,v,w$ with $u\equiv0,1~(\bmod~3)$.
\begin{theorem}
A maximum  AM-OPPS/WP $3$-D $(u\times v\times w,3,1)$-OOC contains $\Psi(u\times v\times w,3,1)$ codewords for positive integers $u,v,w$ with $u\equiv0,1~(\bmod~3)$ and $u,v\geq3$, where
$$\Psi(u\times v\times w,3,1)=\left\{
  \begin{array}{ll}
   6(w-1), & \hbox{if $(u,v)=(3,3)$ and $w\equiv 0~(\bmod~2)$,} \\
   J(u\times v\times w,3,1), & \hbox{otherwise.}
  \end{array}
\right.$$
\end{theorem}

\noindent{\bf Acknowledgements:} The authors would like to thank Prof. Yanxun Chang for many valuable suggestions and comments.

\bigskip\bigskip
{\centering\title{\bf\Large{ Appendix: Supporting Information}}}
\appendix
\renewcommand{\appendixname}{~\Alph{section}}

\section{$K$-SCHGDDs of type $(u,1^v)$}
\begin{lemma}
There exists a $\{3,6\}$-HGDD of type $(6,1^{v})$ for $v\in\{8,14\}$ in which the blocks of size $6$ form a parallel class.
\end{lemma}
\proof
For $v\in\{8,14\}$, the design will be constructed on $X=I_{6}\times Z_{v}$ with group set $\mathcal{G}=\{\{i\}\times Z_{v}:i\in I_{6}\}$ and hole set $\mathcal{H}=\{I_{6}\times \{j\}:j\in Z_{v}\}.$ All blocks are obtained by developing by $(-,+1)\bmod~(-,v)$ successively from the following $5(v-2)$ base blocks of size 3 and one base block of size 6.

$v=8:$ 

\begin{longtable}{llll}


$\{(0,0),(1,1),(3,2)\}$,&
$\{(0,0),(2,1),(5,6)\}$,&
$\{(0,0),(3,1),(2,6)\}$,&
$\{(0,0),(4,1),(5,7)\}$,
\\
$\{(0,0),(5,1),(4,2)\}$,&
$\{(0,0),(1,2),(2,3)\}$,&
$\{(0,0),(5,2),(1,5)\}$,&
$\{(0,0),(1,3),(5,4)\}$,
\\
$\{(0,0),(4,3),(2,5)\}$,&
$\{(0,0),(5,3),(3,6)\}$,&
$\{(0,0),(1,4),(4,5)\}$,&
$\{(0,0),(2,4),(1,7)\}$,
\\
$\{(0,0),(3,4),(2,7)\}$,&
$\{(0,0),(3,5),(4,7)\}$,&
$\{(0,0),(4,6),(3,7)\}$,&
$\{(1,0),(2,2),(3,4)\}$,
\\
$\{(1,0),(3,2),(4,7)\}$,&
$\{(1,0),(4,2),(5,4)\}$,&
$\{(1,0),(5,2),(3,6)\}$,&
$\{(1,0),(2,3),(4,4)\}$,
\\
$\{(1,0),(3,3),(5,6)\}$,&
$\{(1,0),(4,3),(3,7)\}$,&
$\{(1,0),(5,3),(2,7)\}$,&
$\{(1,0),(4,5),(2,6)\}$,
\\
$\{(2,0),(5,1),(4,4)\}$,&
$\{(2,0),(5,2),(3,4)\}$,&
$\{(2,0),(4,3),(5,6)\}$,&
$\{(2,0),(4,5),(3,7)\}$,
\\
$\{(2,0),(3,6),(5,7)\}$,&$\{(3,0),(4,3),(5,7)\};$\smallskip\\
\multicolumn{2}{l}{$\{(0,0),(2,2),(3,3),(4,4),(5,5),(1,6)\}$.}

\end{longtable}

 $v=14:$
\begin{longtable}{llll}
$\{(0,0),(1,1),(3,2)\}$,
&$\{(0,0),(2,1),(5,6)\}$,
&$\{(0,0),(3,1),(5,11)\}$,
&$\{(0,0),(4,1),(5,5)\}$,\\
$\{(0,0),(2,12),(4,13)\}$,&
$\{(0,0),(1,2),(2,3)\}$,&
$\{(0,0),(4,2),(1,3)\}$,&
$\{(0,0),(5,2),(1,9)\}$,\\
$\{(0,0),(4,3),(2,7)\}$,&
$\{(0,0),(5,3),(1,4)\}$,&
$\{(0,0),(2,4),(1,13)\}$,&
$\{(0,0),(3,4),(1,5)\}$,
\\
$\{(0,0),(4,11),(2,13)\}$,&
$\{(0,0),(5,1),(4,6)\}$,&
$\{(0,0),(3,6),(1,12)\}$,&
$\{(0,0),(3,5),(2,6)\}$,
\\
$\{(0,0),(1,7),(2,9)\}$,&
$\{(0,0),(3,7),(5,8)\}$,&
$\{(0,0),(5,7),(1,10)\}$,&
$\{(0,0),(4,7),(2,8)\}$,
\\
$\{(0,0),(3,8),(2,10)\}$,&
$\{(0,0),(4,8),(3,12)\}$,&
$\{(0,0),(3,9),(1,11)\}$,&
$\{(1,0),(2,3),(3,7)\}$,
\\
$\{(0,0),(5,9),(3,10)\}$,&
$\{(0,0),(4,10),(5,12)\}$,&
$\{(2,0),(3,3),(4,5)\}$,&
$\{(2,0),(5,2),(4,9)\}$,
\\
$\{(2,0),(5,4),(4,7)\}$,&
$\{(0,0),(4,12),(5,13)\}$,&
$\{(1,0),(4,1),(2,11)\}$,&
$\{(1,0),(5,1),(4,5)\}$,
\\
$\{(1,0),(3,2),(4,10)\}$,&
$\{(1,0),(4,2),(3,3)\}$,&
$\{(0,0),(4,9),(3,11)\}$,&
$\{(1,0),(5,2),(2,4)\}$,
\\
$\{(1,0),(4,3),(5,8)\}$,&
$\{(1,0),(5,3),(4,4)\}$,&
$\{(1,0),(3,5),(2,12)\}$,&
$\{(1,0),(3,4),(4,8)\}$,
\\
$\{(1,0),(5,5),(3,10)\}$,&
$\{(1,0),(2,6),(5,12)\}$,&
$\{(1,0),(3,6),(4,9)\}$,&
$\{(1,0),(4,6),(2,9)\}$,
\\
$\{(1,0),(5,6),(3,9)\}$,&
$\{(1,0),(2,7),(5,10)\}$,&
$\{(1,0),(4,7),(2,13)\}$,&
$\{(1,0),(2,8),(5,9)\}$,
\\
$\{(0,0),(2,11),(3,13)\}$,&
$\{(2,0),(4,3),(3,10)\}$,&
$\{(0,0),(4,5),(1,8)\}$,&
$\{(0,0),(5,4),(2,5)\}$,
\\
$\{(2,0),(3,5),(5,11)\}$,&
$\{(2,0),(3,6),(5,9)\}$,&
$\{(2,0),(4,6),(3,11)\}$,&
$\{(2,0),(5,7),(3,9)\}$,
\\
$\{(2,0),(3,8),(5,10)\}$,&
$\{(3,0),(5,4),(4,6)\}$,&
$\{(3,0),(5,5),(4,11)\}$,&
$\{(3,0),(4,5),(5,8)\};$\smallskip\\
\multicolumn{2}{l}{$\{(0,0),(2,2),(3,3),(4,4),(1,6),(5,10)\}.$}
\end{longtable}
\qed
\begin{lemma}
There exists a $\{3,5\}$-SCHGDD of type $(8,1^w)$ for $w\in\{9,27\}$.
\end{lemma}
\proof
When $w=9,$ we construct a $\{3,5\}$-SCHGDD of type $(8,1^9)$ on $I_8\times Z_9$ with  group set $\{\{i\}\times Z_{9}:i\in I_{8}\}$ and hole set $\{I_{8}\times \{j\}:j\in Z_{9}\}.$ We list all $68$ base blocks of size $3$  and two base blocks of size $5$ as follows.
\begin{longtable}{llll}
$\{(0,0),(1,1),(3,2)\}$,&
$\{(0,0),(2,1),(1,5)\}$,&
$\{(0,0),(3,1),(4,4)\}$,&
$\{(0,0),(5,1),(6,8)\}$,\\
$\{(0,0),(6,1),(1,2)\}$,&
$\{(0,0),(7,1),(4,2)\}$,&
$\{(0,0),(5,2),(6,3)\}$,&
$\{(0,0),(6,2),(7,3)\}$,\\
$\{(0,0),(7,2),(6,6)\}$,&
$\{(0,0),(1,3),(2,4)\}$,&
$\{(0,0),(2,3),(3,7)\}$,&
$\{(0,0),(4,3),(1,4)\}$,\\
$\{(0,0),(5,3),(1,8)\}$,&
$\{(0,0),(3,4),(2,5)\}$,&
$\{(0,0),(6,4),(3,6)\}$,&
$\{(0,0),(7,4),(3,5)\}$,\\
$\{(0,0),(4,5),(5,7)\}$,&
$\{(0,0),(5,5),(1,6)\}$,&
$\{(0,0),(6,5),(7,7)\}$,&
$\{(0,0),(7,5),(2,6)\}$,\\
$\{(0,0),(4,6),(6,7)\}$,&
$\{(0,0),(5,6),(2,8)\}$,&
$\{(0,0),(7,6),(3,8)\}$,&
$\{(0,0),(1,7),(4,8)\}$,\\
$\{(0,0),(2,7),(7,8)\}$,&
$\{(0,0),(4,7),(5,8)\}$,&
$\{(1,0),(5,1),(6,7)\}$,&
$\{(1,0),(6,1),(2,3)\}$,\\
$\{(1,0),(7,1),(2,6)\}$,&
$\{(1,0),(2,2),(6,6)\}$,&
$\{(1,0),(3,2),(4,3)\}$,&
$\{(1,0),(4,2),(3,3)\}$,\\
$\{(1,0),(5,2),(7,8)\}$,&
$\{(1,0),(7,2),(3,5)\}$,&
$\{(1,0),(6,3),(3,4)\}$,&
$\{(1,0),(7,3),(6,4)\}$,\\
$\{(1,0),(2,4),(4,5)\}$,&
$\{(1,0),(7,4),(4,7)\}$,&
$\{(1,0),(5,5),(3,7)\}$,&
$\{(1,0),(6,5),(3,8)\}$,\\
$\{(1,0),(7,5),(5,7)\}$,&
$\{(1,0),(3,6),(2,8)\}$,&
$\{(1,0),(4,6),(7,7)\}$,&
$\{(1,0),(5,6),(2,7)\}$,\\
$\{(2,0),(5,1),(4,5)\}$,&
$\{(2,0),(6,1),(3,5)\}$,&
$\{(2,0),(3,2),(5,5)\}$,&
$\{(2,0),(4,2),(5,6)\}$,\\
$\{(2,0),(6,2),(4,6)\}$,&
$\{(2,0),(7,2),(4,7)\}$,&
$\{(2,0),(3,3),(7,5)\}$,&
$\{(2,0),(4,3),(3,6)\}$,\\
$\{(2,0),(5,3),(7,7)\}$,&
$\{(2,0),(6,3),(7,6)\}$,&
$\{(2,0),(7,3),(6,5)\}$,&
$\{(2,0),(4,4),(6,8)\}$,\\
$\{(2,0),(5,4),(6,6)\}$,&
$\{(3,0),(6,1),(5,6)\}$,&
$\{(3,0),(7,1),(5,5)\}$,&
$\{(3,0),(4,2),(6,4)\}$,\\
$\{(3,0),(5,2),(7,4)\}$,&
$\{(3,0),(6,2),(5,8)\}$,&
$\{(3,0),(6,3),(4,4)\}$,&
$\{(3,0),(7,3),(4,5)\}$,\\
$\{(3,0),(5,4),(7,5)\}$,&
$\{(4,0),(6,3),(5,7)\}$,&
$\{(4,0),(7,3),(6,6)\}$,&
$\{(4,0),(7,5),(5,6)\};$\smallskip\\
\multicolumn{4}{l}{$\{(0,0),(4,1),(2,2),(3,3),(5,4)\},~~\{(7,1),(1,4),(6,6),(5,7),(4,8)\}$.}
\end{longtable}
When $w=27,$  we construct a $\{3,5\}$-SCHGDD of type $(8,1^{27})$ on $I\times Z_{27}$, where $I=Z_{7}\cup\{\infty\}$, with  group set $\{\{i\}\times Z_{27}:i\in I\}$ and  hole set $\{I\times \{j\}:j\in Z_{27}\}.$ All base blocks can be obtained by developing by $(+1,-)\bmod (7,-)$ successively from the following $28$ initial base blocks of size $3$  and two initial base blocks of $5$. Note that $\infty+1=\infty.$
\small{\begin{longtable}{llll}
$\{(0,0),(4,12),(\infty,26)\}$,&
$\{(0,0),(5,1),(6,6)\}$,&
$\{(0,0),(4,2),(\infty,15)\}$,&
$\{(0,0),(5,2),(3,11)\}$,\\
$\{(0,0),(1,3),(5,21)\}$,&
$\{(0,0),(6,3),(\infty,25)\}$,&$\{(0,0),(\infty,8),(4,15)\}$,
&
$\{(0,0),(5,3),(6,23)\}$,\\
$\{(0,0),(2,4),(4,11)\}$,&
$\{(0,0),(4,4),(5,15)\}$,&
$\{(0,0),(\infty,4),(3,13)\}$,&
$\{(0,0),(3,5),(2,10)\}$,\\
$\{(0,0),(4,5),(1,13)\}$,&
$\{(0,0),(5,5),(2,15)\}$,&
$\{(0,0),(\infty,5),(6,21)\}$,&
$\{(0,0),(3,6),(1,19)\}$,\\
$\{(0,0),(4,6),(2,17)\}$,&
$\{(0,0),(\infty,6),(1,17)\}$,&
$\{(0,0),(1,7),(5,14)\}$,&
$\{(0,0),(3,7),(5,18)\}$,\\
$\{(0,0),(5,7),(1,15)\}$,&
$\{(0,0),(\infty,7),(6,15)\}$,&
$\{(0,0),(1,8),(\infty,17)\}$,&
$\{(0,0),(5,8),(4,17)\}$,\\
$\{(0,0),(\infty,10),(6,13)\}$,&
$\{(0,0),(4,1),(\infty,2)\}$,&$\{(0,0),(6,11),(\infty,23)\}$,
&$\{(0,0),(\infty,3),(4,9)\};$\smallskip\\
\multicolumn{4}{l}{$\{(0,0),(1,1),(3,2),(2,3),(5,4)\}$, $\{(0,0),(3,3),(2,5),(1,9),(4,13)\}$.}
\end{longtable}}
\qed

\section{$w$-cyclic 3-HGDPs of type $(u,w^v)$}
\begin{lemma}
There exists a $3$-HGDP of type $(6,1^{v})$ with $J(6\times v\times1,3,1)$ blocks for $v\in\{6,10\}$.
\end{lemma}
\proof
For $v\in\{6,10\}$, the design will be constructed on $X=I_{6}\times Z_{v}$ with group set $\mathcal{G}=\{\{i\}\times Z_{v}:i\in I_{6}\}$ and hole set $\mathcal{H}=\{I_{6}\times \{j\}:j\in Z_{v}\}$. We list all $5v-6$ base blocks as follows.

$v=6:$

\begin{minipage}[t]{0.45\textwidth}
\centering
\begin{tabular}{lccc}
$\{(0,0),(2,1),(1,2)\}$,&
$\{(0,0),(3,1),(2,2)\}$,&
$\{(0,0),(4,1),(3,2)\}$,&
$\{(0,0),(5,1),(4,2)\}$,
\\
$\{(0,0),(5,2),(1,3)\}$,&
$\{(0,0),(2,3),(3,4)\}$,&
$\{(0,0),(3,3),(1,4)\}$,&
$\{(0,0),(4,3),(1,5)\}$,
\\
$\{(0,0),(5,3),(2,5)\}$,&
$\{(0,0),(2,4),(4,5)\}$,&
$\{(0,0),(4,4),(5,5)\}$,&
$\{(0,0),(5,4),(3,5)\}$,
\\
$\{(1,0),(2,1),(5,3)\}$,&
$\{(1,0),(3,1),(4,3)\}$,&
$\{(1,0),(4,1),(3,3)\}$,&
$\{(1,0),(5,1),(2,4)\}$,
\\
$\{(1,0),(2,2),(4,5)\}$,&
$\{(1,0),(3,2),(5,4)\}$,&
$\{(1,0),(4,2),(2,3)\}$,&
$\{(1,0),(5,2),(3,4)\}$,
\\
$\{(2,0),(5,1),(3,4)\}$,&
$\{(2,0),(4,2),(5,5)\}$,&
$\{(2,0),(3,3),(4,4)\}$,&
$\{(3,0),(5,1),(4,3)\}$.\\
\end{tabular}
\end{minipage}

\makeatletter\def\@captype{table}\makeatother
$v=10:$

\begin{longtable}{lccc}
$\{(0,0),(2,1),(1,2)\}$,&
$\{(0,0),(3,1),(5,6)\}$,&
$\{(0,0),(4,1),(5,7)\}$,&
$\{(0,0),(5,1),(1,6)\}$,
\\
$\{(0,0),(2,2),(4,6)\}$,&
$\{(0,0),(3,2),(5,4)\}$,&
$\{(0,0),(4,2),(1,3)\}$,&
$\{(0,0),(5,2),(2,3)\}$,
\\
$\{(0,0),(3,3),(1,4)\}$,&
$\{(0,0),(4,3),(2,4)\}$,&
$\{(0,0),(5,3),(2,7)\}$,&
$\{(0,0),(3,4),(2,5)\}$,
\\
$\{(0,0),(4,4),(3,5)\}$,&
$\{(0,0),(1,5),(2,6)\}$,&
$\{(0,0),(4,5),(1,7)\}$,&
$\{(0,0),(5,5),(3,6)\}$,
\\
$\{(0,0),(3,7),(4,8)\}$,&
$\{(0,0),(4,7),(2,9)\}$,&
$\{(0,0),(1,8),(3,9)\}$,&
$\{(0,0),(2,8),(4,9)\}$,
\\
$\{(0,0),(3,8),(5,9)\}$,&
$\{(0,0),(5,8),(1,9)\}$,&
$\{(1,0),(4,1),(2,4)\}$,&
$\{(1,0),(5,1),(4,2)\}$,
\\
$\{(1,0),(2,2),(3,6)\}$,&
$\{(1,0),(3,2),(4,4)\}$,&
$\{(1,0),(5,2),(3,4)\}$,&
$\{(1,0),(2,3),(5,6)\}$,
\\
$\{(1,0),(3,3),(2,8)\}$,&
$\{(1,0),(4,3),(3,5)\}$,&
$\{(1,0),(5,3),(4,6)\}$,&
$\{(1,0),(5,4),(2,7)\}$,
\\
$\{(1,0),(2,5),(5,7)\}$,&
$\{(1,0),(4,5),(3,8)\}$,&
$\{(1,0),(2,6),(3,7)\}$,&
$\{(1,0),(4,7),(5,8)\}$,
\\
$\{(2,0),(5,1),(3,7)\}$,&
$\{(2,0),(3,2),(5,5)\}$,&
$\{(2,0),(4,2),(3,6)\}$,&
$\{(2,0),(4,3),(3,8)\}$,
\\
$\{(2,0),(5,4),(4,6)\}$,&
$\{(2,0),(4,5),(5,8)\}$,&
$\{(3,0),(4,3),(5,7)\}$,&
$\{(3,0),(4,4),(5,6)\}$.\\
\end{longtable}

\noindent
All $v(5v-6)$ blocks are obtained by developing these base blocks by $(-,+1)\bmod~ (-,v)$ successively. \qed

\begin{lemma}
There exists a $3$-HGDP of type $(12,1^{4})$ with $J(12\times4\times1,3,1)$  blocks.
\end{lemma}
\proof

 The design will be constructed on  $X=I_{12}\times Z_{4}$ with group set $\mathcal{G}=\{\{i\}\times Z_{4}:i\in I_{12}\}$ and hole set $\mathcal{H}=\{I_{12}\times \{j\}:j\in Z_{4}\}$. We list all $64$ base  blocks as follows.
\begin{longtable}{llll}
$\{(0,0),(1,1),(2,2)\}$,&
$\{(0,0),(2,1),(1,2)\}$,&
$\{(0,0),(3,1),(4,2)\}$,&
$\{(0,0),(4,1),(3,2)\}$,\\
$\{(0,0),(5,1),(6,2)\}$,&
$\{(0,0),(6,1),(5,2)\}$,&
$\{(0,0),(7,1),(1,3)\}$,&
$\{(0,0),(8,1),(2,3)\}$,\\
$\{(1,0),(3,1),(2,2)\}$,&$\{(0,0),(11,2),(9,3)\}$,&
$\{(1,0),(4,1),(5,2)\}$,&
$\{(1,0),(5,1),(3,2)\}$,
\\
$\{(1,0),(7,1),(3,3)\}$,&
$\{(1,0),(8,1),(4,3)\}$,&$\{(1,0),(10,1),(6,3)\}$,&
$\{(1,0),(9,1),(5,3)\}$,\\
$\{(1,0),(11,1),(7,3)\}$,&
$\{(1,0),(6,2),(8,3)\}$,&
$\{(1,0),(9,2),(11,3)\}$,&
$\{(1,0),(8,2),(9,3)\}$,\\
$\{(1,0),(11,2),(10,3)\}$,&
$\{(2,0),(3,1),(6,2)\}$,&
$\{(2,0),(4,1),(7,2)\}$,&
$\{(2,0),(5,1),(4,2)\}$,\\
$\{(2,0),(6,1),(3,2)\}$,&
$\{(2,0),(10,2),(11,3)\}$,&
$\{(2,0),(7,1),(8,3)\}$,&
$\{(2,0),(8,1),(5,3)\}$,\\
$\{(2,0),(10,1),(9,3)\}$,&
$\{(2,0),(9,1),(4,3)\}$,&
$\{(2,0),(11,1),(10,3)\}$,&
$\{(2,0),(9,2),(6,3)\}$,\\
$\{(3,0),(10,1),(4,2)\}$,&
$\{(3,0),(8,2),(10,3)\}$,&
$\{(3,0),(6,2),(11,3)\}$,&$\{(3,0),(9,1),(8,3)\}$,
\\
$\{(5,0),(8,1),(10,3)\}$,&
$\{(5,0),(9,1),(7,3)\}$,&
$\{(5,0),(10,1),(6,2)\}$,&
$\{(5,0),(7,2),(9,3)\}$,\\
$\{(0,0),(9,1),(3,3)\}$,&
$\{(0,0),(10,1),(4,3)\}$,&
$\{(0,0),(11,1),(5,3)\}$,&
$\{(0,0),(7,2),(6,3)\}$,\\
$\{(0,0),(8,2),(7,3)\}$,&
$\{(0,0),(9,2),(10,3)\}$,&
$\{(0,0),(10,2),(8,3)\}$,&
$\{(1,0),(6,1),(4,2)\}$,\\
$\{(2,0),(11,2),(7,3)\}$,&
$\{(3,0),(5,1),(11,2)\}$,&
$\{(3,0),(11,1),(9,3)\}$,&
$\{(3,0),(8,1),(5,2)\}$,\\
$\{(3,0),(10,2),(7,3)\}$,&
$\{(4,0),(8,1),(11,3)\}$,&
$\{(4,0),(10,1),(7,3)\}$,&
$\{(4,0),(9,1),(7,2)\}$,\\
$\{(4,0),(11,1),(5,2)\}$,&
$\{(4,0),(11,2),(8,3)\}$,&
$\{(5,0),(7,1),(10,2)\}$,&
$\{(4,0),(6,2),(9,3)\}$,\\
$\{(6,0),(7,1),(11,2)\}$,&
$\{(6,0),(10,1),(9,2)\}$,&
$\{(6,0),(8,2),(11,3)\}$,&
$\{(6,0),(7,2),(8,3)\}$.\\
\end{longtable}
\noindent
All blocks can be obtained by developing the above  base blocks by $(-,+1)\bmod(-,4)$ successively.\qed
\begin{lemma}
There exists a $3$-HGDP of type $(12,1^{6})$ with $J(12\times6\times1,3,1)$ blocks.
\end{lemma}
\proof
 We construct the desired design on $X=Z_{12}\times I_{6}$ with group set $\mathcal{G}=\{\{i\}\times I_{6}:i\in Z_{12}\}$ and hole set $\mathcal{H}=\{Z_{12}\times \{j\}:j\in I_{6}\}$. All $54$ base blocks are listed as follows.

\begin{longtable}{lllll}
$\{(0,0),(1,1),(2,2)\}$,&
$\{(0,0),(2,1),(7,4)\}$,&
$\{(0,0),(3,1),(4,4)\}$,&
$\{(0,0),(4,1),(1,2)\}$,\\
$\{(0,0),(5,1),(10,2)\}$,&
$\{(0,0),(6,1),(10,3)\}$,&
$\{(0,0),(9,1),(10,4)\}$,&
$\{(0,0),(8,1),(3,4)\}$,\\
$\{(0,0),(10,1),(4,3)\}$,&
$\{(0,0),(11,1),(2,4)\}$,&
$\{(0,0),(3,2),(8,4)\}$,&
$\{(0,0),(4,2),(2,3)\}$,\\
$\{(0,0),(5,2),(8,5)\}$,&
$\{(0,0),(6,2),(8,3)\}$,&
$\{(0,0),(7,2),(5,3)\}$,&
$\{(0,0),(8,2),(9,3)\}$,\\
$\{(0,0),(9,2),(2,5)\}$,&
$\{(0,0),(11,2),(1,4)\}$,&
$\{(0,0),(3,3),(10,5)\}$,&
$\{(0,0),(1,3),(4,5)\}$,\\
$\{(0,0),(6,3),(11,5)\}$,&
$\{(0,0),(7,3),(5,5)\}$,&
$\{(0,0),(11,3),(9,5)\}$,&
$\{(0,0),(5,4),(6,5)\}$,\\
$\{(0,0),(6,4),(7,5)\}$,&
$\{(0,0),(9,4),(3,5)\}$,&
$\{(0,0),(11,4),(1,5)\}$,&
$\{(1,0),(3,1),(6,2)\}$,\\
$\{(1,0),(4,1),(3,2)\}$,&
$\{(1,0),(5,1),(3,4)\}$,&
$\{(1,0),(7,1),(11,2)\}$,&
$\{(1,0),(6,1),(9,2)\}$,\\
$\{(1,0),(8,1),(6,4)\}$,&
$\{(1,0),(10,1),(2,5)\}$,&
$\{(1,0),(11,1),(9,5)\}$,&
$\{(1,0),(9,1),(8,4)\}$,\\
$\{(1,0),(2,2),(8,3)\}$,&
$\{(1,0),(5,2),(6,5)\}$,&
$\{(1,0),(8,2),(10,4)\}$,&
$\{(1,0),(7,2),(9,4)\}$,\\
$\{(1,0),(10,2),(7,4)\}$,&
$\{(1,0),(2,3),(5,4)\}$,&
$\{(1,0),(3,3),(4,4)\}$,&
$\{(1,0),(4,3),(7,5)\}$,\\
$\{(1,0),(6,3),(8,5)\}$,&
$\{(1,0),(9,3),(10,5)\}$,&
$\{(1,0),(10,3),(3,5)\}$,&
$\{(1,0),(7,3),(2,4)\}$,\\
$\{(1,0),(11,3),(5,5)\}$,&
$\{(2,0),(3,1),(6,4)\}$,&
$\{(2,0),(4,1),(5,2)\}$,&
$\{(2,0),(6,1),(7,4)\}$,\\
$\{(2,0),(7,1),(3,5)\}$,&
$\{(2,0),(9,1),(8,3)\}$.&
\end{longtable}
\noindent
All blocks can be obtained by developing the above base blocks by $(+1,-)\bmod(12,-)$  successively.\qed

\begin{lemma}
There exists a $3$-HGDP of type $(12,1^{10})$ with $J(12\times 10\times 1,3,1)$ blocks.
\end{lemma}
\proof
 We construct the desired $3$-HGDP of type $(12,1^{10})$ on $X=I_{120}$ with group set $\mathcal{G}=\{\{12i+j:i\in I_{10}\}:j\in I_{12}\}$ and hole  set $\mathcal{H}=\{\{12i+j:j\in I_{12}\}:i\in I_{10} \}$. We list $96$ base blocks as follows.

\begin{longtable}{lllllll}
$\{0,14,25\}$,&
$\{0,15,41\}$,&
$\{2,36,91\}$,&
$\{0,17,94\}$,&
$\{2,44,83\}$,&
$\{0,19,89\}$,\\
$\{0,29,66\}$,&
$\{0,21,78\}$,&
$\{0,22,79\}$,&
$\{0,23,82\}$,&
$\{1,59,80\}$,&
$\{0,28,65\}$,&
\\$\{2,46,85\}$,
&
$\{2,29,59\}$,&
$\{0,33,55\}$,&
$\{0,35,39\}$,&
$\{1,43,92\}$,&
$\{0,40,62\}$,\\
$\{0,42,50\}$,&
$\{0,43,51\}$,&
$\{1,36,74\}$,&
$\{0,45,85\}$,&
$\{0,46,64\}$,&
$\{0,47,52\}$,\\
$\{0,53,61\}$,&
$\{0,56,67\}$,&
$\{1,17,26\}$,&
$\{0,32,37\}$,&
$\{0,68,76\}$,&
$\{0,69,97\}$,\\
$\{0,70,73\}$,&
$\{0,71,74\}$,&
$\{0,80,86\}$,&
$\{0,83,90\}$,&
$\{1,31,70\}$,&
$\{0,87,98\}$,\\
$\{1,71,89\}$,&
$\{2,47,88\}$,&$\{0,20,26\}$,&
$\{2,64,93\}$,&
$\{0,88,113\}$,&$\{2,70,90\}$,\\
$\{1,12,90\}$,&
$\{1,16,58\}$,&
$\{0,57,103\}$,&
$\{1,22,50\}$,&
$\{1,20,101\}$,&
$\{1,19,66\}$,\\
$\{2,12,37\}$,&
$\{1,21,65\}$,&
$\{1,23,114\}$,&
$\{1,28,56\}$,&
$\{2,68,112\}$,&
$\{1,32,41\}$,\\
$\{1,34,98\}$,&
$\{1,35,63\}$,&
$\{0,44,111\}$,&
$\{2,40,61\}$,&
$\{0,38,116\}$,&
$\{1,44,60\}$,\\
$\{1,45,55\}$,&
$\{1,46,53\}$,&
$\{1,88,104\}$,&
$\{1,52,108\}$,&
$\{0,27,104\}$,&
$\{1,67,94\}$,\\
$\{0,91,110\}$,&
$\{1,83,84\}$,&
$\{1,86,112\}$,&
$\{1,87,119\}$,&
$\{1,95,117\}$,&
$\{1,47,68\}$,\\
$\{1,18,115\}$,&
$\{0,93,114\}$,&
$\{0,16,107\}$,&
$\{1,42,107\}$,&
$\{0,18,100\}$,&
$\{2,45,80\}$,\\
$\{12,40,93\}$,&
$\{0,92,109\}$,&
$\{12,46,86\}$,&
$\{0,95,118\}$,&$\{0,106,119\}$,
&
$\{2,66,89\}$,\\
$\{13,41,86\}$,&
$\{2,63,109\}$,&
$\{2,71,108\}$,&
$\{0,59,117\}$, &$\{0,101,115\}$,&
$\{2,65,92\}$,\\
$\{12,43,87\}$,&$\{12,61,118\}$,&
$\{12,42,89\}$,&
$\{12,62,113\}$,&$\{12,39,119\}$,
&$\{0,30,49\}$,\\
$\{2,67,117\}$,&$\{13,67,113\}$.
\end{longtable}
\noindent
Let $G$ be the group generated by $\alpha$ and $\beta$ where $\alpha=\prod_{i=0}^{23}(i~i+24~i+48~i+72~ i+96)$ and $\beta=\prod_{j=0}^{24}(j~j+3~j+6~ j+9)$. All blocks are obtained by developing the above base blocks under the action of $G$.
\qed

\begin{lemma}
There exists a $3$-HGDP of type $(v,1^{12})$ with $J(v\times 12\times 1,3,1)$ blocks for $v\in\{8,14\}.$
\end{lemma}
\proof
  The design is constructed on $X=I_{12v}$ with group set $\mathcal{G}=\{\{12i+j:j\in I_{12}\}:i\in I_{v} \}$ and hole set $\mathcal{H}=\{\{12i+j:i\in I_{v}\}:j\in I_{12}\}$. Next we list  $11v-12$ base blocks as follows.

$v=8:$

\begin{longtable}{lllllllll}
 $\{ 0, 13, 27 \}$,& $\{ 0, 15, 55 \}$,&
  $\{ 0, 16, 65 \}$,& $\{ 0, 17, 78 \}$,&
   $\{ 0, 20, 79 \}$,& $\{ 0, 19, 90 \}$,&
    $\{ 0, 21, 54 \}$,\\
    $\{ 0, 22, 89 \}$,&
  $\{ 0, 18, 52 \}$,& $\{ 0, 53, 66 \}$,& $\{ 0, 23, 64 \}$,& $\{ 4, 54, 77 \}$,& $\{ 0, 76, 91 \}$,& $\{ 0, 67, 77 \}$,\\
  $\{ 4, 55, 69 \}$,& $\{ 4, 58, 78 \}$,& $\{ 4, 59, 67 \}$,& $\{ 4, 57, 89 \}$,& $\{ 4, 66, 92 \}$,& $\{ 8, 52, 71 \}$,& $\{ 4, 82, 91 \}$,\\
  $\{ 8, 54, 95 \}$,& $\{ 8, 57, 76 \}$,& $\{ 8, 55, 82 \}$,& $\{ 8, 58, 64 \}$,& $\{ 8, 59, 90 \}$,& $\{ 8, 49, 91 \}$,& $\{ 8, 48, 65 \}$,\\
  $\{ 8, 63, 88 \}$,& $\{ 8, 51, 66 \}$,& $\{ 8, 50, 78 \}$,&$\{ 0, 58, 71 \}$ ,& $\{ 0, 68, 82 \}$,& $\{ 4, 17, 42 \}$,& $\{ 8, 77, 85 \}$,\\
  $\{ 4, 20, 45 \}$,& $\{ 0, 32, 45 \}$,& $\{ 0, 33, 44 \}$,& $\{ 4, 21, 35 \}$,& $\{ 0, 28, 46 \}$,& $\{ 0, 34, 87 \}$,& $\{ 0, 35, 56 \}$,\\
  $\{ 0, 47, 49 \}$,& $\{ 4, 23, 75 \}$,& $\{ 4, 32, 63 \}$,& $\{ 4, 33, 61 \}$,& $\{ 4, 34, 51 \}$,& $\{ 4, 44, 60 \}$,& $\{ 8, 69, 74 \}$,\\
  $\{ 4, 47, 62 \}$,& $\{ 4, 46, 80 \}$,& $\{ 8, 70, 93 \}$,& $\{ 8, 83, 86 \}$,& $\{ 0, 29, 50 \}$,& $\{ 0, 31, 70 \}$,& $\{ 0, 51, 61 \}$,\\
  $\{ 0, 42, 80 \}$,& $\{ 0, 62, 83 \}$,& $\{ 0, 30, 94 \}$,& $\{ 0, 40, 81 \}$,& $\{ 0, 41, 86 \}$,& $\{ 4, 56, 74 \}$,& $\{ 4, 71, 84 \}$,\\
  $\{ 0, 57, 85 \}$,& $\{ 4, 50, 83 \}$,& $\{ 4, 81, 85 \}$,& $\{ 4, 19, 87 \}$,& $\{ 0, 43, 73 \}$,& $\{ 4, 48, 70 \}$,& $\{ 0, 63, 95 \}$,\\
  $\{ 0, 74, 93 \}$,& $\{ 0, 69, 75 \}$,& $\{ 0, 59, 92 \}$,& $\{ 48, 63, 79 \}$,& $\{ 48, 61, 91 \}$,& $\{ 48, 66, 73\}.$
\end{longtable}

 Let $G$ be the group generated by $\alpha$ and $\beta$ where $$\alpha=\prod_{i=0}^{11}\prod_{l=0}^{1}(i+48l~i+48vl+12~i+48vl+24~i+48vl+36),$$ $$\beta=\prod_{j=0}^{7}\prod_{l=0}^{2}(12j+4l~12j+4l+1~12j+4l+2~12j+4l+3).$$ All blocks are obtained by developing the above base blocks under the action of $G$.

$v=14$:

\begin{longtable}{llllllllllll}
 $\{ 0, 13, 27 \}$, & $\{ 84, 101, 152 \}$,&$\{ 8, 33, 151 \}$,& $\{ 0, 17, 104 \}$,& $\{ 0, 19, 143 \}$,& $\{ 0, 18, 41 \}$,\\
 $\{ 0, 20, 154 \}$,&  $\{ 4, 18, 149 \}$,& $\{ 0, 28, 167 \}$,& $\{ 0, 23, 117 \}$,& $\{ 0, 29, 131 \}$,& $\{ 0, 22, 38 \}$,\\
   $\{ 0, 30, 129 \}$,&  $\{ 4, 116, 127 \}$,&
  $\{ 0, 33, 164 \}$,& $\{ 84, 119, 153 \}$,& $\{ 0, 35, 106 \}$,& $\{ 0, 34, 42 \}$,\\ $\{ 0, 40, 116 \}$,& $\{ 4, 102, 166 \}$,& $\{ 0, 44, 153 \}$,&
   $\{ 0, 113, 124 \}$,& $\{ 0, 123, 140 \}$,& $\{ 0, 47, 78 \}$,\\
  $\{ 0, 55, 128 \}$,&   $\{ 8, 100, 119 \}$,& $\{ 0, 56, 142 \}$,&
  $\{ 0, 66, 105 \}$,& $\{ 0, 65, 155 \}$,& $\{ 0, 52, 77 \}$,\\
  $\{ 0, 57, 119 \}$,& $\{ 0, 69, 130 \}$,& $\{ 0, 59, 134 \}$,&$\{ 0, 71, 159 \}$,& $\{ 0, 76, 157 \}$,&
  $\{ 0, 80, 87 \}$,\\
  $\{ 0, 81, 101 \}$,& $\{ 0, 82, 114 \}$,& $\{ 0, 83, 139 \}$,& $\{ 0, 85, 151 \}$,& $\{ 0, 89, 122 \}$,& $\{ 0, 86, 97 \}$,\\
  $\{ 0, 92, 109 \}$,& $\{ 0, 91, 135 \}$,& $\{ 0, 93, 158 \}$,& $\{ 0, 100, 146 \}$,& $\{ 0, 98, 121 \}$,&$\{ 0, 90, 99 \}$,\\ $\{ 0, 102, 145 \}$,&
  $\{ 0, 103, 161 \}$,& $\{ 0, 107, 148 \}$,& $\{ 0, 110, 138 \}$,& $\{ 0, 111, 137 \}$,& $\{0,46,53\},$\\
  $\{ 0, 118, 162 \}$,& $\{ 0, 127, 133 \}$,& $\{ 0, 115, 147 \}$,& $\{ 0, 126, 149 \}$,& $\{ 0, 125, 160 \}$,&
  $\{0,45,67\}$,\\
   $\{ 0, 152, 166 \}$,& $\{ 4, 92, 119 \}$,& $\{ 4, 89, 167 \}$,& $\{ 4, 93, 118 \}$,& $\{ 4, 104, 130 \}$,&
   $\{0,31,94\}$,
     \\
  $\{ 4, 140, 153 \}$,& $\{ 4, 103, 117 \}$,& $\{ 8, 95, 130 \}$,& $\{ 4, 101, 154 \}$,&
  $\{ 4, 113, 141 \}$,& $\{0,43,54\}$,
  \\
  $\{ 4, 114, 152 \}$,& $\{ 8, 88, 129 \}$,& $\{ 8, 93, 166 \}$,& $\{ 8, 105, 135 \}$,&
  $\{ 8, 98, 118 \}$,& $\{0,58,64\},$
 \\
  $\{ 8, 131, 150 \}$,& $\{ 8, 120, 165 \}$,& $\{ 8, 89, 141 \}$,& $\{ 84, 97, 167 \}$,&
  $\{ 84, 98, 131 \}$, & $\{ 0, 16, 25 \}$, \\
  $\{ 84, 106, 146 \}$, &$\{ 84, 109, 164 \}$,& $\{ 84, 127, 140 \}$,&
  $\{ 84, 100, 130 \}$,& $\{ 84, 103, 166 \}$,& $\{0,32,95\}$,\\
  $\{ 88, 102, 164 \}$,& $\{ 4, 90, 115 \}$,&
  $\{ 4, 91, 125 \}$,& $\{ 4, 59, 126 \}$,& $\{ 4, 137, 150 \}$,& $\{0,21,39\}$,\\
  $\{ 4, 108, 161 \}$,&
  $\{ 8, 21, 138 \}$,& $\{ 4, 47, 162 \}$,& $\{ 8, 22, 101 \}$,& $\{ 8, 23, 113 \}$,& $\{0,15,37\}$,\\
   $\{ 8, 86, 149 \}$, &$\{ 8, 34, 137 \}$,&$\{ 8, 45, 114 \}$,& $\{ 8, 103, 111 \}$,& $\{ 8, 90, 145 \}$,&  $\{ 4, 23, 99 \}$, \\
   $\{ 84, 115, 122 \}$,&
  $\{ 4, 20, 145 \}$,& $\{ 4, 21, 134 \}$,& $\{ 4, 22, 83 \}$,& $\{ 8, 136, 156 \}$,
  & $\{ 8, 46, 97 \}$,\\
  $\{ 4, 32, 146 \}$,&
  $\{ 4, 33, 144 \}$,& $\{ 4, 35, 135 \}$,& $\{ 4, 58, 121 \}$,& $\{ 4, 46, 111 \}$, & $\{ 4, 44, 96 \}$,\\ $\{ 4, 71, 84 \}$,&
  $\{ 4, 69, 87 \}$,& $\{ 4, 56, 110 \}$,& $\{ 4, 34, 158 \}$, & $\{ 4, 123, 156 \}$,& $\{ 4, 43, 98 \}$,\\ $\{ 4, 17, 122 \}$,&
  $\{ 4, 30, 159 \}$,&$\{ 84, 102, 118 \}$,& $\{ 4, 41, 157 \}$,&  $\{ 84, 116, 123 \}$, & $\{ 8, 47, 85 \}$,\\
   $\{ 0, 112, 163 \}$,& $\{ 8, 107, 160 \}$,& $\{ 4, 68, 138 \}$,& $\{ 4, 97, 151 \}$,& $\{ 8, 139, 158 \}$,& $\{ 4, 42, 86 \}$,\\

  $\{ 0, 141, 150 \}$,& $\{ 0, 136, 165 \}$,& $\{ 4, 105, 142 \}$,& $\{ 4, 129, 163 \}$.

\end{longtable}
\noindent
 Let $G$ be the group generated by $\alpha$ and $\beta$ where $$\alpha=\prod_{i=0}^{11}\prod_{l=0}^{1}(i+84l~i+84vl+12~i+84vl+24~i+84vl+36~i+84vl+48~i+84vl+60~i+84vl+72),$$ $$\beta=\prod_{j=0}^{13}\prod_{l=0}^{2}(12j+4l~12j+4l+1~12j+4l+2~12j+4l+3).$$  All blocks are obtained by developing the above base blocks under the action of $G$.
\qed
\begin{lemma}
There exists a $3$-HGDP of type $(4,1^8)$ with $J(4\times8\times1,3,1)$ blocks.
\end{lemma}
\proof
The design will be constructed on $X=I_{32}$ with group set $\mathcal{G}=\{\{4i+j:i\in I_{8}\}:j\in I_{4}\}$ and hole set $\mathcal{H}=\{\{4j+i:i\in I_{4}\}:j\in I_{8}\}$. All $106$ blocks are listed as follows.

\centering
\begin{longtable}{lllllll}
 $\{ 0, 7, 26 \}$,& $\{ 0, 8, 22 \}$,& $\{ 0, 7, 21 \}$,& $\{ 0, 9, 23 \}$,& $\{ 0, 9, 33 \}$,& $\{ 0, 11, 18 \}$, & $\{ 0, 9, 13 \}$, \\$\{ 0, 11, 28 \}$,&
  $\{ 1, 12, 23 \}$,& $\{ 2, 6, 24 \}$,& $\{ 0, 24, 38 \}$,& $\{ 2, 6, 18 \}$,& $\{ 2, 8, 35 \}$,& $\{ 2, 6, 39 \}$, \\$\{ 2, 9, 38 \}$,& $\{ 2, 11, 19 \}$,&
  $\{ 2, 10, 26 \}$,& $\{ 0, 19, 28 \}$,& $\{ 1, 22, 28 \}$,& $\{ 3, 10, 21 \}$, & $\{ 3, 6, 20 \}$,\\$\{ 4, 10, 18 \}$,&
  $\{ 2, 15, 36 \}$,& $\{ 4, 10, 16 \}$,& $\{ 0, 16, 32 \}$,& $\{ 2, 29, 36 \}$,& $\{ 1, 15, 37 \}$, & $\{ 1, 7, 35 \}$,\\

  $\{ 1, 10, 19 \}$,& $\{ 1, 13, 29 \}$,&
  $\{ 2, 5, 26 \}$,& $\{ 3, 7, 25 \}$,& $\{ 0, 18, 27 \}$,& $\{ 4, 13, 17 \}$, & $\{ 4, 5, 12 \}$,\\
  $\{ 4, 17, 33 \}$,& $\{ 1, 8, 20 \}$,& $\{ 2, 13, 39 \}$,&
  $\{ 3, 17, 36 \}$,& $\{ 2, 26, 39 \}$,& $\{ 2, 16, 33 \}$,& $\{ 3, 6, 25 \}$, \\$\{ 4, 7, 20 \}$,& $\{ 3, 11, 17 \}$,& $\{ 4, 5, 36 \}$,&
  $\{ 4, 16, 28 \}$,& $\{ 3, 21, 25 \}$,& $\{ 0, 17, 33 \}$, & $\{ 1, 5, 22 \}$,\\$\{ 3, 22, 30 \}$,& $\{ 1, 23, 34 \}$,& $\{ 3, 24, 26 \}$,&
  $\{ 1, 14, 17 \}$,& $\{ 1, 15, 28 \}$,& $\{ 1, 18, 24 \}$,&
  $\{ 4, 8, 10 \}$, \\
  $\{ 1, 33, 35 \}$,& $\{ 3, 27, 35 \}$,& $\{ 3, 12, 15 \}$,&
  $\{ 4, 12, 16 \}$,& $\{ 3, 17, 20 \}$,&$\{ 2, 18, 35 \}$, &
   $\{ 4, 8, 20 \}$, \\

  $\{ 0, 22, 31 \}$,& $\{ 0, 23, 34 \}$,&$\{ 2, 29, 38 \}$,&
  $\{ 4, 5, 18 \}$,& $\{ 3, 16, 24 \}$,& $\{ 0, 13, 22 \}$,& $\{ 4, 6, 33 \}$, \\$\{ 0, 12, 23 \}$,& $\{ 0, 13, 34 \}$,& $\{ 0, 21, 33 \}$,&
  $\{ 2, 25, 33 \}$,& $\{ 2, 29, 31 \}$,& $\{ 1, 20, 28 \}$, & $\{ 3, 7, 11 \}$,\\$\{ 4, 13, 32 \}$,& $\{ 3, 5, 29 \}$,& $\{ 2, 14, 36 \}$,&
  $\{ 3, 5, 17 \}$,& $\{ 1, 24, 38 \}$,& $\{ 2, 5, 23 \}$,& $\{ 2, 5, 14 \}$, \\$\{ 2, 11, 38 \}$,& $\{ 1, 13, 15 \}$,& $\{ 2, 14, 26 \}$,&
  $\{ 2, 13, 30 \}$,& $\{ 3, 29, 36 \}$,& $\{ 0, 27, 31 \}$,&
   $\{ 1, 9, 25 \}$, \\
  $\{ 2, 24, 30 \}$,& $\{ 3, 14, 21 \}$,& $\{ 3, 17, 31 \}$,&
  $\{ 3, 27, 30 \}$,& $\{ 1, 8, 39 \}$,& $\{ 0, 31, 38 \}$,&
  $\{ 1, 9, 33 \}$, \\$\{ 2, 25, 38 \}$,&$\{ 0, 34, 36 \}$,  & $\{ 1, 19, 28 \}$,&
  $\{ 1, 14, 25 \}$,& $\{ 3, 12, 39 \}$,& $\{ 4, 12, 18 \}$,& $\{ 3, 9, 30 \}$, \\
  $\{ 3, 11, 35 \}$.&&&&&&\qed
 \end{longtable}

\begin{lemma}
There exists a $3$-HGDP of type $(4,1^{14})$ with $J(4\times14\times1,3,1)$  blocks.
\end{lemma}
\proof
The design will be constructed on $X=I_{56}$ with group set $\mathcal{G}=\{\{4i+j:i\in I_{14}\}:j\in I_{4}\}$ and hole set $\mathcal{H}=\{\{4j+i:i\in I_{4}\}:j\in I_{14}\}$. All $354$ blocks are listed as follows.

\centering
\begin{longtable}{llllllll}
$\{ 0, 5, 14 \}$,& $\{ 0, 6, 19 \}$,& $\{ 0, 15, 17 \}$,& $\{ 0, 13, 18 \}$,& $\{ 5, 8, 22 \}$,&$\{ 17, 30, 32 \}$,  & $\{ 7, 8, 29 \}$, \\
$\{ 8, 21, 30 \}$,&
  $\{ 4, 13, 30 \}$,& $\{ 4, 22, 31 \}$,& $\{ 4, 23, 29 \}$,& $\{ 4, 15, 21 \}$,& $\{ 18, 31, 32 \}$,& $\{ 6, 8, 23 \}$,\\
  $\{ 14, 19, 32 \}$,&
  $\{ 15, 29, 32 \}$,& $\{ 6, 15, 44 \}$,& $\{ 7, 21, 44 \}$,& $\{ 13, 22, 44 \}$,& $\{ 14, 23, 44 \}$,& $\{ 5, 18, 52 \}$, \\$\{ 7, 17, 52 \}$,&
  $\{ 19, 29, 52 \}$,& $\{ 6, 31, 52 \}$,& $\{ 7, 13, 48 \}$,& $\{ 5, 30, 48 \}$,& $\{ 14, 31, 48 \}$,& $\{ 0, 25, 50 \}$, \\$\{ 0, 26, 51 \}$,&
  $\{ 0, 27, 49 \}$,& $\{ 0, 29, 10 \}$,& $\{ 0, 30, 11 \}$,& $\{ 0, 31, 9 \}$,& $\{ 0, 45, 34 \}$,& $\{ 0, 46, 35 \}$, \\$\{ 0, 47, 33 \}$,&
  $\{ 0, 41, 54 \}$,& $\{ 0, 42, 55 \}$,& $\{ 0, 43, 53 \}$,& $\{ 0, 37, 22 \}$,& $\{ 0, 38, 23 \}$,& $\{ 0, 39, 21 \}$, \\$\{ 4, 37, 2 \}$,&
  $\{ 4, 38, 3 \}$,& $\{ 4, 39, 1 \}$,& $\{ 4, 53, 10 \}$,& $\{ 4, 54, 11 \}$,& $\{ 4, 55, 9 \}$,& $\{ 4, 33, 50 \}$, \\$\{ 4, 34, 51 \}$,&
  $\{ 4, 35, 49 \}$,& $\{ 4, 17, 42 \}$,& $\{ 4, 18, 43 \}$,& $\{ 4, 19, 41 \}$,& $\{ 4, 25, 46 \}$,& $\{ 4, 26, 47 \}$, \\$\{ 4, 27, 45 \}$,&
  $\{ 8, 2, 33 \}$,& $\{ 8, 34, 3 \}$,& $\{ 8, 35, 1 \}$,& $\{ 8, 49, 14 \}$,& $\{ 8, 50, 15 \}$,& $\{ 8, 51, 13 \}$, \\$\{ 8, 45, 38 \}$,&
  $\{ 8, 46, 39 \}$,& $\{ 8, 47, 37 \}$,& $\{ 8, 25, 42 \}$,& $\{ 8, 26, 43 \}$,& $\{ 8, 27, 41 \}$,& $\{ 8, 17, 54 \}$, \\$\{ 8, 18, 55 \}$,& $\{ 12, 6, 45 \}$,& $\{ 12, 5, 47 \}$,& $\{ 12, 21, 54 \}$,& $\{ 12, 22, 55 \}$, & $\{ 12, 51, 29 \}$, &
  $\{ 8, 19, 53 \}$,\\
  $\{ 12, 23, 53 \}$,&
  $\{ 12, 37, 42 \}$,& $\{ 12, 38, 43 \}$,& $\{ 12, 39, 41 \}$,& $\{ 12, 49, 30 \}$,& $\{ 12, 50, 31 \}$,& $\{ 12, 7, 46 \}$,\\

  $\{ 12, 33, 26 \}$,&
  $\{ 12, 34, 27 \}$,& $\{ 12, 35, 25 \}$,& $\{ 16, 41, 6 \}$,& $\{ 16, 42, 7 \}$,& $\{ 16, 43, 5 \}$,& $\{ 16, 9, 50 \}$, \\$\{ 16, 10, 51 \}$,&
  $\{ 29, 26, 3 \}$,& $\{ 31, 25, 2 \}$,& $\{ 29, 14, 47 \}$,& $\{ 30, 15, 45 \}$,& $\{ 31, 13, 46 \}$, & $\{ 30, 27, 1 \}$, \\$\{ 29, 6, 43 \}$,&
  $\{ 30, 7, 41 \}$,& $\{ 29, 34, 39 \}$,& $\{ 30, 35, 37 \}$,& $\{ 31, 33, 38 \}$,& $\{ 29, 46, 11 \}$, & $\{ 31, 5, 42 \}$, \\$\{ 30, 47, 9 \}$,&
  $\{ 31, 45, 10 \}$,& $\{ 29, 50, 35 \}$,& $\{ 30, 51, 33 \}$,& $\{ 31, 49, 34 \}$,& $\{ 33, 6, 55 \}$,& $\{ 34, 7, 53 \}$, \\$\{ 35, 5, 54 \}$,& $\{ 34, 15, 5 \}$,& $\{ 35, 13, 6 \}$,& $\{ 33, 18, 15 \}$,& $\{ 34, 19, 13 \}$,& $\{ 35, 17, 14 \}$, &
  $\{ 33, 14, 7 \}$, \\$\{ 33, 54, 39 \}$,&
  $\{ 34, 55, 37 \}$,& $\{ 35, 53, 38 \}$,& $\{ 37, 10, 43 \}$,& $\{ 38, 11, 41 \}$,& $\{ 39, 9, 42 \}$,& $\{ 37, 14, 3 \}$, \\$\{ 38, 15, 1 \}$,& $\{ 37, 18, 11 \}$,& $\{ 38, 19, 9 \}$,& $\{ 39, 17, 10 \}$,& $\{ 37, 50, 19 \}$,& $\{ 38, 51, 17 \}$, &
  $\{ 39, 13, 2 \}$, \\$\{ 39, 49, 18 \}$,& $\{ 42, 47, 1 \}$,& $\{ 43, 45, 2 \}$,& $\{ 45, 54, 23 \}$,& $\{ 46, 55, 21 \}$,& $\{ 47, 53, 22 \}$,&
  $\{ 41, 46, 3 \}$, \\
  $\{ 45, 26, 55 \}$,&
  $\{ 46, 27, 53 \}$,& $\{ 47, 25, 54 \}$,& $\{ 9, 46, 43 \}$,& $\{ 10, 47, 41 \}$,& $\{ 11, 45, 42 \}$,& $\{ 13, 10, 3 \}$, \\$\{ 18, 23, 41 \}$,&
  $\{ 15, 9, 2 \}$,& $\{ 17, 26, 31 \}$,& $\{ 18, 27, 29 \}$,& $\{ 19, 25, 30 \}$,& $\{ 17, 22, 43 \}$,& $\{ 14, 11, 1 \}$, \\$\{ 19, 21, 42 \}$,&
  $\{ 20, 53, 2 \}$,& $\{ 20, 54, 3 \}$,& $\{ 20, 45, 50 \}$, & $\{ 20, 25, 6 \}$,& $\{ 20, 26, 7 \}$,& $\{ 20, 27, 5 \}$, \\$\{ 20, 41, 14 \}$,&
  $\{ 20, 42, 15 \}$,& $\{ 20, 43, 13 \}$,& $\{ 20, 29, 38 \}$,& $\{ 20, 30, 39 \}$,& $\{ 20, 31, 37 \}$,& $\{ 20, 55, 1 \}$, \\
  $\{ 20, 46, 51 \}$,&
  $\{ 20, 47, 49 \}$,& $\{ 24, 45, 14 \}$,& $\{ 24, 46, 15 \}$,& $\{ 24, 47, 13 \}$,& $\{ 24, 17, 6 \}$,& $\{ 24, 18, 7 \}$, \\$\{ 24, 23, 37 \}$,&
  $\{ 24, 49, 42 \}$,& $\{ 24, 50, 43 \}$,& $\{ 24, 51, 41 \}$,& $\{ 24, 21, 38 \}$,& $\{ 24, 22, 39 \}$,&  $\{ 24, 19, 5 \}$,\\
  $\{ 24, 29, 54 \}$,&
  $\{ 24, 30, 55 \}$,& $\{ 24, 31, 53 \}$,& $\{ 48, 9, 54 \}$,& $\{ 48, 10, 55 \}$,& $\{ 48, 11, 53 \}$,& $\{ 48, 41, 2 \}$,\\
  $\{ 48, 18, 47 \}$,&
  $\{ 48, 19, 45 \}$,& $\{ 48, 17, 46 \}$, & $\{ 48, 43, 1 \}$,& $\{ 48, 25, 38 \}$,& $\{ 48, 26, 39 \}$, & $\{ 48, 42, 3 \}$,\\

  $\{ 48, 27, 37 \}$,&
  $\{ 48, 33, 22 \}$,& $\{ 48, 34, 23 \}$,& $\{ 48, 35, 21 \}$,& $\{ 40, 29, 2 \}$,& $\{ 40, 30, 3 \}$,& $\{ 40, 31, 1 \}$, \\$\{ 40, 37, 6 \}$,& $\{ 40, 39, 5 \}$,& $\{ 40, 53, 14 \}$,& $\{ 40, 54, 15 \}$,& $\{ 40, 55, 13 \}$,& $\{ 40, 45, 18 \}$, &
  $\{ 40, 38, 7 \}$, \\$\{ 40, 46, 19 \}$,&
  $\{ 40, 47, 17 \}$,& $\{ 40, 21, 50 \}$,& $\{ 40, 22, 51 \}$,& $\{ 40, 23, 49 \}$,& $\{ 44, 33, 10 \}$,&$\{ 44, 17, 2 \}$,\\$\{ 44, 35, 9 \}$,& $\{ 44, 34, 11 \}$,
  & $\{ 44, 18, 3 \}$,& $\{ 44, 49, 38 \}$,& $\{ 44, 50, 39 \}$,& $\{ 44, 51, 37 \}$, & $\{ 44, 19, 1 \}$,\\$\{ 44, 29, 42 \}$,&
  $\{ 44, 30, 43 \}$,& $\{ 44, 31, 41 \}$,& $\{ 44, 53, 26 \}$,& $\{ 44, 54, 27 \}$,& $\{ 44, 55, 25 \}$,& $\{ 1, 10, 12 \}$,\\
  $\{ 1, 24, 34 \}$,& $\{ 10, 24, 35 \}$,& $\{ 18, 20, 35 \}$,&$\{ 3, 9, 24 \}$, & $\{ 9, 12, 18 \}$,& $\{ 11, 12, 17 \}$, & $\{ 2, 12, 19 \}$,\\
  $\{ 11, 20, 33 \}$,&
  $\{ 27, 33, 40 \}$,& $\{ 10, 19, 20 \}$,& $\{ 17, 20, 34 \}$,& $\{ 26, 35, 40 \}$,& $\{ 11, 25, 40 \}$, & $\{ 9, 34, 40 \}$,\\
  $\{ 1, 18, 36 \}$,& $\{ 1, 6, 51 \}$,& $\{ 5, 28, 38 \}$,& $\{ 27, 32, 38 \}$,& $\{ 13, 38, 52 \}$,& $\{ 16, 38, 47 \}$,&
  $\{ 6, 28, 39 \}$, \\
  $\{ 25, 32, 39 \}$,& $\{ 14, 39, 52 \}$,& $\{ 16, 39, 45 \}$,& $\{ 26, 32, 37 \}$,& $\{ 15, 37, 52 \}$,& $\{ 16, 37, 46 \}$, & $\{ 7, 28, 37 \}$,\\$\{ 16, 22, 29 \}$,&
  $\{ 29, 36, 55 \}$,& $\{ 14, 16, 55 \}$,& $\{ 16, 23, 30 \}$,& $\{ 15, 16, 53 \}$,& $\{ 16, 25, 34 \}$, & $\{ 3, 16, 33 \}$, \\$\{ 30, 36, 53 \}$,& $\{ 1, 16, 26 \}$,& $\{ 1, 7, 54 \}$,& $\{ 1, 22, 52 \}$,& $\{ 33, 46, 52 \}$,& $\{ 23, 28, 46 \}$,&
  $\{ 5, 36, 46 \}$, \\$\{ 1, 32, 46 \}$,& $\{ 1, 23, 50 \}$,& $\{ 3, 5, 50 \}$,& $\{ 3, 17, 36 \}$,& $\{ 17, 27, 28 \}$,& $\{ 17, 50, 55 \}$,& $\{ 5, 10, 23 \}$, \\$\{ 5, 11, 26 \}$,& $\{ 2, 5, 55 \}$,& $\{ 3, 6, 53 \}$,& $\{ 32, 42, 53 \}$,& $\{ 28, 50, 53 \}$,& $\{ 18, 51, 53 \}$,&
  $\{ 5, 32, 51 \}$, \\$\{ 18, 25, 28 \}$,& $\{ 35, 45, 52 \}$,& $\{ 23, 25, 36 \}$,& $\{ 23, 33, 42 \}$,& $\{ 42, 51, 52 \}$,&
  $\{ 13, 23, 26 \}$, & $\{ 3, 25, 52 \}$,\\$\{ 14, 25, 51 \}$,& $\{ 2, 23, 52 \}$,& $\{ 2, 28, 35 \}$,& $\{ 13, 28, 42 \}$,& $\{ 35, 36, 42 \}$,& $\{ 22, 35, 41 \}$, & $\{ 9, 23, 32 \}$,\\$\{ 2, 21, 51 \}$,&
  $\{ 2, 16, 27 \}$,& $\{ 11, 16, 49 \}$,& $\{ 13, 16, 54 \}$,& $\{ 31, 36, 54 \}$,& $\{ 16, 21, 31 \}$,& $\{ 9, 36, 51 \}$, \\$\{ 6, 9, 27 \}$,& $\{ 9, 14, 28 \}$,& $\{ 7, 9, 22 \}$,& $\{ 7, 10, 25 \}$,& $\{ 15, 22, 25 \}$&
  $\{ 22, 27, 36 \}$,,& $\{ 9, 26, 52 \}$, \\$\{ 14, 21, 27 \}$, & $\{ 2, 7, 49 \}$,& $\{ 13, 27, 50 \}$,& $\{ 11, 13, 36 \}$,& $\{ 11, 50, 52 \}$,& $\{ 36, 41, 50 \}$, & $\{ 7, 32, 50 \}$,\\
  $\{ 15, 28, 41 \}$,&
  $\{ 34, 41, 52 \}$,& $\{ 21, 47, 52 \}$,& $\{ 28, 34, 47 \}$,& $\{ 21, 34, 43 \}$,& $\{ 14, 36, 43 \}$,& $\{ 7, 36, 45 \}$,\\

  $\{ 10, 27, 52 \}$,& $\{ 43, 49, 52 \}$,& $\{ 32, 43, 54 \}$,& $\{ 28, 33, 43 \}$,& $\{ 19, 33, 36 \}$,&
  $\{ 2, 11, 32 \}$,& $\{ 2, 36, 47 \}$, \\$\{ 6, 32, 47 \}$,& $\{ 6, 36, 49 \}$,& $\{ 19, 26, 28 \}$, & $\{ 11, 22, 28 \}$,& $\{ 3, 28, 45 \}$,& $\{ 3, 21, 32 \}$,&
  $\{ 3, 22, 49 \}$, \\$\{ 22, 32, 45 \}$,& $\{ 10, 32, 49 \}$,& $\{ 15, 26, 49 \}$,& $\{ 10, 15, 36 \}$,& $\{ 21, 26, 36 \}$,& $\{ 10, 21, 28 \}$,& $\{ 6, 11, 21 \}$, \\
  $\{ 19, 49, 54 \}$,& $\{ 28, 49, 55 \}$,& $\{ 28, 51, 54 \}$,&
  $\{ 32, 41, 55 \}$.&&&\qed
  \end{longtable}

\begin{example}
We construct a $3$-cyclic $3$-HGDP of type $(6,3^6)$ with $J(6\times6\times3,3,1)$ base blocks on $Z_{6}\times I_{6}\times Z_{3}$ with group set $\{\{i\}\times I_{6}\times Z_{3}:i\in Z_{6}\}$ and hole set $\{Z_{6}\times I_{6}\times Z_{3}:i\in I_{6}\}$. Only $74$ initial base blocks are listed:

\begin{longtable}{c c cc}

$\{(0,0,0),(1,2,0),(2,1,0)\}$,&
$\{(0,0,0),(1,3,0),(2,1,1)\}$,&
$\{(0,0,0),(1,4,0),(5,5,2)\}$,\\
$\{(0,0,0),(1,5,0),(5,4,1)\}$,&
$\{(0,0,0),(2,2,0),(1,1,1)\}$,&
$\{(0,0,0),(2,3,0),(3,2,1)\}$,\\
$\{(0,0,0),(2,4,0),(5,1,1)\}$,&
$\{(0,0,0),(2,5,0),(5,1,0)\}$,&
$\{(0,0,0),(3,1,0),(4,3,2)\}$,\\
$\{(0,0,0),(3,2,0),(5,5,1)\}$,&
$\{(0,0,0),(3,3,0),(4,4,1)\}$,&
$\{(0,0,0),(3,4,0),(2,5,2)\}$,\\
$\{(0,0,0),(3,5,0),(5,2,1)\}$,&
$\{(0,0,0),(4,1,0),(5,3,1)\}$,&
$\{(0,0,0),(4,2,0),(2,3,2)\}$,\\
$\{(0,0,0),(4,3,0),(5,2,0)\}$,&
$\{(0,0,0),(4,4,0),(3,5,1)\}$,&
$\{(0,0,0),(4,5,0),(1,1,2)\}$,\\
$\{(0,0,0),(5,3,0),(4,2,1)\}$,&
$\{(0,0,0),(5,4,0),(1,2,2)\}$,&
$\{(0,0,0),(5,5,0),(1,2,1)\}$,\\
$\{(0,0,0),(1,3,1),(4,4,2)\}$,&
$\{(0,0,0),(1,4,1),(2,2,1)\}$,&
$\{(0,0,0),(1,5,1),(2,3,1)\}$,\\
$\{(0,0,0),(2,4,1),(3,1,1)\}$,&
$\{(0,0,0),(2,5,1),(1,3,2)\}$,&
$\{(0,0,0),(3,3,1),(1,4,2)\}$,\\
$\{(0,0,0),(3,4,1),(4,3,1)\}$,&
$\{(0,0,0),(4,1,1),(1,5,2)\}$,&
$\{(0,0,0),(4,5,1),(2,1,2)\}$,\\
$\{(0,0,0),(2,2,2),(3,4,2)\}$,&
$\{(0,0,0),(2,4,2),(4,1,2)\}$,&
$\{(0,0,0),(3,1,2),(4,2,2)\}$,\\
$\{(0,0,0),(3,2,2),(5,1,2)\}$,&
$\{(0,0,0),(3,3,2),(5,2,2)\}$,&
$\{(0,0,0),(3,5,2),(5,3,2)\}$,\\
$\{(0,0,0),(4,5,2),(5,4,2)\}$,&
$\{(0,1,0),(1,3,0),(2,5,1)\}$,&
$\{(0,1,0),(1,4,0),(4,5,1)\}$,\\
$\{(0,1,0),(1,5,0),(2,2,2)\}$,&
$\{(0,1,0),(2,2,0),(4,3,1)\}$,&
$\{(0,1,0),(2,3,0),(5,5,2)\}$,\\
$\{(0,1,0),(2,4,0),(3,2,1)\}$,&
$\{(0,1,0),(2,5,0),(5,4,0)\}$,&
$\{(0,1,0),(3,2,0),(4,4,1)\}$,\\
$\{(0,1,0),(3,3,0),(4,4,0)\}$,&
$\{(0,1,0),(3,4,0),(4,3,2)\}$,&
$\{(0,1,0),(4,2,0),(5,5,0)\}$,\\
$\{(0,1,0),(4,3,0),(1,5,1)\}$,&
$\{(0,1,0),(4,5,0),(5,2,0)\}$,&
$\{(0,1,0),(5,3,0),(2,4,1)\}$,\\
$\{(0,1,0),(1,2,1),(3,4,2)\}$,&
$\{(0,1,0),(1,4,1),(3,3,1)\}$,&
$\{(0,1,0),(2,2,1),(1,4,2)\}$,\\
$\{(0,1,0),(2,3,1),(5,2,1)\}$,&
$\{(0,1,0),(4,2,1),(1,5,2)\}$,&
$\{(0,1,0),(5,4,1),(3,2,2)\}$,\\
$\{(0,1,0),(5,5,1),(3,3,2)\}$,&
$\{(0,1,0),(2,3,2),(4,5,2)\}$,&
$\{(0,1,0),(2,4,2),(5,3,2)\}$,\\
$\{(0,2,0),(2,3,0),(3,4,2)\}$,&
$\{(0,2,0),(2,4,0),(3,3,1)\}$,&
$\{(0,2,0),(2,5,0),(1,4,2)\}$,\\
$\{(0,2,0),(3,3,0),(4,5,0)\}$,&
$\{(0,2,0),(3,4,0),(1,5,1)\}$,&
$\{(0,2,0),(4,3,0),(1,5,2)\}$,\\
$\{(0,2,0),(4,4,0),(2,5,2)\}$,&
$\{(0,2,0),(5,3,0),(3,5,1)\}$,&
$\{(0,2,0),(1,3,1),(3,5,2)\}$,\\
$\{(0,2,0),(3,4,1),(5,5,1)\}$,&
$\{(0,2,0),(2,3,2),(4,4,2)\}$,&
$\{(0,3,0),(5,5,0),(4,4,1)\}$,\\
$\{(0,3,0),(2,4,1),(5,5,1)\}$,&
$\{(0,3,0),(4,5,1),(2,4,2)\}$.\\
\end{longtable}

All $444$ base blocks are obtained by developing these initial base blocks by $(+1,-,-) ~\bmod$ $ (6,-,-)$ successively.\qed
\end{example}

\section{3-IHGDPs of type $(u,t,1^v)$}
\begin{lemma}
There exists a $3$-IHGDP of type $(u,4,1^{4})$ with $\Theta(u,4,4,1)$ blocks for $u\in\{10,16\}$.
\end{lemma}
\proof
  The design will be constructed on $X=I_{u}\times Z_{4}$ with group set 
 $\mathcal{G}=\{\{i\}\times Z_{4}:i\in I_{u-4}\}\cup\{\{u-4,u-3,u-2,u-1\}\times Z_{4}\}$ and hole set $\mathcal{H}=\{I_{u}\times\{j\}:j\in Z_{4}\}$. We list all $\frac{(u-4)(3u+8)}{6}$ base blocks.

$u=10:$
\begin{longtable}{llllllll}
$\{(0,0),(1,1),(2,2)\}$,&
$\{(0,0),(1,3),(6,1)\}$,&
$\{(0,0),(2,1),(9,2)\}$,&
$\{(0,0),(2,3),(6,2)\}$,\\
$\{(0,0),(3,1),(8,3)\}$,&
$\{(0,0),(3,2),(6,3)\}$,&
$\{(0,0),(3,3),(7,1)\}$,&
$\{(0,0),(4,1),(7,2)\}$,\\
$\{(0,0),(4,2),(8,1)\}$,&
$\{(0,0),(4,3),(9,1)\}$,&
$\{(0,0),(5,1),(9,3)\}$,&
$\{(0,0),(5,2),(7,3)\}$,\\
$\{(0,0),(5,3),(8,2)\}$,&
$\{(1,0),(2,2),(7,3)\}$,&
$\{(1,0),(2,3),(6,1)\}$,&
$\{(1,0),(3,1),(9,3)\}$,\\
$\{(1,0),(3,2),(8,3)\}$,&
$\{(1,0),(3,3),(8,2)\}$,&
$\{(1,0),(4,1),(9,2)\}$,&
$\{(1,0),(4,2),(9,1)\}$,\\
$\{(1,0),(4,3),(7,2)\}$,&
$\{(1,0),(5,1),(6,3)\}$,&
$\{(1,0),(5,2),(7,1)\}$,&
$\{(1,0),(5,3),(8,1)\}$,\\
$\{(2,0),(3,2),(9,3)\}$,&
$\{(2,0),(3,3),(7,2)\}$,&
$\{(2,0),(4,1),(8,2)\}$,&
$\{(2,0),(4,2),(6,1)\}$,\\
$\{(2,0),(4,3),(8,1)\}$,&
$\{(2,0),(5,1),(7,3)\}$,&
$\{(2,0),(5,2),(8,3)\}$,&
$\{(2,0),(5,3),(9,2)\}$,\\
$\{(3,0),(4,1),(5,3)\}$,&
$\{(3,0),(4,2),(6,3)\}$,&
$\{(3,0),(4,3),(7,1)\}$,&
$\{(3,0),(5,1),(6,2)\}$,\\
$\{(3,0),(5,2),(9,3)\}$,&
$\{(4,0),(5,3),(6,2)\}$.&
\end{longtable}
$u = 16:$
\begin{longtable}{lllllllll}
$\{(0,0),(1,1),(2,2)\}$,&
$\{(0,0),(1,2),(10,1)\}$,&
$\{(0,0),(2,1),(15,2)\}$,&
$\{(0,0),(1,3),(9,1)\}$,\\
$\{(3,0),(9,1),(15,3)\}$,
&
$\{(1,0),(11,2),(14,3)\}$,&
$\{(0,0),(3,3),(11,1)\}$,&
$\{(0,0),(4,1),(5,2)\}$,\\
$\{(0,0),(4,2),(12,1)\}$,&
$\{(0,0),(10,2),(14,3)\}$,&
$\{(0,0),(8,2),(13,3)\}$,&
$\{(0,0),(5,3),(7,1)\}$,\\
$\{(0,0),(6,3),(12,2)\}$,&
$\{(0,0),(7,2),(12,3)\}$,&
$\{(0,0),(7,3),(15,1)\}$,&$\{(0,0),(5,1),(6,2)\}$,
\\
$\{(0,0),(8,3),(13,1)\}$,&
$\{(0,0),(9,2),(14,1)\}$,&
$\{(0,0),(9,3),(13,2)\}$,&$\{(0,0),(4,3),(6,1)\}$,
\\
$\{(0,0),(10,3),(14,2)\}$,&
$\{(0,0),(11,2),(15,3)\}$,&
$\{(1,0),(2,2),(10,1)\}$,&
$\{(1,0),(2,3),(9,1)\}$,\\
$\{(1,0),(4,1),(12,3)\}$,&
$\{(1,0),(7,3),(14,2)\}$,&
$\{(1,0),(3,3),(10,2)\}$,&$\{(1,0),(3,1),(4,2)\}$,
\\
$\{(1,0),(4,3),(14,1)\}$,&
$\{(1,0),(5,1),(13,2)\}$,&
$\{(1,0),(5,3),(11,1)\}$,&
$\{(1,0),(5,2),(6,1)\}$,\\
$\{(1,0),(6,2),(13,1)\}$,&
$\{(1,0),(7,1),(13,3)\}$,&
$\{(1,0),(7,2),(12,1)\}$,&$\{(1,0),(3,2),(6,3)\}$,
\\
$\{(1,0),(8,1),(12,2)\}$,&
$\{(1,0),(8,2),(15,3)\}$,&
$\{(1,0),(8,3),(15,1)\}$,&$\{(0,0),(3,2),(8,1)\}$,
\\
$\{(1,0),(11,3),(15,2)\}$,&
$\{(2,0),(3,1),(11,2)\}$,&
$\{(2,0),(4,2),(15,3)\}$,&
$\{(2,0),(3,3),(4,1)\}$,\\
$\{(2,0),(4,3),(14,2)\}$,&$\{(3,0),(7,1),(13,2)\}$,
&
$\{(2,0),(5,3),(12,1)\}$,&
$\{(2,0),(5,2),(7,1)\}$,\\
$\{(2,0),(6,1),(14,3)\}$,&
$\{(3,0),(11,3),(15,1)\}$,&
$\{(2,0),(7,3),(13,2)\}$,&
$\{(2,0),(7,2),(8,3)\}$,\\
$\{(2,0),(9,1),(12,3)\}$,&
$\{(2,0),(9,3),(15,2)\}$,&
$\{(2,0),(10,1),(12,2)\}$,
&
$\{(3,0),(4,3),(5,2)\}$,
\\
$\{(2,0),(11,1),(13,3)\}$,&
$\{(2,0),(11,3),(14,1)\}$,
&
$\{(2,0),(10,2),(13,1)\}$,
&$\{(3,0),(5,1),(7,2)\}$,\\
$\{(3,0),(5,3),(12,2)\}$,&
$\{(3,0),(6,2),(12,3)\}$,&
$\{(3,0),(6,3),(12,1)\}$,&$\{(2,0),(5,1),(6,3)\}$,
\\
$\{(3,0),(7,3),(15,2)\}$,&
$\{(3,0),(8,1),(14,2)\}$,&
$\{(3,0),(8,2),(13,1)\}$,&
$\{(0,0),(2,3),(3,1)\}$,
\\
$\{(3,0),(9,2),(14,3)\}$,&
$\{(3,0),(9,3),(14,1)\}$,&
$\{(3,0),(10,2),(13,3)\}$,&
$\{(2,0),(6,2),(8,1)\}$,

\\
$\{(4,0),(6,1),(13,2)\}$,&
$\{(4,0),(6,3),(15,2)\}$,&
$\{(4,0),(7,2),(15,3)\}$,&
$\{(4,0),(7,1),(9,2)\}$,\\
$\{(4,0),(7,3),(10,1)\}$,&
$\{(4,0),(8,3),(10,2)\}$,&
$\{(4,0),(8,2),(11,1)\}$,&$\{(4,0),(8,1),(9,3)\}$,
\\
$\{(4,0),(11,3),(12,1)\}$,&
$\{(4,0),(10,3),(13,1)\}$,&
$\{(4,0),(11,2),(14,1)\}$,&$\{(4,0),(9,1),(13,3)\}$,
\\
$\{(5,0),(8,1),(14,3)\}$,&
$\{(5,0),(8,2),(12,1)\}$,&
$\{(5,0),(8,3),(14,2)\}$,&
$\{(5,0),(9,1),(15,2)\}$,\\
$\{(5,0),(10,2),(15,1)\}$,&
$\{(5,0),(9,3),(11,1)\}$,&
$\{(5,0),(10,1),(15,3)\}$,&$\{(5,0),(9,2),(13,3)\}$,
\\
$\{(5,0),(10,3),(14,1)\}$,&
$\{(5,0),(11,3),(13,2)\}$,&
$\{(6,0),(7,2),(14,3)\}$,&
$\{(6,0),(7,3),(14,1)\}$,\\
$\{(6,0),(8,1),(10,2)\}$,&
$\{(6,0),(8,2),(15,1)\}$,&
$\{(6,0),(9,1),(10,3)\}$,&
$\{(6,0),(9,2),(11,3)\}$,\\
$\{(6,0),(9,3),(11,2)\}$,&
$\{(6,0),(10,1),(15,2)\}$,&
$\{(6,0),(11,1),(13,2)\}$,&
$\{(7,0),(8,2),(11,3)\}$,\\
$\{(7,0),(8,3),(11,1)\}$,&
$\{(7,0),(9,2),(10,3)\}$,&
$\{(7,0),(10,1),(11,2)\}$,&
$\{(7,0),(9,3),(12,2)\}$,\\
$\{(8,0),(9,1),(12,2)\}$,&
$\{(10,0),(11,3),(12,2)\}$,&
$\{(10,0),(11,2),(12,3)\}$,&$\{(8,0),(9,3),(10,2)\}$
.\\
\end{longtable}

All  blocks are obtained by developing these initial base blocks by $(-,+1) \bmod(-,4)$ successively. \qed

\begin{lemma}
There exists  a $3$-IHGDP of type $(10,4,1^{10})$ with $\Theta(10,4,10,1)$ blocks.
\end{lemma}
\proof
The design will be constructed on  $X=I_{100}$ with group set
$\mathcal{G}=
\{\{i,i+10,\ldots,i+90\}:i\in I_{6}\}\cup\{10j+i:j\in I_{10}, i\in\{6,7,8,9\}\}$ and hole set $\mathcal{H}=\{\{10j+i:i\in I_{10}\}:j\in I_{10}\}$. All $116$ base blocks are listed as follows.
\begin{longtable}{llllllll}
 $\{ 0, 12, 26 \}$,& $\{ 0, 16, 32 \}$,& $\{ 0, 22, 56 \}$,& $\{ 0, 36, 72 \}$,& $\{ 0, 21, 96 \}$,& $\{ 0, 23, 66 \}$,& $\{ 0, 34, 76 \}$,\\
  $\{ 0, 46, 94 \}$,&
  $\{ 0, 65, 86 \}$,& $\{ 1, 15, 26 \}$,& $\{ 1, 16, 35 \}$,& $\{ 1, 36, 85 \}$,& $\{ 1, 23, 86 \}$,& $\{ 1, 13, 46 \}$,\\ $\{ 1, 56, 64 \}$,& $\{ 1, 24, 96 \}$,&
   $\{ 1, 66, 95 \}$,& $\{ 2, 14, 96 \}$,& $\{ 2, 26, 33 \}$,& $\{ 2, 43, 56 \}$,& $\{ 2, 46, 84 \}$,\\ $\{ 3, 15, 76 \}$,& $\{ 2, 35, 76 \}$,& $\{ 3, 14, 26 \}$,&
   $\{ 3, 24, 56 \}$,& $\{ 3, 55, 86 \}$,& $\{ 4, 26, 35 \}$,& $\{ 0, 13, 47 \}$,\\ $\{ 0, 14, 87 \}$,& $\{ 0, 17, 54 \}$,& $\{ 0, 27, 71 \}$,& $\{ 0, 37, 83 \}$,&
    $\{ 0, 15, 97 \}$,& $\{ 0, 45, 57 \}$,& $\{ 0, 33, 77 \}$,\\
  $\{ 0, 42, 67 \}$,& $\{ 1, 14, 37 \}$,& $\{ 1, 27, 74 \}$,& $\{ 1, 17, 22 \}$,& $\{ 1, 25, 87 \}$,&
   $\{ 1, 47, 75 \}$,& $\{ 1, 67, 92 \}$,\\ $\{ 1, 42, 97 \}$,& $\{ 1, 77, 83 \}$,& $\{ 2, 17, 24 \}$,& $\{ 2, 15, 47 \}$,& $\{ 2, 63, 87 \}$,& $\{ 2, 45, 67 \}$,&
    $\{ 2, 37, 73 \}$,\\ $\{ 3, 17, 34 \}$,& $\{ 3, 44, 87 \}$,& $\{ 3, 35, 77 \}$,& $\{ 4, 17, 25 \}$,& $\{ 4, 37, 85 \}$,& $\{ 0, 18, 61 \}$,& $\{ 0, 38, 82 \}$,\\
  $\{ 0, 28, 63 \}$,& $\{ 0, 24, 88 \}$,& $\{ 0, 25, 48 \}$,& $\{ 0, 31, 58 \}$,& $\{ 0, 35, 98 \}$,& $\{ 0, 19, 74 \}$,& $\{ 0, 29, 75 \}$,\\
  $\{ 0, 39, 64 \}$,&
   $\{ 0, 41, 95 \}$,& $\{ 0, 43, 59 \}$,& $\{ 0, 51, 89 \}$,& $\{ 0, 44, 53 \}$,& $\{ 0, 52, 99 \}$,& $\{ 0, 49, 73 \}$,\\ $\{ 0, 55, 93 \}$,& $\{ 0, 62, 78 \}$,&
  $\{ 0, 68, 81 \}$,& $\{ 0, 69, 85 \}$,& $\{ 0, 79, 91 \}$,& $\{ 0, 84, 92 \}$,& $\{ 1, 19, 93 \}$,\\
  $\{ 1, 29, 73 \}$,& $\{ 1, 32, 59 \}$,& $\{ 1, 49, 72 \}$,&
  $\{ 1, 44, 99 \}$,& $\{ 1, 33, 68 \}$,& $\{ 1, 12, 69 \}$,& $\{ 1, 18, 63 \}$,\\ $\{ 1, 34, 48 \}$,& $\{ 1, 38, 53 \}$,& $\{ 1, 45, 84 \}$,& $\{ 1, 52, 98 \}$,&
  $\{ 1, 54, 78 \}$,& $\{ 1, 65, 82 \}$,& $\{ 1, 43, 62 \}$,\\
  $\{ 1, 79, 94 \}$,& $\{ 2, 55, 89 \}$,& $\{ 2, 38, 64 \}$,& $\{ 2, 25, 99 \}$,& $\{ 3, 74, 99 \}$,&
  $\{ 3, 75, 89 \}$,& $\{ 3, 39, 45 \}$,\\ $\{ 2, 65, 98 \}$,& $\{ 2, 23, 68 \}$,& $\{ 3, 25, 98 \}$,& $\{ 3, 18, 64 \}$,& $\{ 2, 75, 88 \}$,& $\{ 4, 45, 98 \}$,&
  $\{ 3, 78, 95 \}$,\\ $\{ 2, 28, 44 \}$,& $\{ 4, 38, 95 \}$,& $\{ 3, 49, 85 \}$,& $\{ 2, 74, 93 \}$,& $\{ 2, 19, 54 \}$,& $\{ 2, 39, 95 \}$,& $\{ 3, 54, 69 \}$,\\
  $\{ 2, 53, 78 \}$,& $\{ 4, 48, 55 \}$,& $\{ 2, 34, 69 \}$,& $\{ 4, 75, 99 \}$.
\end{longtable}

 Let $G$ be the group generated by $\alpha=\prod_{i=0}^{9}(i\ 10+i \cdots 90+i)$. All blocks are obtained by
  developing the above base blocks under the action of $G.$\qed

\begin{lemma}
There exists a $3$-IHGDP of type $(10,4,1^{v})$ with $\Theta(10,4,v,1)$ blocks for $v\in\{8,14\}$.
\end{lemma}
\proof
 We construct the design on $X=I_{10v}$ with group set $\mathcal{G}=\{\{10i+j:i\in I_{v}\}:j\in \{4,5,\ldots,9\}\}\cup\{10i+j:i\in I_{v},~j\in I_{4}\}$  and hole set $\mathcal{H}=\{\{10i+j:j\in I_{10}\}:i\in I_{v}\}$. Next we list  $\frac{v(13v-14)}{3}$ base blocks.

$v = 8 :$
\begin{longtable}{llllllll}
$\{ 0, 14, 25 \}$,& $\{ 0, 17, 58 \}$,& $\{ 0, 27, 44 \}$,& $\{ 0, 34, 77 \}$,& $\{ 0, 37, 64 \}$,& $\{ 0, 47, 54 \}$,& $\{ 0, 67, 75 \}$,\\
  $\{ 1, 17, 26 \}$, &$\{ 1, 14, 38 \}$,& $\{ 1, 27, 55 \}$, &$\{ 1, 34, 75 \}$, &$\{ 1, 44, 65 \}$, &$\{ 1, 47, 58 \}$,& $\{ 1, 67, 79 \}$,\\
  $\{ 2, 14, 39 \}$, &$\{ 2, 24, 75 \}$, &$\{ 2, 27, 58 \}$,& $\{ 2, 17, 46 \}$, &$\{ 2, 34, 56 \}$,& $\{ 2, 47, 68 \}$, &$\{ 2, 64, 79 \}$,\\
  $\{ 3, 27, 76 \}$,& $\{ 4, 17, 56 \}$, &$\{ 3, 44, 79 \}$,& $\{ 3, 14, 59 \}$, &$\{ 3, 24, 36 \}$,& $\{ 3, 37, 56 \}$, &$\{ 3, 47, 65 \}$,\\
  $\{ 3, 17, 68 \}$,& $\{ 4, 18, 37 \}$.
\end{longtable}

$v = 14:$
\begin{longtable}{lllllll}
   $\{ 0, 14, 25 \}$,& $\{ 0, 17, 38 \}$,& $\{ 0, 27, 44 \}$,& $\{ 0, 34, 65 \}$,& $\{ 0, 47, 58 \}$,& $\{ 0, 54, 86 \}$, \\
  $\{ 4, 18, 66 \}$,& $\{ 4, 16, 28 \}$,& $\{ 2, 14, 65 \}$,& $\{ 3, 14, 35 \}$,& $\{ 3, 27, 105 \}$, &$\{ 0, 67, 98 \}$,\\ $\{ 4, 88, 59 \}$,& $\{ 1, 14, 55 \}$,&
  $\{ 1, 17, 68 \}$,& $\{ 1, 24, 95 \}$,& $\{ 1, 27, 108 \}$,& $\{ 1, 34, 76 \}$, \\
  $\{ 1, 37, 89 \}$,& $\{ 1, 44, 135 \}$,&
  $\{ 3, 17, 75 \}$,& $\{ 2, 17, 78 \}$,& $\{ 1, 97, 125 \}$,& $\{ 1, 47, 59 \}$, \\

  $\{ 2, 54, 99 \}$,& $\{ 3, 57, 128 \}$,& $\{ 0, 87,124 \}$,& $\{ 1, 77, 115 \}$,&
  $\{ 2, 27, 135 \}$,& $\{ 2, 24, 88 \}$, \\
   $\{ 2, 34, 47 \}$,& $\{ 2, 44, 117 \}$,& $\{ 2, 57, 125 \}$,& $\{ 2, 67, 106 \}$,&
  $\{ 2, 74, 127 \}$, & $\{ 2, 37, 94 \}$,\\
  $\{ 1, 84, 117 \}$,& $\{ 4, 27, 129 \}$,& $\{ 2, 84, 138 \}$,& $\{ 0, 94, 138 \}$,& $\{ 3, 54, 115 \}$,& $\{ 3, 24, 46 \}$,\\
  $\{ 0, 77, 104 \}$,& $\{ 0, 74, 117 \}$,& $\{ 1, 64, 127 \}$,& $\{ 3, 84, 109 \}$,& $\{ 3, 77, 126 \}$, & $\{ 3, 67, 96 \}$, \\
   $\{ 3, 37, 138 \}$,&
  $\{ 2, 107, 116 \}$,& $\{ 3, 87, 134 \}$,& $\{ 3, 97,119 \}$,& $\{ 0, 107, 114 \}$, & $\{ 3, 47, 66 \}$,\\ $\{ 0, 127, 135 \}$,& $\{ 1, 104, 138 \}$.
\end{longtable}

  Let $G$ be the group generated by $\alpha$ and $\beta$ where $$\alpha=\prod_{i=0}^{9}(i\ 10+i\cdots 10(v-1)+i),~ \beta=\prod_{j=0}^{v-1}(10j+4~10j+5~10j+6)(10j+7~10j+8~10j+9).$$ All blocks are obtained by
  developing the above base blocks under the action of $G.$\qed
\begin{lemma}
There exists a $3$-IHGDP of type $(16,4,1^{v})$  with $\Theta(16,4,v,1)$ blocks for $v\in\{8,14\}$.
\end{lemma}
\proof
First let
$v=8 $.  The design will be constructed on $X=I_{128}$ with group set $\mathcal{G}=\{\{16i+j:i\in I_{8}\}:j\in\{4,5,\ldots,15\}\}\cup\{16i+j:i\in I_{8},j\in I_{4}\}$ and hole set $\mathcal{H}=\{\{16i+j:j\in I_{16}\}:i\in I_{8}\}$. We list $44$ base blocks as follows.
\begin{longtable}{llllll}
$\{ 0, 20, 37 \}$,& $\{ 0, 26, 59 \}$,& $\{ 0, 42, 68 \}$,& $\{ 0, 52, 101 \}$,& $\{ 0, 74, 91 \}$,& $\{ 0, 84, 122 \}$, \\$\{ 4, 61, 69 \}$,&
  $\{ 1, 20, 59 \}$,& $\{ 1, 36, 76 \}$,& $\{ 1, 26, 117 \}$,& $\{ 1, 42, 91 \}$,& $\{ 1, 52, 106 \}$, \\ $\{ 1, 68, 86 \}$,& $\{ 4, 59, 94 \}$,&
  $\{ 2, 20, 94 \}$,& $\{ 2, 26, 107 \}$,& $\{ 2, 36, 55 \}$,& $\{ 2, 42, 124 \}$, \\ $\{ 2, 58, 69 \}$,&$\{ 4, 25, 55 \}$, &
  $\{ 3, 20, 93 \}$,& $\{ 2, 84, 117 \}$,& $\{ 3, 36, 88 \}$,& $\{ 3, 52, 105 \}$, \\ $\{ 3, 26, 123 \}$,& $\{ 3, 42, 70 \}$,&
  $\{ 4, 27, 46 \}$,& $\{ 3, 58, 108 \}$,& $\{ 4, 24, 76 \}$,& $\{ 3, 74, 116 \}$, \\$\{ 2, 74, 101 \}$, & $\{ 4, 29, 54 \}$,& $\{ 4, 30, 70 \}$,& $\{ 10, 28, 79 \}$,& $\{ 4, 38, 91 \}$,&
  $\{ 4, 60, 126 \}$, \\ $\{ 4, 39, 109 \}$,& $\{ 4, 75, 95 \}$,& $\{ 4, 28, 41 \}$,& $\{ 4, 31, 123 \}$,& $\{ 4, 45, 79 \}$,& $\{ 4, 93, 108 \}$,\\$\{ 0, 106, 116 \}$,& $\{ 1, 100, 122 \}$.
\end{longtable}

 Let $G$ be the group generated by $\alpha$ and $\beta$ where $$\alpha=\prod_{i=0}^{15}(i\ 16+i \ 32+i \cdots 112+i),$$
  $$\beta=\prod_{i=0}^{7}(16i+4\ 16i+5\ \cdots\ 16i+9)(16i+10\ 16i+11\ \cdots \ 16i+15).$$
 All  blocks are obtained by developing the above base blocks under the action of $G.$

\indent Then let
$v=14.$ The design will be constructed on $X=I_{324}$ with group set $\mathcal{G}=\{\{16i+j:i\in I_{14}\}:j\in I_{12}\}\cup\{16i+j:i\in I_{14},j\in \{12,13,14,15\}\}$ and hole set $\mathcal{H}=\{\{16i+j:j\in I_{16}\}:i\in I_{14}\}$. We list $82$ base blocks as follows.

\begin{longtable}{llllll}
$\{ 0, 18, 47 \}$,& $\{ 0, 17, 78 \}$,& $\{ 0, 19, 94 \}$,& $\{ 0, 73, 82 \}$,& $\{ 0, 37, 92 \}$,& $\{ 0, 22, 191 \}$,\\  $\{ 0, 35, 79 \}$,&
  $\{ 0, 24, 46 \}$,& $\{ 0, 69, 98 \}$,  & $\{ 0, 25, 172 \}$,& $\{ 1, 41, 95 \}$,& $\{ 0, 30, 198 \}$,\\$\{ 1, 46, 55 \}$, & $\{ 0, 34, 95 \}$,&$\{ 0, 59, 97 \}$,
 & $\{ 0, 39, 143 \}$,& $\{ 0, 77, 88 \}$,  & $\{ 0, 40, 220 \}$,\\  $\{ 0, 23, 206 \}$, & $\{ 1, 30, 37 \}$,&  $\{ 1, 78, 83 \}$, &
  $\{ 0, 44, 193 \}$,& $\{ 0, 49, 142 \}$,& $\{ 0, 41, 223 \}$,\\ $\{ 0, 50, 110 \}$,& $\{ 0, 52, 158 \}$,& $\{ 1, 19, 47 \}$,& $\{ 0, 76, 101 \}$,&
  $\{ 1, 21, 172 \}$,& $\{ 1, 23, 142 \}$,\\ $\{ 1, 25, 206 \}$,&$\{ 0, 42, 207 \}$,& $\{ 1, 31, 87 \}$,& $\{ 1, 27, 140 \}$,&$\{ 0, 27, 127 \}$,& $\{ 1, 35, 204 \}$,\\$\{ 0, 33, 140 \}$, & $\{ 1, 39, 156 \}$,& $\{ 0, 43, 126 \}$,
  & $\{ 0, 29, 163 \}$,& $\{ 0, 45, 105 \}$,& $\{ 0, 51, 107 \}$,\\
  $\{ 0, 53, 116 \}$,& $\{ 0, 54, 167 \}$,& $\{ 0, 55, 171 \}$,& $\{ 0, 56, 155 \}$,& $\{ 0, 58, 114 \}$,& $\{ 0, 57, 103 \}$,\\ $\{ 0, 61, 165 \}$,&
  $\{ 0, 65, 187 \}$,&   $\{ 0, 36, 188 \}$,      & $\{ 0, 66, 185 \}$,& $\{ 0, 70, 179 \}$,& $\{ 0, 81, 136 \}$,\\ $\{ 0, 67, 104 \}$,& $\{ 0, 71, 106 \}$,&$\{ 0, 93, 135 \}$,
  &$\{ 0, 20, 124 \}$,& $\{ 0, 74, 146 \}$,& $\{ 0, 86, 213 \}$,\\ $\{ 0, 75, 133 \}$,& $\{ 0, 83, 102 \}$,& $\{ 0, 26, 222 \}$,& $\{ 0, 115, 203 \}$,& $\{ 0, 21, 156 \}$,&
  $\{ 0, 85, 169 \}$,\\
   $\{ 0, 89, 129 \}$,& $\{ 0, 72, 197 \}$,& $\{ 0, 109, 177 \}$,& $\{ 0, 131, 205 \}$,&
  $\{ 0, 117, 157 \}$,& $\{ 0, 91, 183 \}$,\\ $\{ 0, 87, 189 \}$,& $\{ 1, 61, 147 \}$,& $\{ 0, 123, 209 \}$,& $\{ 0, 145, 211 \}$,&
  $\{ 0, 151, 219 \}$,& $\{ 1, 43, 131 \}$,\\
  $\{ 1, 73, 135 \}$,& $\{ 0, 125, 147 \}$,& $\{ 0, 141, 149 \}$,& $\{ 0, 153, 173 \}$.
\end{longtable}
\noindent
Let $G$ be the group generated by $\alpha$ and $\beta$ where $$\alpha=\prod_{i=0}^{15}(i\ 16+i \ 32+i \cdots 208+i),$$  $$\beta=\prod_{i=0}^{13}(16i\ 16i+2\  \cdots\ 16i+10)(16i+1\ 16i+3\ \cdots \ 16i+11)(16i+12\ 16i+15\ 16i+ 14).$$ All blocks are obtained by developing the above base blocks under the action of $G.$ \qed
\section{$K$-SCGDDs of type $w^v$}
\begin{lemma}
 For $v\in\{14,20\}$, there exists a $\{3,4,6\}$-SCGDD of type $3^{v}$ containing two base blocks of size $4$, in which all blocks of size $4$ and size $6$ form a parallel class.
\end{lemma}
\proof
We construct the desired design on $X= I_{3v}$ with group set  $\{\{i,v+i,2v+i\}:i\in I_{3v}\}$. All $\frac{3v^2-8v+16}{6}$ base blocks of size $3$, two base blocks of size $4$ and $\frac{v-8}{6}$ base blocks of size $6$ are listed as follows.

 $v=14:$

\begin{longtable}{llllllll}
$\{ 0, 4, 8 \}$,& $\{ 0, 7, 39 \}$,& $\{ 0, 5, 37 \}$,& $\{ 1, 4, 23 \}$,& $\{ 0, 9, 36 \}$,& $\{ 1, 8, 37 \}$, &$\{ 0, 18, 23 \}$,\\
$\{ 1, 6, 9 \}$,&
  $\{ 2, 7, 37 \}$,& $\{ 2, 9, 35 \}$, &$\{ 2, 10, 23 \}$, &$\{ 3, 9, 25 \}$, &$\{ 1, 7, 34 \}$,& $\{ 4, 26, 37 \}$, \\
  $\{ 6, 8, 40 \}$,& $\{ 5, 9, 20 \}$,&
  $\{ 1, 5, 18 \}$,& $\{ 5, 23, 38 \}$, &$\{ 7, 9, 27 \}$, &$\{ 2, 4, 22 \}$,& $\{ 1, 20, 36 \}$,\\
   $\{ 2, 6, 38 \}$, &$\{ 2, 8, 17 \}$,& $\{ 3, 8, 24 \}$,&
  $\{ 2, 27, 36 \}$,& $\{ 3, 6, 10 \}$,& $\{ 5, 22, 34 \}$, &$\{ 5, 11, 36 \}$,\\
  $\{ 2, 5, 41 \}$, &$\{ 7, 12, 36 \}$,& $\{ 7, 10, 22 \}$,&
  $\{ 8, 27, 39 \}$,& $\{ 0, 6, 30 \}$,& $\{ 2, 11, 34 \}$, &$\{ 0, 20, 31 \}$,\\ $\{ 0, 12, 16 \}$,& $\{ 0, 22, 32 \}$, &$\{ 0, 10, 35 \}$,&
  $\{ 0, 11, 26 \}$, &$\{ 0, 15, 27 \}$,& $\{ 0, 13, 38 \}$,& $\{ 0, 17, 19 \}$, \\$\{ 0, 24, 29 \}$,& $\{ 0, 21, 40 \}$, &$\{ 0, 25, 33 \}$,&
  $\{ 0, 34, 41 \}$, &$\{ 3, 19, 22 \}$, &$\{ 3, 12, 37 \}$,& $\{ 3, 20, 41 \}$,\\
  $\{ 1, 21, 22 \}$, &$\{ 4, 34, 40 \}$,& $\{ 6, 26, 35 \}$,&
  $\{ 4, 13, 20 \}$, &$\{ 6, 24, 39 \}$, &$\{ 1, 31, 35 \}$,& $\{ 3, 21, 27 \}$, \\$\{ 2, 21, 33 \}$,& $\{ 4, 27, 35 \}$,& $\{ 5, 10, 21 \}$,&
  $\{ 3, 35, 39 \}$,& $\{ 4, 21, 39 \}$,& $\{ 3, 18, 26 \}$, &$\{ 4, 11, 24 \}$,\\ $\{ 3, 32, 38 \}$,& $\{ 2, 31, 32 \}$, &$\{ 1, 17, 25 \}$,&
  $\{ 3, 33, 40 \}$,& $\{ 1, 24, 32 \}$,& $\{ 1, 10, 40 \}$,& $\{ 5, 24, 40 \}$,\\ $\{ 2, 13, 24 \}$,& $\{ 2, 18, 39 \}$, &$\{ 6, 11, 23 \}$,&
  $\{ 4, 19, 41 \}$,& $\{ 1, 16, 33 \}$,& $\{ 1, 11, 27 \}$,& $\{ 5, 13, 26 \}$,\\ $\{ 3, 13, 23\}$,
  &$\{ 2, 12, 25 \}$, &$\{ 1, 26, 41 \}$,&
  $\{ 1, 12, 30 \}$,& $\{ 1, 19, 39 \};$ \smallskip\\$\{0,1,2,3\}$,& $\{4,5,6,7\};$
  \smallskip\\
  \multicolumn{3}{l}{$\{8,9,10,11,12,13\}.$}
  \end{longtable}
$v=20:$
\begin{longtable}{llllllll}
 $\{ 0, 8, 14 \}$,& $\{ 4, 9, 28 \}$,& $\{ 2, 6, 8 \}$,& $\{ 7, 9, 48 \}$,& $\{ 0, 4, 52 \}$,& $\{ 0, 5, 28 \}$,& $\{ 2, 12, 54 \}$,\\
  $\{ 3, 5, 10 \}$,& $\{ 2, 4, 34 \}$,& $\{ 5, 8, 33 \}$,& $\{ 0, 9, 56 \}$,& $\{ 3, 6, 50 \}$,& $\{ 0, 7, 12 \}$,& $\{ 5, 34, 46 \}$,\\
  $\{ 0, 6, 41 \}$,& $\{ 2, 5, 55 \}$,& $\{ 4, 16, 54 \}$,& $\{ 6, 17, 34 \}$,& $\{ 3, 7, 51 \}$,& $\{ 2, 7, 15 \}$,& $\{ 7, 29, 34 \}$,\\
  $\{ 1, 5, 49 \}$,& $\{ 2, 9,53 \}$,& $\{ 2, 19, 57 \}$,& $\{ 5, 17, 47 \}$,& $\{ 1, 8, 42 \}$,& $\{ 1, 9, 39 \}$,&$\{ 9, 51, 54 \}$,\\
  $\{ 1, 4, 48 \}$,& $\{ 3, 9, 32 \}$,& $\{ 0, 29, 54 \}$,& $\{ 2, 48, 56 \}$,& $\{ 3, 8, 31 \}$,& $\{ 5, 9, 30 \}$,& $\{ 3, 19, 28 \}$, \\ $\{ 4, 8, 46 \}$,&
  $\{ 0, 34, 48 \}$,& $\{ 3, 14, 55 \}$,& $\{ 1, 10, 34 \}$,& $\{ 1, 7, 32 \}$,& $\{ 7, 8, 38 \}$, & $\{ 6, 48, 55 \}$,\\ $\{ 7, 28, 59 \}$,&
  $\{ 8, 17, 30 \}$,& $\{ 8, 18, 56 \}$,& $\{ 8, 19, 50 \}$,& $\{ 8, 32, 51 \}$,& $\{ 1, 6, 16 \}$,& $\{ 8, 53, 57 \}$, \\ $\{ 6, 14, 27 \}$,&
  $\{ 0, 11, 26 \}$,& $\{ 0, 19, 46 \}$,& $\{ 2, 26, 30 \}$,& $\{ 3, 25, 34 \}$,& $\{ 6, 9, 52 \}$,& $\{ 1, 11, 46 \}$,\\ $\{ 2, 36, 46 \}$,&
  $\{ 3, 13, 26 \}$,& $\{ 3, 46, 53 \}$,& $\{ 4, 26, 49 \}$,& $\{ 5, 26, 58 \}$,& $\{ 8, 36, 55 \}$,& $\{ 6, 11, 18 \}$, \\ $\{ 6, 15, 32 \}$,&
  $\{ 6, 12, 37 \}$,& $\{ 6, 33, 59 \}$,& $\{ 6, 19, 58 \}$,& $\{ 6, 30, 36 \}$,& $\{ 6, 35, 47 \}$,& $\{ 6, 49, 57 \}$, \\ $\{ 0, 21, 53 \}$,&
  $\{ 0, 13, 36 \}$,& $\{ 1, 53, 56 \}$,& $\{ 3, 15, 54 \}$,& $\{ 1, 13, 59 \}$,& $\{ 0, 33, 47 \}$,& $\{ 2, 13, 32 \}$, \\ $\{ 2, 10, 33 \}$,&
  $\{ 0, 39, 42 \}$,& $\{ 4, 15, 53 \}$,& $\{ 3, 33, 38 \}$,& $\{ 4, 18, 33 \}$,& $\{ 5, 18, 53 \}$,& $\{ 7, 13, 49 \}$, \\ $\{ 7, 33, 35 \}$,&
  $\{ 3, 11, 57 \}$,& $\{ 5, 36, 51 \}$,& $\{ 2, 43, 44 \}$,& $\{ 7, 14, 58 \}$,& $\{ 1, 14, 31 \}$,& $\{ 1, 12, 58 \}$, \\ $\{ 2, 14, 52 \}$,&
  $\{ 1, 17, 35 \}$,& $\{ 1, 15, 29 \}$,& $\{ 1, 19, 38 \}$,& $\{ 5,  32,48 \}$,& $\{ 1, 18, 23 \}$,& $\{ 1, 24, 30 \}$, \\ $\{ 1, 37, 52 \}$,&
  $\{ 1, 22, 45 \}$,& $\{ 1, 25, 47 \}$,& $\{ 1, 27, 43 \}$,& $\{ 1, 36, 57 \}$,& $\{ 1, 44, 51 \}$,& $\{ 0, 18, 49 \}$, \\ $\{ 2, 18, 29 \}$,&
  $\{ 1, 28, 55 \}$,& $\{ 3, 29, 36 \}$,& $\{ 2, 17, 49 \}$,& $\{ 4, 29, 35 \}$,& $\{ 3, 17, 48 \}$, & $\{ 3, 39, 49 \}$,\\ $\{ 5, 29, 56 \}$,&
  $\{ 9, 19, 31 \}$,& $\{ 9, 18, 35 \}$,& $\{ 9, 50, 57 \}$,& $\{ 0, 10, 38 \}$,& $\{ 0, 35, 58 \}$,& $\{ 3, 12, 18 \}$, \\ $\{ 2, 38, 50 \}$,&
  $\{ 2, 11, 58 \}$,& $\{ 4, 38, 56 \}$,& $\{ 4, 25, 58 \}$,& $\{ 4, 27, 37 \}$,& $\{ 0, 17, 27 \}$,& $\{ 7, 18, 52 \}$, \\ $\{ 0, 15, 30 \}$,&
  $\{ 5, 13, 54 \}$,& $\{ 0, 22, 55 \}$,& $\{ 4, 32, 55 \}$,& $\{ 3, 30, 35 \}$,& $\{ 8, 34, 58 \}$,& $\{ 5, 15, 37 \}$, \\ $\{ 5, 31, 35 \}$,&
  $\{ 7, 39, 56 \}$,& $\{ 2, 24, 51 \}$,& $\{ 4, 19, 51 \}$,& $\{ 4, 30, 45 \}$,& $\{ 0, 32, 50 \}$,& $\{ 4, 12, 50 \}$, \\ $\{ 4, 47, 59 \}$,&
  $\{ 7, 10, 36 \}$,& $\{ 7, 16, 30 \}$,& $\{ 7, 31, 50 \}$& $\{ 1, 50, 54 \}$,& $\{ 4, 13, 14 \}$,& $\{ 0, 16, 51 \}$,\\ $\{ 5, 12, 39 \}$,&
  $\{ 11, 32, 37 \}$,& $\{ 3, 16, 52 \}$,& $\{ 5, 52, 59 \}$,& $\{ 3, 45, 56 \}$,& $\{ 0, 25, 31 \}$,& $\{ 2, 16, 39 \}$, \\ $\{ 4, 17, 36 \}$,&
  $\{ 10, 35, 59 \}$,& $\{ 11, 35, 53 \}$,& $\{ 13, 34, 57 \}$,& $\{ 0, 45, 59 \}$,& $\{ 2, 37, 45 \}$,& $\{ 0, 23, 57 \}$, \\ $\{ 2, 27, 31 \}$,&
  $\{ 10, 31, 53 \}$,& $\{ 11, 16, 55 \}$,& $\{ 12, 16, 33 \}$,& $\{ 13, 19, 37 \}$,& $\{ 3, 24, 59 \}$,& $\{ 0, 24, 43 \}$, \\ $\{ 2, 23, 47 \}$,&
  $\{ 11, 19, 34 \}$,& $\{ 12, 14, 56 \}$,& $\{ 12, 15, 59 \}$,& $\{ 11, 17, 38 \}$,& $\{ 10, 18, 37 \}$,& $\{ 0, 37, 44 \}$, \\
  $\{ 10, 19, 54 \};$\smallskip\\
 $\{0,1,2,3\}$,& $\{4,5,6,7\};$\smallskip\\
\multicolumn{6}{l}{ $\{8,9,10,11,12,13\}$, $\{14,15,16,17,18,19\}$.}
 \end{longtable}

Let $G$ be the group generated by $\alpha=\prod_{i=0}^{v-1}(i\ v+i\ 2v+i)$. All blocks are obtained by developing these base blocks under the action of $G.$ Obviously, it is isomorphic to a $\{3,4,6\}$-SCGDD of type $3^{v}$.
\qed
\end{document}